\documentclass[a4paper,11pt]{article}
\usepackage{graphicx}
\usepackage{amssymb}
\usepackage{undertilde}
\usepackage{color}
\usepackage{url}
\catcode`\@=11 \@addtoreset{equation}{section}

\catcode`\@=12

\newtheorem{Theorem}{Theorem}[section]
\newtheorem{Proposition}{Proposition}[section]
\newtheorem{Lemma}{Lemma}[section]
\newtheorem{Corollary}{Corollary}[section]
\newtheorem{Remark}{Remark}[section]
\newtheorem{Definition}{Definition}[section]

\newcommand{\bTheorem}[1]{
\begin{Theorem} \label{T#1} }
\newcommand{\eT}{\end{Theorem}}

\newcommand{\bProposition}[1]{
\begin{Proposition} \label{P#1}}
\newcommand{\eP}{\end{Proposition}}

\newcommand{\bLemma}[1]{
\begin{Lemma} \label{L#1} }
\newcommand{\eL}{\end{Lemma}}

\newcommand{\bCorollary}[1]{
\begin{Corollary} \label{C#1} }
\newcommand{\eC}{\end{Corollary}}

\newcommand{\bRemark}[1]{
\begin{Remark} \label{R#1} }
\newcommand{\eR}{\end{Remark}}

\newcommand{\bDefinition}[1]{
\begin{Definition} \label{D#1} }
\newcommand{\eD}{\end{Definition}}

\newcommand{\bFormula}[1]{
\begin{equation} \label{#1}}
\newcommand{\eF}{\end{equation}}

\newcommand{\Divh}{{\rm div}_h}
\newcommand{\Gradh}{\nabla_h}

\newcommand{\Ov}[1]{\overline{#1}}
\newcommand{\av}[1]{ \left\{ #1 \right\}}

\newcommand{\aleq}{\stackrel{<}{\sim}}

\newcommand{\vr}{\varrho}

\newcommand{\vu}{\vc{u}}
\newcommand{\vc}[1]{{\bf #1}}

\newcommand{\Div}{{\rm div}_x}
\newcommand{\Grad}{\nabla_x}

\newcommand{\tn}[1]{\mbox {\F #1}}
\newcommand{\dx}{{\rm d} {x}}
\newcommand{\dt}{{\rm d} t }

\newcommand{\ju}[1]{[[ #1 ]]}

\newcommand{\intO}[1]{\int_{\Omega} #1 \ \dx}

\newcommand{\vv}{\vc{v}}

\newcommand{\ep}{\varepsilon}

\font\F=msbm10 scaled 1000
\newcommand{\R}{\mbox{\F R}}

\definecolor{Cgrey}{rgb}{0.85,0.85,0.85}
\definecolor{Cblue}{rgb}{0.50,0.85,0.85}
\definecolor{Cred}{rgb}{1,0,0}
\definecolor{fancy}{rgb}{0.10,0.85,0.10}

\newcommand\Cbox[2]{%
    \newbox\contentbox%
    \newbox\bkgdbox%
    \setbox\contentbox\hbox to \hsize{%
        \vtop{
            \kern\columnsep
            \hbox to \hsize{%
                \kern\columnsep%
                \advance\hsize by -2\columnsep%
                \setlength{\textwidth}{\hsize}%
                \vbox{
                    \parskip=\baselineskip
                    \parindent=0bp
                    #2
                }%
                \kern\columnsep%
            }%
            \kern\columnsep%
        }%
    }%
    \setbox\bkgdbox\vbox{
        \color{#1}
        \hrule width  \wd\contentbox %
               height \ht\contentbox %
               depth  \dp\contentbox
        \color{black}
    }%
    \wd\bkgdbox=0bp%
    \vbox{\hbox to \hsize{\box\bkgdbox\box\contentbox}}%
    \vskip\baselineskip%
}

\date{}

\begin{document}


\title{Consistency, convergence and error estimates for a mixed finite element--finite volume scheme 
to compressible Navier-Stokes equations with general inflow/outflow boundary data}

\author{Young-Sam Kwon
\thanks{ The work of the first author was partially supported by NRF-
2017R1D1A1B03030249 and NRF-2019H1D3A2A01101128.}   \and
Antonin Novotn{y}\thanks{The work of the second author was supported by Brain Pool program funded by the Ministry of Science and ICT through the National Research Foundation of Korea, NRF-2019H1D3A2A01101128.}}

\maketitle

\bigskip

\centerline{Department of Mathematics, Dong-A University} 
\centerline{Busan 604-714, Republic of Korea, ykwon@dau.ac.kr}

\medskip
\centerline{University of Toulon, IMATH, EA 2134,  BP 20139} 
\centerline{839 57 La Garde, France, novotny@univ-tln.fr}

\maketitle

\bigskip

\begin{abstract}
We study  convergence  of a mixed finite element-finite volume scheme for the compressible Navier-Stokes equations in the isentropic regime
under the full range $1<\gamma<\infty$ of the adiabatic coefficient $\gamma$ for the problem with general non zero inflow-outflow boundary conditions.
We propose a modification of Karper's scheme [Numer. Math. 125:441-510, 2013] in order to accommodate the non zero boundary data, prove existence of its 
solutions, establish the stability and uniform estimates, derive a convenient consistency formulation of the balance laws and use it
to show the weak convergence of the numerical solutions to a dissipative solution with the Reynolds defect introduced in Abbatiello et al. [Preprint Arxiv 1912.12896].
If the target system admits a strong solution then the convergence is strong towards the strong solution. Moreover, we establish the convergence rate
of the strong convergence in terms of the size of the space discretization $h$ (which is supposed to be comparable with the time step $\Delta t$).
In the case of non zero inflow-outflow boundary data, all results are new. The latter result is new also for the no-slip boundary conditions.
\end{abstract}
\\ \\
\noindent
{\bf AMS classifiaction:} 35Q30, 65N12, 65N30, 76N10, 76N15, 76M10, 76M12
\tableofcontents

\section{Introduction}
\subsection{State of the art and goals}
Evolution of the density $\mathfrak{r}=\mathfrak{r}(t,x)$ and velocity $\mathfrak{u}=\mathfrak{u}(t,x)$ through the time interval $[0,T)$, $T>0$, $t\in [0,T)$
in a bounded (Lipschitz) domain $\Omega$, $x\in \Omega$ of a viscous compressible fluid can be described by the
Navier-Stokes equations. {\em Unlike most of the theoretical  literature}, in  this paper, we consider the Navier-Stokes system with the {\em nonzero inflow-outflow 
boundary conditions}, which is, from the point of view of applications, a more appropriate setting, than the "standard" setting with the no-slip or
Navier boudary conditions. The equations read: 
\begin{equation} \label{NS1}
\partial_t {\mathfrak{r}} + \Div ({\mathfrak{r}} {\mathfrak{u}}) = 0,
\end{equation}
$$
\partial_t ({\mathfrak{r}} {\mathfrak{u}}) + \Div ({\mathfrak{r}} {\mathfrak{u}} \otimes {\mathfrak{u}}) + \Grad p({\mathfrak{r}}) = 
\Div \mathbb{S}(\Grad\mathfrak{u}), 
$$
$$
\mathbb{S}(\Grad{\mathfrak{u}})=\mu\Grad{\mathfrak{u}} +(\mu+\lambda){\rm div}{\mathfrak{u}}\mathbb{I},\ \mu>0,\, \lambda+\frac 23\mu>0 \footnote{
One usually writes $\mathbb{S}(\Grad\vu)$ in its frame indifferent form, $\mathbb{S}(\Grad\vu)=\mu(\Grad\mathfrak{u}+(\Grad\mathfrak{u})^T) +\lambda{\rm div}{\mathfrak{u}}\mathbb{I}$. Both writing are equivalent in the strong formulation but not in the weak formulation later and in the numerical scheme}
$$
with initial and boundary conditions,
\bFormula{ibc}
\mathfrak{r}(0)=\mathfrak{r}_0,\ \mathfrak{ u}(0)={\mathfrak{ u}}_0,\ \mathfrak{ u}|_{\partial\Omega}=\mathfrak{ u}_B,\
\mathfrak{ r}|_{\Gamma_{\rm in}}=\mathfrak{r}_B,
\eF
where (for the case of simplicity),\footnote{The regularity of initial conditions could be relaxed up to
$$
0<\mathfrak{r}_0\in L^\gamma(\Omega),\;{\bf{\mathfrak{ u}}_0}\in L^1(\Omega),\;\mathfrak{r}_0\mathfrak{u}_0^2\in L^1(\Omega).
$$
}
\bFormula{ru0}
0<\mathfrak{r}_0\in C(\overline\Omega),\;{\bf{\mathfrak{ u}}_0}\in C(\overline\Omega),
\eF
\bFormula{ruB}
0<\mathfrak{r}_B\in C^1(\overline\Omega),\;{\bf{\mathfrak{ u}}_B}\in C^2(\overline\Omega)
\eF
are given initial
 and boundary data, and
\bFormula{inB}
\Gamma^{\rm in}=\partial\Omega\setminus \Gamma^{\rm out},\;\Gamma^{\rm out}=\overline{\{x\in \partial\Omega|\;\mathfrak{u}_B(x)\cdot\vc n(x)\ge 0
\;\mbox{or $\vc n(x)$ is not defined}\}}
\eF
(where $\vc n$ is the outer normal to $\partial\Omega$),
is the inflow boundary. Here, "$\overline{\{\ldots\}}$" means the closure with respect to the trace topology of $R^3$ on $\partial\Omega$. 
 
In the above the pressure $p$
is a given function of density $\mathfrak{r}$ which will be specified later. A typical example which will be treated in the paper
includes pressure functions 
$$
p(\mathfrak{r})\approx \mathfrak{r}^\gamma,\;\mbox{for large values of $\mathfrak{r}$},\; \gamma>1.
$$

In the mathematical literature there is a large variety of numerical methods for solving efficiently compressible Navier-Stokes or Euler equations,
see \cite{Feist}, \cite{Kr}, \cite{TaZh}, \cite{DoFeFeKl}, \cite{FeistFeLu}, \cite{FeistFeLuWa}, \cite{GaGaLaHe}, \cite{GaHeLa}, \cite{EyGaHe} -- for the finite volume and/or a combination of finite volume and the finite element methods, e.g.
\cite{FeistFeSt}, \cite{FeistCeKu} -- for the dicontinuous Galerkin methods and e.g. \cite{HeLaMiTh}, \cite{HeKhLa}, \cite{GaHeLaMa},\cite{GHMN-MAC} - for the finite volume/finite difference schemes --  and the references quoted there. Although these methods give quite convincing results in numerical experiments and engineering applications, their mathematical analysis is less advanced, see the works of Tadmor et al. \cite{Ta0}, \cite{Ta1}, \cite{TaZh},
Gallouet et al. \cite{GaGaLaHe}, \cite{GaHeLa}, \cite{EyGaHeLa}, \cite{GaHeMaNo}, Jovanovic, Rohde \cite{JoRo}, \cite{Jo} for some pertinent results in this direction. Namely, the rigorous results
about the convergence of the numerical solutions to the (weak, Young measure-valued, strong ) solutions is in a short supply, and if available, then
only for the flows with the no slip or Navier boundary conditions. 

We name the work of Karper \cite{Ka} (and subsequent papers by Feireisl et al.
\cite{FeKaNo}, \cite{FeKaMi}) which provide the convergence for one particular (implicit) scheme (proposed originally by Karlsen and Karper \cite{KaKa0},\cite{KaKa1}, \cite{KaKa2} and which we will call Karper's scheme)
to a weak solution for large adiabatic coefficient ($\gamma>3$)--and so far no other scheme was proved to have this convergence property, meaning that  the proof is strongly dependent on the structure of the scheme. The convergence to a larger class of dissipative Young measure valued solutions (introduced
in Feireisl et al. in \cite{FeGwGwWi}) is available for a larger class of schemes and, under some circumstances, for the whole range of adiabatic exponents $\gamma>1$,
see \cite{FeLu}, \cite{FeLuMi}, \cite{FeLuMiSh}, \cite{HoSh}, and the proofs, again, are quite scheme structure dependent. Moreover, in these papers, the weak strong uniqueness property of the class of Young-measure valued solutions translates to the proofs of strong convergence of the numerical solutions to a strong solution, provided
the latter exists. This approach, however, does not provide the rate of convergence in terms of the size of the discretization, the so called error estimates. 

In 2016, Gallouet et al. \cite{GaHeMaNo} established the rigorous error estimates for the Karper's scheme and the adiabatic exponents $\gamma >3/2$. These
works were followed by other works establishing error estimates for a few other schemes, see \cite{GHMN-MAC}, \cite{MaNo}, \cite{FeHoMaNo}, \cite{FeLuNeNoSh}. We remark that the
proofs of the error estimates are much more "structure dependent" than the convergence proofs and that they are available only for the adiabatic 
exponents $\gamma>3/2$.

Let us notice that none of those results is available for the flows with the non zero inflow-outflow boundary conditions which we consider in this paper. No error estimates are known under the threshold $\gamma\le 3/2$ even for the no-slip boundary value case.

Recently, the authors of \cite{AbFeNo} suggested a notion of the dissipative solutions with a Reynolds defect, which is still consistent with the notion of the classical solution (any dissipative solution which is regular is a classical solution), whose class is larger than the class of Young measure valued solutions
suggested in \cite{FeGwGwWi}, but which still possesses the weak strong uniqueness property. Moreover, this notion was built within the class of flows 
allowing the non zero inflow-outflow boundary conditions. As claimed by Fjordholm et al. \cite{FjMiTa0}, \cite{FjMiTa1}, \cite{FjKaMiTa} at least in the context of hyperbolic conservation laws, this type of solutions may represent a more appropriate concept of solutions than the concept of weak solutions. This
may be the case for the compressible Navier-Stokes equations as well.

Note that existence of weak solutions to the Navier-Stokes system (\ref{NS1}--\ref{inB}) is known only in the case $\gamma>3/2$, and it is nowadays a standard
result, see \cite{FNP}
(and monographs by Lions \cite{Li}, Feireisl \cite{Fe},  Novotny, Straskraba \cite{NoSt}) for the no-slip boundary conditions. The same result for the
nonzero inflow-outflow boundary conditions is more recent, see \cite{ChJiNo}, \cite{ChoNoY}, \cite{KwNo} preceded by Girinon \cite{Gi}.
This is just to say that the proof of the convergence of any numerical scheme to a weak solution below the treshold $\gamma=3/2$ is nowadays unreachable, and any proof of this type below "Karper's" treshold $\gamma=3$ would be a tremendous progress. Even for the Karper's scheme
and adiabatic coefficients $\gamma>3$ the convergence of numerical solutions to weak solutions in the non zero inflow/outflow setting is an open problem. We do not wish to address it here, however, the methodology we develop would allow to address it in an efficient way.

Our goal in this paper is to {\em adapt the Karper's numerical scheme to the nonzero inflow-outflow boundary conditions, to prove existence of numerical solutions
for this scheme on the polygonal domains and to study the convergence of its solutions to the dissipative solutions with Reynolds defect} 
introduced in \cite{AbFeNo} for the {\em whole
scale of adiabatic exponents} $1<\gamma<\infty$. We did not choose  Karper's scheme for its computational performance, but rather because
it is recognized for its structural closeness to the original system which allows to adapt to it in a large extend the well known techniques
from the "continuous case". However, when doing this we have always in mind the universality of the approach and the transportability of the "main steps" of the method to other numerical schemes.

We consider the Karper's scheme in several versions: 1) The "classical" version with artificial density dissipation. 2) Without artificial density dissipation but with artificial pressure perturbation. 3) Combining both features mentioned above. 

Our approach is based on the following steps:

\begin{enumerate}
\item We modify Karper's scheme in order to be able to accommodate the {\em non-zero inflow/outflow boundary data}, we prove {\em existence
of numerical solutions} with the positive density.
\item We deduce from the scheme the balance of mass and the balance of energy and we use the balance of energy
to derive the uniform estimates.
\item We rewrite the {\em numerical continuity equation} and the {\em numerical momentum equation} in their {\em variational forms} letting
appear "remainders" vanishing as space $h$ and time $\Delta t$ discretizations tend to zero.
\item We rewrite the {\em numerical balance of energy} in a {\em consistent form compatible with the variational formulations of the
balance of mass and momentum}.
\item Using the {\em consistency formulation of the balance of mass, momentum, and energy}, and employing the uniform estimates, we show 
that the {\em numerical solution generates a dissipative solution with the Reynolds defect} in the sense of \cite{AbFeNo}. This part
provides {\em another proof of the existence of the dissipative solutions via a numerical scheme}, and by itself, it is of independent interest.
\item Last but not least, by using the weak strong uniqueness principle established in \cite{AbFeNo} for the limiting dissipative solution
we show that the {\em whole sequence of numerical solutions converges strongly  to the pointwise classical solution} of the original problem provided the latter exists.
\item Finally using again the consistency formulations of the continuity and momentum, and employing in addition the consistency formulation of the energy balance, we derive, by mimicking the continuous case, {\em the relative energy inequality} between the numerical and strong solutions
and deduce the {\em error estimates.}
\end{enumerate}

All convergence results are {\em "unconditional"}: to obtain them {\em we use only a priory estimates derived from the numerical scheme}. No a posteriori bounds
are needed, in contrast with the most mathematical literature about the subject.

To the best of our knowledge, {\em all results listed above} are {\em new in the case of general non-zero inflow/outflow boundary conditions.}  They hold true
also {\em in the case of zero boundary conditions}, in which situation {\em the error estimates} mentioned in the point 7 {\em is a new result}, as well. In the latter case,
the convergence to dissipative solutions with Reynolds stress is reminiscent to the convergence to Young measure valued solutions proved in
Feireisl et al. \cite{FeGwGwWi}. The advantage of the dissipative solutions with  the Reynolds stress may dwell in the fact that the convergence proofs
in this setting require slightly less from the structure of the scheme and thus, may be applicable to a larger class of numerical methods. The pertinence of this
statement, however, remains to be confirmed. All results in this paper are obtained through {\em an extensive  application of the functional analysis and theory of PDEs in the numerical
analysis which goes far beyond the standard approaches current in this domain of research.}

The paper is organized as follows. In the Introduction, we fix the structural assumptions, define the dissipative solutions with the Reynolds defect and report the weak-strong uniqueness principle in this class. In Section \ref{Nn} we describe the numerical setting and introduce the discrete functional spaces related to it. In Section \ref{NN} we suggest the adaptation of Karper's numerical scheme, see (\ref{N1}--\ref{N2}), and state the main results. Theorem \ref{TN1} about the
the existence of numerical solutions and Theorem \ref{TN2} consisting of three parts: 1) Convergence to dissipative solutions (item 1). 2) Strong convergence to strong solutions
(item 2a). 3) Error estimates related to the strong convergence (item 2b). Section \ref{Prel} is reporting necessary numerical and analytical preliminaries for the
treatement of the problem. Then, Theorem \ref{TN1} is proved in Section \ref{DP}, and the proof of Theorem \ref{TN2} is established through Sections \ref{SEst}--\ref{RE}. More exactly, we derive energy balance and uniform estimates for the numerical solutions in Section \ref{SEst} in Lemmas \ref{Lebalance} and 
\ref{Lestimates}. In Section \ref{NC} in Lemmas
\ref{Lc-consistency}, \ref{Lm-consistency} and in formulas (\ref{ee-v}--\ref{enb3}) we derive the convenient consistency formulations of continuity, momentum and energy balance.  
In Section \ref{SConv} we use the consistency formulation and uniform estimates to establish the weak convergence to dissipative solutions. This finishes the proof of the first part of Theorem \ref{TN2}. In Section \ref{SSt}, we use the weak strong uniqueness principle in the class of dissipative solutions to prove the strong convergence to a strong solution. This concludes
the proofs of part 2 of Theorem \ref{TN2}. Finally, in Section \ref{RE}, we establish the relative energy inequality between the numerical and strong solutions and conclude the proof of part 3 of Theorem \ref{TN2} by Lemma \ref{Error} which states the error estimates. We finish the paper by several concluding remarks in Section \ref{CR}.

We conclude this introductory part by a remark concerning the notation: The special functional spaces are always defined in the text. For the classical Lebesgue, Sobolev, Bochner spaces and their duals, we use the standard notation, see e.g. Evans \cite{Ev}. Strong convergence in a Banach space is always denoted "$\to$", while "$\rightharpoonup$" means the weak convergence and "$\rightharpoonup_*$" means the star-weak convergence.

\subsection{Some basic notions and assumptions}

We suppose
\bFormula{pres1}
p\in C^1[0,\infty),\; p(0)=0,\; p'(\vr)>0\,\mbox{for all $\vr>0$}
\eF
(and we extend $p$ by zero to the negative real line so that $p\in C(R)$ if needed).  
We associate to $p$ its Helmholtz function $H$,
\bFormula{H}
H(\vr)=\vr\int_1^\vr\frac{p(z)}{z^2}{\rm d}z\;\mbox{in particular}\;{ \vr H'(\vr)-H(\vr)=p(\vr).}
\eF
and we shall additionally suppose that\footnote{We remark that (\ref{H}) in combination with (\ref{pres1}) yields, in particular, $H(0)=0$.} 
\bFormula{pres4}
 H-\underline a p,\quad \overline a p-H\;\mbox{are convex functions for some $0<\underline a<\overline a$}.
\eF

It is to be noticed that an iconic example of the isentropic pressure $p(\vr)=a\vr^\gamma$, $a>0$, $\gamma>1$ complies with (\ref{pres1}--\ref{pres4}).

For a function $B\in C^1([0,\infty))$, we define
$$
E_B(\vr|r)=B(\vr)-B'(r)(\vr-r)-B(r)
$$
and we write $E_H=E$ for the particular $B=H$. For the further use, we also introduce the relative energy functional
\bFormula{calE}
{\cal E}(\vr,\vu|r,\vc U)=\intO{\Big(\frac 12\vr|\vu-\vc U|^2{\rm d}x+E(\vr|r)\Big)}
\eF
to measure a "distance" between vector fields $\vu$, $\vc U$ and positive scalar fields $\vr$ and $r$.

\subsection{Dissipative solutions with Reynolds defect}

\begin{Definition}{\rm [Dissipative solution with Reynolds defect] }\label{DD1}

The quantity $[\mathfrak{r}, \mathfrak{u}]$ is a \emph{dissipative solution} of the problem (\ref{NS1}--\ref{inB}) in $(0,T) \times \Omega$
if the following is satisfied:
\begin{enumerate}
\item  
\[
0\le\mathfrak{r} \in C_{\rm weak}([0,T]; L^\gamma(\Omega)) \cap L^\gamma(0,T; L^\gamma(\partial \Omega; |\mathfrak{u}_B \cdot \vc{n}|
{\rm d}S_x)),\ \ \mbox{for some}\ \gamma > 1, \footnote{We say that $r\in C_{\rm weak}([0,T];X)$, $X$ a Banach space, if
$r:[0,T]\mapsto X$ is defined on $[0,T]$, $r\in L^\infty(0,T;X)$ and
$<{\cal F},r>_{X^*,X}\in C[0,T]$ with any ${\cal F}$ in the dual space $X^*$ to $X$.}
\]
\[
 \mathfrak{v}=\mathfrak{u} - \mathfrak{u}_B  \in L^2(0,T; W^{1,2}_0(\Omega; R^3)),
\mathfrak{m} \in C_{\rm weak}([0,T]; L^{\frac{2 \gamma}{\gamma + 1}}(\Omega; R^3)),
\]
\bFormula{D1}
\mathfrak{m}=\mathfrak{r}\mathfrak{u}\;\mbox{a.e. in $(0,T)\times\Omega$},\; 
\mathfrak{r}{\mathfrak{u}}^2\in L^\infty(0,T;L^1(\Omega)),\;p(\mathfrak{r})\in L^1((0,T)\times\Omega));
\eF

\item The continuity equation
\begin{equation} \label{D2}
\left[ \intO{ {\mathfrak{r}} \varphi } \right]_{t = 0}^{t = \tau} +
\int_0^\tau \int_{\Gamma_{\rm out}} \varphi {\mathfrak{r}} {\mathfrak{u}}_B \cdot \vc{n} \ {\rm d}  S_x
+
\int_0^\tau \int_{\Gamma_{\rm in}} \varphi {\mathfrak{r}}_B {\mathfrak{u}}_B \cdot \vc{n} \ {\rm d} S_x
\end{equation}
$$
=
\int_0^\tau \intO{ \Big[ {\mathfrak{r}} \partial_t \varphi + {\mathfrak{r}} {\mathfrak{u}} \cdot \Grad \varphi \Big] } \dt,\ {\mathfrak{r}}(0,\cdot) = {\mathfrak{r}}_{0},
$$
holds for any $0 \leq \tau \leq T$, and any test function $\varphi \in C^1([0,T] \times \Ov{\Omega})$;

\item
There is a tensor
\[
 \mathfrak{R} \in L^\infty(0,T; \mathcal{M}^+(\Ov{\Omega}; R^{3 \times 3}_{\rm sym})), \footnote{This means that for any $\xi\in R^3$ and a. a. $t\in (0,T)$, $\xi^T\mathfrak{R}(t)\xi$ is a positive Radon measure on $\overline\Omega$
and $\forall\varphi\in C(\overline\Omega)$, $\int_{\overline\Omega}\varphi(x){\rm d}\mathfrak{R}(\cdot)
\in L^\infty((0,T))$.}
\]
where $ \mathcal{M}^+(\Ov{\Omega}; R^{3 \times 3}_{\rm sym})$ denotes the space of positively semi-definite tensor valued (Radon) measures on $\overline\Omega$,
 such that the integral identity
$$
\left[ \intO{ {\mathfrak{m}}  \cdot \phi } \right]_{t=0}^{t = \tau} =
\int_0^\tau \intO{ \Big[ {\mathfrak{r}} {\mathfrak{u}} \cdot \partial_t \phi + {\mathfrak{r}} {\mathfrak{u}} \otimes {\mathfrak{u}} : \Grad \phi
+ p({\mathfrak{r}}) \Div \phi - \mathbb{S}(\Grad{\mathfrak{u}}) : \Grad \phi \Big] }{\rm d}t
$$
\bFormula{D3}
+ \int_0^\tau \int_{{\Omega}} \Grad \phi : {\rm d} \mathfrak{R}(t) \ \dt,\; 
{\mathfrak{m}}(0,\cdot)= {\mathfrak{r}}_{0}{\mathfrak{u}}_0
\eF
holds for any $0 \leq \tau \leq T$ and any $\phi \in C^1([0,T] \times \Ov{\Omega}; R^3)$, $\phi|_{\partial \Omega} = 0$;
\item
The energy inequality
\begin{equation} \label{D4}
\intO{\left[ \frac{1}{2} {\mathfrak{r}} |{\mathfrak{v}}|^2 + H({\mathfrak{r}}) \right](\tau) } + D \int_{\Ov{\Omega}} {\rm d} {\rm Tr}[\mathfrak{R}(\tau)] 
+
\int_0^\tau \intO{ \mathbb{S}(\Grad {\mathfrak{u}}):\Grad\mathfrak{u} } \dt
\end{equation}
$$
+\int_0^\tau  \int_{\Gamma_{\rm out}} H({\mathfrak{r}})  {\mathfrak{u}}_B \cdot \vc{n} \ {\rm d} S_x \dt 
+
\int_0^\tau\int_{\Gamma_{\rm in}} H({\mathfrak{r}}_B)  {\mathfrak{u}}_B \cdot \vc{n} \ {\rm d}S_x \dt	
\leq 
\intO{\left[ \frac{1}{2} {\mathfrak{r}}_0 |{\mathfrak{v}}_0|^2 + H({\mathfrak{r}}_0) \right] } 
$$
$$
-\int_0^\tau \intO{ \left[ {\mathfrak{r}} {\mathfrak{u}} \otimes {\mathfrak{u}} + p({\mathfrak{r}}) \mathbb{I} \right]  :  \Grad {\mathfrak{u}}_B } \dt
+ \int_0^\tau   \intO{ {{\mathfrak{r}}} {\mathfrak{u}} \cdot \Grad {\mathfrak{u}}_B   \cdot {\mathfrak{u}}_B  }
\dt 
$$
$$
+ \int_0^\tau  \intO{ \mathbb{S}(\Grad\mathfrak{u}) : \Grad {\mathfrak{u}}_B } \dt -
\int_0^\tau  \int_{\Ov{\Omega}} \Grad {\mathfrak{u}}_B : {\rm d} \mathfrak{R}(t) \dt
$$
holds 
for a.a. $\tau\in (0,T)$ with $D = \min \left\{ {1}/{2}, {\underline{a}}/{3} \right\}$, 
${\mathfrak{v}}_0=\mathfrak{u}_0 - {\mathfrak{u}}_B$.
\end{enumerate}
\end{Definition}

Global existence of dissipative solutions with Reynolds defect  has been proved in \cite[Theorem 3.8]{AbFeNo}, in particular, under assumptions (\ref{ru0}--\ref{ruB}), (\ref{pres1}), (\ref{pres4}). The compatibility of these solutions with the classical formulation has been proved in \cite[Theorem 4.1]{AbFeNo}. For the purpose of this paper, we report the following weak--strong uniqueness theorem in the class of these solutions, see \cite[Theorem 6.3]{AbFeNo}.
\begin{Lemma} {\rm [{\rm Weak--strong uniqueness for dissipative solutions}]} \label{WUT1}
Let $\Omega \subset R^3$ be a bounded Lipschitz domain. 
Suppose that $p$  satisfies the hypotheses (\ref{pres1}--\ref{pres4}) and the initial and boundary conditions verify (\ref{ru0}--\ref{ruB}).
Let $[\mathfrak{r}, \mathfrak{u}]$ be a dissipative solution of problem (\ref{NS1}--\ref{ibc}) 
in the sense of Definition \ref{DD1}, and let $[r, \vc U=\vc V+\mathfrak{u}_B]$ be a strong solution of the same problem belonging to the 
class 
\[
\vc V \in C^1([0,T] \times \Ov{\Omega}; R^3), 
\, \Grad^2\vc V \in C([0,T] \times \Ov{\Omega}; R^{3}),
\]
\bFormula{tests}
\vc V|_{\partial\Omega}=0,\
r \in C^1([0, T] \times \Ov{\Omega}), \ \underline r:=\inf_{(0,T) \times \Omega} r > 0 
\eF
emanating from the same initial and boundary data. Then 
\[
\mathfrak{r} = r, \ \mathfrak{u} = \vc U \ \mbox{in}\ (0,T) \times \Omega,\ \mathfrak{R} = 0.
\]
\end{Lemma}

\section{Numerical setting}
\label{Nn}

\subsection{Mesh}
\label{mesh}

We suppose that the physical space is a \emph{polyhedral bounded domain} $\Omega \subset R^3$ that admits a \emph{tetrahedral} mesh ${\cal T}={\cal T}_h$;
the individual elements in the mesh will be denoted by $K=K_h \in {\cal T}$ (closed sets) and their gravity centers by $x_K$. Faces in the mesh are denoted as $\sigma=\sigma_h$ (close sets in $\R^2$) and their gravity centers by $x_\sigma$, whereas ${\cal E}={\cal E}_h$ is the set of all faces.\footnote{In the sequel
we shall omit in the notation the dependence on the "size" $h$ whenever there is no danger of confusion.} We also denote by ${\cal E}(K)$ the set of all faces
of $K\in {\cal T}$. The set of faces ${\sigma} \subset \partial \Omega$ is denoted ${\cal E}_{{\rm ext}}$, while ${\cal E}_{{\rm int}} = {\cal E} \setminus {\cal E}_{{\rm ext}}$. Finally, we suppose that the mesh fits to the inflow-outflow boundaries, meaning that
\bFormula{fit}
{\cal E}_{\rm ext}={\cal E}^{\rm in}\cup {\cal E}^{\rm out},\;
\overline{\Gamma^{\rm in}}=\cup_{\sigma\in {\cal E}^{\rm in}}\sigma,\;
\overline{\Gamma^{\rm out}}=\cup_{\sigma\in {\cal E}^{\rm out}}\sigma,
\eF
where ${\cal E}^{\rm in}$, ${\cal E}^{\rm out}$ are defined in (\ref{in-out}) later.

We denote by $h_K$ the diameter of $K$ and by $\mathfrak{h}_K$ the radius of the largest ball included in $K$. We call $h= \sup_{K\in{\cal T}}h_K$ the size of the mesh  and denote
$\mathfrak{h}= \inf_{K\in{\cal T}}\mathfrak{h}_K$. 

For two numerical quantities $a$, $b$, we shall write
\[
a \aleq b \ \mbox{if} \ a \leq c b, \ c > 0 \ \mbox{a constant}, \ a \approx b \ \mbox{if} \ a \aleq b \ \mbox{and} \ b \aleq a.
\]
Here, ``constant'' typically means a generic quantity independent of the size $h$ of the mesh and the time step $\Delta t$ used in the numerical scheme as well as other parameters as the case may be.

In addition, we require the mesh to be admissible in the sense of Eymard et al. \cite[Definition 2.1]{EyGaHe}:
\begin{enumerate}
\item For $K, L \in {\cal T}$, $K \ne L$, the intersection $K \cap L$, if non-empty, is either a vertex, or an edge, or a face $\sigma \in {\cal E}$. In the latter case, we write
$\sigma=K|L$.
\item There holds
$$
h\approx \mathfrak{h}
$$
\end{enumerate}

We denote by $\vc{n}_{\sigma, K}$ the unit normal to the face $\sigma\in {\cal E}(K)$
outwards to $K$. On the other hand we associate to each element $\sigma \in {\cal E}$ 
a fixed normal vector $\vc{n}={\vc n}_\sigma$. If $\sigma\in{\cal E}_{\rm ext}$ then
$\vc n_\sigma$ is always the outer normal to $\partial\Omega$.

\subsection{Piecewise constant finite elements}\label{Q}

We introduce the space
\[
Q (\Omega)=Q_h(\Omega) = \left\{ g \in L^1(\Omega) \ \Big| \ 
g|_K = a_K \in \tn R \right\}
\]
of piecewise constant functions along with the associated projection
\[
\Pi^Q=\Pi_h^Q: L^1(\Omega) \to Q(\Omega), \ \widehat g|_K:=\Pi^Q [g]|_K 
= g_K:= \frac{1}{|K|} \int_K{g}{\rm d}x.
\]

For a function $g\in Q(\Omega)$ and any $\sigma\in {\cal E}_{\rm int}$, we denote
\bFormula{gpm}
g^+_\sigma:=g^+_{\vc n_\sigma} = \lim_{\delta \to 0+ } g(x_\sigma + \delta \vc{n} ),\
g^-_\sigma:=g_{\vc n_\sigma} = \lim_{\delta \to 0+ } g(x_\sigma - \delta \vc{n} ).
\eF
Further, we define  the jumps and mean values over $\sigma$ 
(relative to $\vc n_\sigma$),
\begin{equation}\label{jump}
\ju{g}_\sigma=\ju{g}_{\sigma,\vc{n}}:= g^+_\sigma - g^-_\sigma,
\av{g}_\sigma:= \frac{1}{2} \left( g^+ + g^- \right).
\end{equation}

\subsection{Crouzeix-Raviart finite elements}\label{V}

A differential operator $D$ acting on the $x-$variable will be discretized as
\[
D_h v|_K = D (v|_K) \ \mbox{for any}\ v \ \mbox{differentiable on each element}\ K\in {\cal T}.
\]

The \emph{Crouzeix-Raviart finite element spaces} (see Brezzi and Fortin \cite{BreFor}, among others) are defined as
\bFormula{n1}
V (\Omega)=V_h(\Omega): = \Big\{ {v} \in L^1 (\Omega) \ \Big| \ {v}|_K = \mbox{affine function},
\eF
$$
 \ \int_\sigma v|_K{\rm dS}_x= \int_\sigma v|_L{\rm dS}_x \ \mbox{for any}\ \sigma =K|L\in {\cal E}_{\rm int} \Big\},
$$
together with
\bFormula{n1+}
V_{0} (\Omega)=V_{h,0} = \left\{ {v} \in V(\Omega) \ \Big| \ \int_{\sigma} {v} \ {\rm dS}_x = 0\ \mbox{for any}\ \sigma \in {\cal E}_{\rm ext} \right\}.
\eF

Next,
we introduce the associated projection
\begin{equation}\label{PiV}
\Pi_V=\Pi_h^V : W^{1,1}(\Omega) \to V (\Omega), \ {\tilde v}= \Pi^V[v]:=\sum_{\sigma\in{\cal E}}
v_\sigma\phi_\sigma,
\end{equation}
where
\[
v_\sigma= 
\frac 1{|\sigma|}\int_\sigma{ v }{\rm d}S_x 
\]
and $\{\phi_\sigma\}_{\sigma\in {\cal E}}\subset V(\Omega)$ is a basis in $V(\Omega)$
given by
$$
\frac 1{|\sigma'|}\int_{\sigma'}\phi_\sigma=\delta_{\sigma,\sigma'},\;(\sigma,\sigma')\in {\cal E}^2.
$$

\subsection{Convective terms, upwinds}\label{Upw}

Suppose, that we are given
 $\vu_B\in V(\Omega,R^3)$. We define
\begin{equation}\label{in-out}
{\cal E}^{\rm in}=\{\sigma\in {\cal E}_{\rm ext}\,|\,\vu_{B,\sigma}\cdot\vc n<0\},\quad
{\cal E}^{\rm out}=\{\sigma\in {\cal E}_{\rm ext}\,|\,\vu_{B,\sigma}\cdot\vc n\ge 0\}
\end{equation}

We  define for  any $\sigma\in{\cal E}_{\rm int}$ for any
$g\in Q(\Omega)$  its upwind value on any $\sigma\in{\cal E}_{\rm int}$,
\begin{equation}\label{gup}
g_\sigma^{\rm up}=\left\{
\begin{array}{c}
g_K\;\mbox{if $\sigma= K|L\in {\cal E}_{\rm int}$ and 
 $\vc u_\sigma\cdot\vc n_{\sigma,K}\ge 0$},\\
 g_L\;\mbox{if $\sigma= K|L\in {\cal E}_{\rm int}$ and 
 $\vc u_\sigma\cdot\vc n_{\sigma,K}< 0$}
\end{array}
\right\}.
\end{equation}

Finally, we associate to any face $\sigma\in{\cal E}_{\rm int}$
the \emph{upwind} operator 
${\rm Up}_{\sigma}[{ g}, \vu]:={\rm Up}_{\sigma,\vc n}[{ g}, \vu]$ defined as
\bFormula{Up}
{\rm Up}_{\sigma,\vc n}[{ g} ,\vu] = { g}^- [\vu_\sigma \cdot \vc{n}]^+ + { g}^+ [\vu_\sigma \cdot \vc{n}]^-,\ \mbox{where}\ [c]^+ = \max \{ c, 0 \}, \ [c]^-= \min \{ c, 0 \},
\eF
and to any face $\sigma\in {\cal E}(K)$ the specific flux $F_{\sigma,K}$ (outwards the element $K$) defined as
\bFormula{F}
F_{\sigma,K}[g,\vu]=g_\sigma^{\rm up}\vu_\sigma\cdot{\vc n}_{\sigma,K}.
\eF

\subsection{Time discretization}
For simplicity, we shall consider the constant time step $\Delta t>0$ where $T=N\Delta t$, $N\in \tn{N}$ and
we set 
\bFormula{Ik}
I_k=(\tau_{k-1}, \tau_k],\; \tau_k=k\Delta t,\; k\in \mathbb{Z}.
\eF
 Suppose that we have functions $v^k:\Omega\to R$,
$k=0,\ldots, N$. For convenience, we set $v^k(x)=v^0(x)$ if $k\le 0$ and $v^k(x)=v^{N}(x)$ if $k>N$, and
we introduce numbers
$$
D_tv^k(x)=\frac{v^k(x)-v^{k-1}(x)}{\Delta t}, \; k\in \mathbb{Z}.
$$
Finally, we define
\bFormula{pwt}
v(t,x)=\sum_{k\in\mathbb{Z}}1_{I_k}(t)v^k(x),\; \quad D_tv(t,x)=\sum_{k\in\mathbb{Z}}1_{I_k}(t)D_tv^k(x),
\eF
{
\bFormula{pwt+}
\utilde{v}(t,x)=\sum_{k\in\mathbb{Z}} 1_{I_k}(t)\Big(v^{k-1}(x)+(t-(k-1)\Delta t)D_tv^k(x)\Big),\;\mbox{so that
$\partial_t\utilde{v}(t,x)=D_tv(x)$. }
\eF 
}
In the sequel, we denote by $L_{\Delta t}(0,T; Q_h(\Omega)):=L(0,T;Q(\Omega))$ resp.
$L_{\Delta t}(0,T; V_h(\Omega)):=L(0,T;V(\Omega))$ the spaces of piecewise constant functions from $[0,T]$
to $Q(\Omega)$ and $V(\Omega)$, respectively (constant on each $I_k$ and extended by the value in $I_0$
to the negative real axes and by the value in $I_{N}$ to $[T,\infty)$).

\section{Numerical scheme, main results}
\label{NN}

Having collected the necessary preliminary material, we are ready to introduce the numerical scheme to solve the compressible Navier-Stokes system with the nonhomogenous boundary conditions.

\subsection{Numerical scheme}
We  are given the approximations of the initial and boundary conditions
\bFormula{N1}
\vr^0=\vr^0_{h} = \Pi^Q_h [\mathfrak{r}_0], \ \vu^0=\vu^0_{h} = \Pi^V[{\bf \mathfrak{u}}_0], \ \vr_B=\vr_{B,h}=\Pi^Q[\mathfrak{r}_B],\
\vu_B=\vu_{B,h}=\Pi^V[{\bf \mathfrak{u}}_B]
\eF
We are searching for
$$
\vr_{h,\Delta t}(t,x)=\sum_{k=1}^{N}\sum_{K\in{\cal T}}1_{I_k}(t)1_K(x)\vr^k_{K,h,\Delta t},
\;
\vu_{h,\Delta t}(t,x)=\sum_{k=1}^{N}\sum_{K\in{\cal T}}1_{I_k}(t)1_K(x)\vu^k_{K,h,\Delta t}
$$
where
\bFormula{N1+}
\vr^k\in Q(\Omega),
\vr^k_{h,\Delta t}> 0,\;\vu^k_{h,\Delta t}\in V(\Omega;R^3),\;\vc v^k=\vu^k-\vu_B\in
V_0(\Omega;R^3),\; k=1,\ldots,N
\eF
such that the following algebraic equations
(for the unknowns $\vr^k_K$, $\vu^k_\sigma$, $k=1,\ldots, N$, $K\in {\cal T}$,
$\sigma\in {\cal E}$) are satisfied:\footnote{In what follows,
we omit the indexes ``$h$'' and/or ``$\Delta t$'' and write simply $\vr^k$ instead of $\vr^k_{h,\Delta t}$,
$\vr$ instead of $\vr_{h,\Delta t}$,
etc.,  in order to avoid the cumbersome notation, whenever there is no danger of confusion.}
\begin{enumerate}
\item{\textsc{Approximation of the continuity equation}}
\bFormula{N2}
\intO{ D_t \vr^k \phi } + \sum_{K\in {\cal T}}\sum_{\sigma \in {\cal E}(K)\cap{\cal E}_{\rm int}}
\int_\sigma F_{\sigma,K}(\vr^k,\vu^k)\phi{\rm d}S_x
+
\sum_{\sigma\in{\cal E}^{\rm out}}\int_\sigma \vr^{k}\vu_{B,\sigma}\cdot\vc n_\sigma
\phi{\rm d}S_x
\eF
$$
+\sum_{\sigma\in {\cal E}^{\rm in}}\int_{\sigma}
\vr_{B}\vu_{B,\sigma}\cdot\vc n_{\sigma}\phi{\rm d}S_x
{ + \kappa h^\omega\sum_{\sigma\in{\cal E}_{\rm int}}\int_\sigma[[\vr^k]]_{\sigma}[[\phi]]_\sigma{\rm d}S_x}=0
$$
for all $\phi \in Q (\Omega)$, { where $\kappa=1$ or $\kappa=0$ and $\omega>0$.}
\item{\textsc{Approximation of the momentum equation}}
\bFormula{N3}
\intO{ D_t (\vr^k \widehat{\vv}^k) \cdot \widehat\phi }
+\sum_{K\in {\cal T}}\sum_{\sigma \in {\cal E}(K)\cap{\cal E}_{\rm int}}\int_\sigma F_{\sigma,K}(\vr^k\widehat{\vv}^k,\vu^k)\cdot\widehat\phi{\rm d}S_x
\eF
$$
{ + \sum_{K\in {\cal T}}\sum_{\sigma \in {\cal E}(K)}\int_\sigma\vr^k\widehat\vu^k\cdot{\vc n}_{\sigma,K}\vu_B\cdot\widehat\phi
{\rm d} S_x} 
+
\sum_{\sigma\in{\cal E}^{\rm out}}\int_\sigma \vr^{k}\vu_{B,\sigma}\cdot\vc n_\sigma\widehat\vv^k\cdot\widehat\phi{\rm d}S_x
+
\sum_{\sigma\in {\cal E}^{\rm in}}\int_\sigma \vr_B\vu_{B,\sigma}\cdot\vc n_\sigma\widehat\vv^k\cdot\widehat\phi{\rm d}S_x
$$
$$
+ \intO{ \Big(\mathbb{S}(\nabla_h\vu^k):\Gradh \phi 
-
p_h(\vr^k) \Divh \phi\Big) }
 + \kappa h^\omega \sum_{\sigma\in{\cal E}_{\rm int}}\int_\sigma[[\vr^k]]_{\sigma}\{\widehat\vv^k\}_{\sigma}[[\widehat\phi]]_\sigma{\rm d}S_x
 =0
$$
for any $\phi \in V_{0}(\Omega;R^3)$,\footnote{In agreement with (\ref{NS1}), here and in the sequel, $\intO{\mathbb{S}(\nabla_h\vu^k):\nabla_h \phi}$
means exactly $\intO{(\mu\nabla_h\vu^k:\nabla_h\phi+(\mu+\lambda){\rm div}_h\vu^k{\rm div}_h\phi)}.$ This form is important for the estimates, since the Korn
inequality does not hold in the  Crouzeix-Raviart finite element space.}
{ where
\bFormula{ph}
p_h(\vr)=p(\vr)+\tilde\kappa h^\eta\vr^2,\;\tilde\kappa\in\{0,1\},\;\eta>0.
\eF
}
\end{enumerate}

\bRemark{tr}
\begin{enumerate}
\item It is to be noticed that the background linear momentum $\vr\vu_B$ in the momentum equation is not "upwinded". If it were "upwinded"
we would loose the derivation of uniform estimates from the energy balance. This observation seems to have an universal character valid for 
the discretizations of the inflow/outflow problems via the finite volume methods, in general.
\item { The term   
$
{ \kappa h^\omega\sum_{\sigma\in{\cal E}_{\rm int}}\int_\sigma[[\vr^k]]_{\sigma}[[\phi]]_\sigma{\rm d}S_x}
$
and the last term at the last line
of the momentum equation create an artificial diffusion in the density field. The perturbation
$\tilde\kappa h^\eta \vr^2$ represents a "regularization" of pressure. Under certain assumptions on
$\eta$, the pressure regularization itself guarantees already the convergence of the scheme to the dissipative solutions with Reynolds defect 
(with or without
the artificial density diffusion) for the whole range $1<\gamma<\infty$. We have chosen this
particular setting to state the convergence theorem (Theorem \ref{TN2}). However, under certain hypotheses on $p$ and  for some ranges of $\gamma$'s, the artificial pressure perturbation and even the artificial diffusion are not needed. We shall discuss these situations in the Remark \ref{Rrest} after the main theorem.} 
\end{enumerate}
\eR

\subsection{Main results}

The first result postulates the existence of numerical solutions.

\bTheorem{N1}{\rm[ Existence of numerical solutions.]}
Under assumptions { (\ref{pres1}--\ref{pres4})}, (\ref{ru0}--\ref{ruB}),
the algebraic system (\ref{N1}--\ref{N3}) admits at least one solution $[\vr_{h,\Delta t},\vu_{h,\Delta t}]\in L_{\Delta t}(0,T;Q_h(\Omega)) \times 
L_{\Delta t}(0,T; V_h(\Omega))$. Any of its solution has a strictly positive density.
\eT

Theorem \ref{TN1} will be proved in Section \ref{DP}. Let us notice that uniqueness of numerical solutions to this scheme is not known.

The main result deals with the case $h \approx\Delta t$. It guarantees a (weak){\it convergence of a  subsequence of numerical solutions to a dissipative solution with Reynolds defect.} Moreover,
if the original problem admits a strong solution, than {\it the weak convergence of a subsequence becomes a strong convergence of the whole sequence}
to the strong solution. Even more, in the latter case the {\it convergence rate of the strong convergence} can be evaluated in terms
of a power of the discretization parameter $h$: in the other words we have an {\it error estimate.} 
\bTheorem{N2} {\rm [Main theorem: Convergence and error estimates.]}
Let $h=\Delta t$ and
{
\bFormula{range1}
\tilde\kappa=1,\;\mbox{and}\; \eta\in (0,2/3) \, \mbox{if $\kappa=0$}\quad or\quad \eta\in (0,\min\{2\omega,2/3\})\,\mbox{if $\kappa=1$}.
\eF
 }
Suppose that the pressure satisfies assumptions (\ref{pres1}--\ref{pres4}) and that the initial and boundary conditions verify
(\ref{ru0}--\ref{ruB}). Consider a sequence of numerical solutions $[\vr_h,\vu_h]\in L_{h}(0,T; $ $Q_h(\Omega))\times L_{h}(0,T; V_h(\Omega))$ of the problem (\ref{N1}--\ref{N3}).\footnote{To simplify notation, we denote
$\vr_{h,h}$ by $\vr_h$ or even $\vr$, $\vu_{h,h}$ by $\vu_h$ or even $\vu$ etc., if there is no danger of confusion.} Then we have:
\begin{enumerate}
\item There exists a subsequence $[\vr_h,\vu_h]$ (not relabeled) such that
{ \bFormula{A1}
\vr_h\rightharpoonup\mathfrak{r}\;\mbox{in $L^\infty(0,T;L^\gamma(\Omega))$},\; \gamma=1+ 1/\overline a,
\eF
\bFormula{A2}
\vv_h\rightharpoonup\mathfrak{v}\;\mbox{in $L^2(0,T;L^6(\Omega))$},
\eF
\bFormula{A3}
\vr_h\widehat\vu_h\rightharpoonup_*\mathfrak{r}\mathfrak{u}\;\mbox{in $L^\infty(0,T;L^{\frac{ 2\gamma}{\gamma+1}}(\Omega))$}
\eF
(where $\mathfrak{r}\mathfrak{u}\otimes\mathfrak{u}, p(\mathfrak{r})\in L^\infty(0,T;L^1(\Omega))$, $\mathfrak{u}=\mathfrak{v}+\mathfrak{u}_B$)
and
\bFormula{A4}
\vr_h\widehat\vu_h\otimes\widehat\vu_h+p(\vr_h)\mathbb{I}\rightharpoonup_*\overline{(\mathfrak{r}\mathfrak{u}\otimes\mathfrak{u}+p(\mathfrak{r})\mathbb{I})}
\;\mbox{in $L^\infty(0,T;{\cal M}(\overline\Omega; R^{3\times 3}_{\rm sym}))$}
\eF
where 
}
$$
[\mathfrak{r},\mathfrak{u}]\;\mbox{ and }\;\mathfrak{R}=\overline{(\mathfrak{r}\mathfrak{u}\otimes\mathfrak{u}+p(\mathfrak{r})\mathbb{I})}
- (\mathfrak{r}\mathfrak{u}\otimes\mathfrak{u}+p(\mathfrak{r})\mathbb{I})
$$
is a dissipative solution of the problem (\ref{NS1}--\ref{inB}) in the sense of Definition \ref{DD1}.
\item  { If the problem (\ref{NS1}--\ref{inB}) admits a strong solution $[r,\vc U]$ in the class (\ref{tests}) then
there holds:}
\begin{enumerate}
 \item The (weak) limit
$[\mathfrak{r},\mathfrak{u}]$ in (\ref{A1}--\ref{A2}) is equal to $[r,\vc U=\vc V+\mathfrak{u}_B]$ and  the defect $\mathfrak{R}=0$. In this case
{\it the whole sequence} $[\vr_h,\vu_h]$ converges to $[r,\vc U]$ in the sense (\ref{A1}--\ref{A3}),
and, moreover,
{
\bFormula{A4+}
p_h(\vr_h)\rightharpoonup p({r})\;\mbox{ in $L^1((0,T)\times\Omega)$},
\eF
\bFormula{A5}
\vr_h\widehat\vu_h^2\rightharpoonup { r}{{\vc U }}^2\;\mbox{ in $L^1((0,T)\times\Omega)$}.
\eF 
}

\item  
 There exists
 $\alpha=\alpha(\eta)>0$ (in the case $\kappa=0$) and $\alpha=\alpha(\eta,\omega)>0$ (in the case $\kappa=1$),
and a positive number $c$ dependent on 
\bFormula{C}
\underline r:=\inf_{(0,T)\times\Omega}r,\;\overline r:=\sup_{(0,T)\times\Omega}r,\; \|\Grad r,\partial_t r, \vc V,\Grad\vc V, \Grad^2\vc V,\partial_t \vc V, \, \mathfrak{u}_B,
\Grad\mathfrak{u}_B, \Grad^2\mathfrak{u}_B\|_{C([0,T]\times\overline\Omega)}
\eF
such that
\bFormula{error1}
\Big[{\cal E}\Big(\vr_h,\widehat\vv_h| r,\vc V\Big)\Big]_0^\tau +\int_0^\tau\Big(\|\vu_h-\tilde{\vc U}\|_{L^2(\Omega)}^2+
\|\nabla_h(\vu_h-\tilde{\vc U})\|_{L^2(\Omega)}^2\Big){\rm d}t \aleq c \,h^\alpha. 
\eF
\end{enumerate}
\end{enumerate}
\eT
{
\bRemark{tr1}
\begin{enumerate}
\item Theorem \ref{TN2} holds, in particular, for the isentropic pressure
$p(\vr)=a\vr^\gamma$, $a>0$, $1<\gamma<\infty$.
\item
{ Formulas (\ref{A4+}) and (\ref{A5}) imply strong convergence in many particular pertinent situations: For example, for the "iconic" isentropic pressure $p(\vr)=a\vr^\gamma$, $a>0$ we have 
$$
\vr_h\to r\;\mbox{in $L^\gamma((0,T)\times\Omega)$},
$$
$$
\sqrt{\vr_h}\widehat\vu_h\to \sqrt r \vc U\;\mbox{in $L^2((0,T)\times\Omega)$}
$$
for the whole sequence of numerical approximations.} Moreover, formula (\ref{error1}) implies, in addition,
$$
\vu_h\to\vc U,\;\Gradh\vu_h\to\Grad\vc U\;\mbox{in $L^2((0,T)\times\Omega)$}
$$
for the whole sequence.
\item { Suppose now that $\tilde\kappa=0$. In this situation, we shall suppose for the pressure $p$ in addition to (\ref{pres1}--\ref{pres4}) also 
\bFormula{pres2}
\forall \vr>0,\;\underline\pi\vr^{\mathfrak{a}-1}+\underline p\vr^{\gamma-1} \le  p'(\vr),
\; 1<\mathfrak a\le \min\{\gamma,2\}, \;\gamma>1,\;\mbox{where $0\le\underline\pi$, $0<\underline p$.}
\eF
{\it In this setting, Theorem \ref{TN2} continues still to hold in the following situations:}}
\begin{enumerate}
\item If in (\ref{pres2}), $\underline\pi>0$, and
\bFormula{range2}
\kappa\in \{0,1\},\;\omega>0, \;2\le \gamma< \infty.
\eF
In this case, if  $\kappa=0$, the number $\alpha>0$ in (\ref{error1}) depends solely on $\gamma$.
\item If  $3/2<\gamma<2$, and in addition to (\ref{pres1}--\ref{pres4}), (\ref{pres2})
\bFormula{pres3}
\vr\mapsto \frac{p'(\vr)}\vr\;\mbox{ is nonincreasing on $(0,\infty)$},
\eF
provided
\bFormula{range3}
\kappa\in \{0,1\}, \;\omega>0.
\eF
In this case, in (\ref{error1}), $\alpha=\alpha(\gamma)$ if $\kappa=0$ and
$\alpha=\alpha(\gamma,\omega)$ if $\kappa=1$. 

We notice that this
case includes the isentropic pressure $p(\vr)=a\vr^\gamma$, $a>0$.
\item If $ 6/5< \gamma \le 3/2$ and if, in addition to (\ref{pres1}--\ref{pres4}), the condition (\ref{pres3})
holds, provided
\bFormula{range4}
\kappa=1,\; \omega\in (0,2(\gamma-1)).
\eF
In this case, in (\ref{error1}), $\alpha=\alpha(\gamma,\omega)$.

We notice that  this
case includes the isentropic pressure $p(\vr)=a\vr^\gamma$, $a>0$.
\end{enumerate}
\item The error estimates derived in Theorem \ref{TN2} require $\Delta t\approx h$, and they may not be optimal. They certainly are not optimal
in the range   $\gamma>3/2$ when, at least in the case of no-slip boundary conditions, "better" estimates (which do not require
the condition $\Delta t \approx h$) can be derived by using another approach, see Gallouet at al. \cite{GaHeMaNo}. For the general non zero inflow/outflow setting,
this is still an interesting open problem. We refer the reader to the paper \cite{MaNo} for the discussion on the
"optimality" of the error estimates.
\item Local in time existence of strong solutions to problem (\ref{NS1}--\ref{inB}) notably with non zero inflow outflow
boundary conditions is discussed in Valli, Zajaczkowski \cite{VaZa}.
\end{enumerate}
\eR
}

Theorem \ref{TN2} will be proved through Sections \ref{SEst}--\ref{RE}. Some hints how to modify the proof of its variants
mentioned in the third item of Remark \ref{Rtr1} are available in Remarks \ref{Rrest}, \ref{RRc-consistency}, \ref{RRm-consistency}.

\section{Preliminaries}\label{Prel}
\subsection{Preliminaries from numerical analysis} \label{Preln}

 In this part, we recall  several classical inequalities related to the discrete functional spaces which will be used throughout the paper.

\subsubsection{Some useful elementary inequalities}

We recall Jensen's inequalities
\bFormula{JensenV}
\begin{array}{c}
\|\widehat v\|_{L^q(K)}\aleq \|v\|_{L^q(K)}\;\mbox{for all $v\in L^q(K)$, $1\le q\le\infty$},\\
\|v_\sigma\|_{L^q(\sigma)}\aleq \|v\|_{L^q(\sigma)}\;\mbox{for all $v\in L^q(\sigma)$, $1\le q\le\infty$}.
\end{array}
\eF
together with the error estimate
\bFormula{errorQ}
\forall v\in W^{1,q}(\Omega),\; 
\left\{\begin{array}{c}
        \left\| v - \widehat v \right\|_{L^q(K)} \aleq h \| \Grad v \|_{L^q(K; R^3)}\;\forall v\in W^{1,q}(K)\\
       \left\| v - \widehat v \right\|_{L^q(\Omega)} \aleq h \| \Grad v \|_{L^q(\Omega; R^3)}\; \forall v\in W^{1,q}(\Omega)
       \end{array}
       \right\},\; 1\le q\le \infty.
\eF

We also recall the Poincar\'e type inequalities on the  mesh elements,
\bFormula{PoincareV-} 
\begin{array}{c}
\|v-v_\sigma\|_{L^q(K)}\aleq h \| \Grad v \|_{L^q(K; R^{3})}, \ \forall\sigma\in {\cal E}(K),\,v\in W^{1,q}(K),\,1\le q\le\infty,\\
\|v-v_K\|_{L^q(K)}\aleq h \| \Grad v \|_{L^q(K; R^{3})}, \ \forall K\in {\cal T},\,v\in W^{1,q}(K),\,1\le q\le\infty.
\end{array}
\eF

\subsubsection{Properties of piecewise constant and  Crouzeix-Raviart finite elements}\label{Prop}

We report the estimates of jumps on mesh elements,
\bFormula{PoincareV} 
\begin{array}{c}
\|[[\widehat v]]_{\sigma=K|L}\|_{L^q(\sigma)}\aleq h \| \nabla_h v \|_{L^q(K\cup L; R^{3})}, \ \forall v\in V(\Omega),\,1\le q\le\infty,\\
\|v|_K-v|_L \|_{L^q(K\cup L)}\aleq h \| \nabla_h v \|_{L^q(K\cup L; R^{3})}, \ \forall \sigma=K|L,\, v\in V(\Omega),\,1\le q\le\infty,
\end{array}
\eF
see  Gallouet et al. \cite[Lemma 2.2]{GaHeLa}.

Further we recall a global version of the Poincar\'e inequality on $V(\Omega)$
\bFormula{globalV}
\|v-\widehat v\|_{L^q(\Omega)}\aleq h\|\nabla_h v\|_{L^q(\Omega)}
\;\mbox{for all $v\in V(\Omega)$, $1\le q<\infty$},
\eF  
along with the global error estimate 
\bFormula{errorV}
\left\| v - \tilde v  \right\|_{L^q(\Omega)} +
h \left\| \Gradh \left( v - \tilde v \right) \right\|_{L^q(\Omega;R^3)}    \aleq  h^m \left\| v \right\|_{W^{m,q}
(\Omega)}
\eF 
for any $v \in W^{m,q} (\Omega)$, $m = 1,2,\  1 \le q \le \infty$,
see Karper \cite[Lemma 2.7]{Ka} or Crouzeix and Raviart \cite{CrRa}.

Next, we shall deal with the Sobolev properties of piecewise constant and Crouzeix-Raviart finite elements.
To this end we introduce a discrete (so called broken) Sobolev $H^{1,p}$-(semi)norm on $Q(\Omega)$,
\bFormula{?}
\| g \|_{Q^{1,p}(\Omega)}^p = \sum_{\sigma \in {\cal E}_{\rm int}} \int_{\sigma} \frac{ \ju{g}_\sigma^p }{h^{p-1}} \ {\rm dS}_x,\; 1\le p<\infty.
\eF
and a discrete (so called broken) Sobolev $H^{1,p}$-(semi)norm on $V(\Omega)$,
\bFormula{EM2}
\| {v} \|_{V^{1,p}(\Omega)}^p = \intO{  |\Gradh {v} |^p }.
\eF

 Related to the $Q^{1,p}$-norm, we report the following discrete Sobolev and Poincar\'e type inequalities: We have
\bFormula{SobolevQ}
\forall g\in Q(\Omega),\;\| g \|_{L^{p^*}(\Omega)} \aleq \| {g} \|_{Q^{1,p}(\Omega)} 
+ \| g \|_{L^p(\Omega)},
\eF
where $1\le p^*\le \frac{3p}{3-p}$ if $1\le p<3$, and $1\le p^*<\infty$, if $p\ge 3$,
see {Bessemoulin-Chatard et al. \cite[Theorem 3]{ChHiDr}.} 
Likewise, we have the Discrete Sobolev inequality,
\bFormula{SobolevV}
\forall v\in V(\Omega),\,
\left\{\begin{array}{c}
\| {v} \|_{L^{p^*}(K)} \aleq \| \Grad{v} \|_{L^p(K)}+
\|v\|_{L^p(K)},\\
\| {v} \|_{L^{p^*}(\Omega)} \aleq \| {v} \|_{V^{1,p} (\Omega)}+
\|v\|_{L^p(\Omega)}
\end{array}\right\},\; 1\le p\le\infty,
\eF
\begin{equation}\label{SobolevV0}
\forall v\in V_0(\Omega),\,\| {v} \|_{L^{p^*}(\Omega)} \aleq \| {v} \|_{V^{1,p}(\Omega)},
\end{equation}
 see e.g. \cite[Lemma 9.3]{GaHeMaNo}. 

Finally, we report the well-known identities for the Crouzeix-Raviart finite elements,
\bFormula{vv1}
\intO{ \Divh \tilde{\vc{u} }\ w } = \intO{ \Div \vc{u} \ {w} } \ \mbox{for any} \ {w} \in Q(\Omega)\ \mbox{and}\,\vc{u}\in V(\Omega;R^3),
\eF
{
\bFormula{proj}
\intO{ \Gradh v \otimes \Gradh \tilde\varphi } = \intO{\Gradh v \otimes \Grad \varphi } \ \mbox{for all}\ v \in V_{h}(\Omega),\
\varphi \in W^{1,1} (\Omega),
\eF
see \cite[Lemma 2.11]{Ka}.
}

\subsubsection{Trace and ``negative'' $L^p-L^q$ estimates for finite elements}\label{TNE}

We start by the classical trace estimate,
\bFormula{trace}
\| v \|_{L^q(\partial K)} \aleq  \frac{1}{h^{1/q}} \left( \| v \|_{L^q(K)} + h \| \Grad v \|_{L^q(K; R^3)} \right),\
1 { \leq q} \le \infty
\ \mbox{for any}\ v \in C^1({ \overline K}),
\eF
The following can be easily  obtained from the previous one by means of the scaling arguments:
\bFormula{traceFE}
\| w \|_{L^q(\partial K)} \aleq \frac{1}{h^{1/q}} \| w \|_{L^q(K)} \ \mbox{for any}\ 1 \leq q { \le} \infty,
\  w \in P_m,
\eF
where $P_m$ denotes the space of polynomials of order $m$.

In a similar way, from the local estimate
\bFormula{inter}
\| w \|_{L^p(K)} \aleq h^{3 \left( \frac{1}{p} - \frac{1}{q} \right) } \| w \|_{L^q(K)} \ 1 \leq q < p \leq  \infty, \ w \in P_m,
\eF
we deduce the global version
\bFormula{interG}
\| w \|_{L^p(\Omega)} \leq c h^{3 \left( \frac{1}{p} - \frac{1}{q} \right) } \| w \|_{L^q(\Omega)} \ 1 \leq q < p \leq  \infty,
\ \mbox{for any} \ w|_K \in P_m (K), \ K \in {\cal T}.
\eF
In particular, for a piecewise constant function in $(0,T=N\Delta t]$, $a(t)=\sum_{k=1}^{N} a_k 1_{I_k}(t)$, where $I_k=((k-1)\Delta t,k\Delta t]$ we have
\bFormula{time}
\| w \|_{L^p(I_k)} \aleq (\Delta t)^{ \left( \frac{1}{p} - \frac{1}{q} \right) } \| w \|_{L^q(I_k)} \ 1 \leq q < p \leq  \infty,
\eF
and
\bFormula{timeG}
\| w \|_{L^p(0,T)} \aleq  (\Delta t)^{ \left( \frac{1}{p} - \frac{1}{q} \right) } \| w\|_{L^q(0,T)} \ 1 \leq q < p \leq  \infty.
\eF
\subsubsection{Some formulas related to upwinding}
The following formulas can be easily verified by a direct calculation, see e.g. \cite[Section 2.4]{FeKaNo}
\begin{enumerate}
 \item Local conservation of fluxes:

\begin{equation}\label{Up1}
 \begin{array}{c}
  \forall \sigma=K|L\in {\cal E}_{\rm int},\;F_{\sigma,K}[g,\vu]=-F_{\sigma,L}[g,\vu],\\
  \forall \sigma \in {\cal E}_{\rm int},\ {\rm Up}_{\sigma,\vc n}[{ g}, \vu]=-
  {\rm Up}_{\sigma,-\vc n}[{ g}, \vu]\ \mbox{and}\ [[g]]_{\vc n_\sigma}=-[[g]]_{-\vc n_\sigma}.
 \end{array}
\end{equation}
\item For all $r,g\in Q(\Omega)$, $\vu\in V(\Omega,R^3)$,
\begin{equation}\label{Up2}
 \sum_{K\in {\cal T}}r_K\sum_{\sigma\in {\cal E}_{\rm int}}|\sigma| F_{\sigma,K}[g,\vu]=
  -\sum_{\sigma\in{\cal E}_{\rm int}} {\rm Up}_{\sigma}[g,\vu][[r]]_\sigma.
  \end{equation}
\item For all $r,g  \in Q(\Omega)$, $\vu,\vu_B\in V(\Omega;R^3)$,  $\vu-\vu_B\in V_0(\Omega;R^3)$,
 $\phi \in C^1(\Ov{\Omega})$,
\bFormula{Up3}
\intO{ g \vu \cdot \Grad \phi } = -\sum_{K\in {\cal T}}\sum_{\sigma\in {\cal E}(K)\cap{\cal E}_{\rm int}} F_{\sigma,K}[g,\vu] r
\eF
$$
+ \sum_{K\in{\cal T}} \sum_{\sigma\in {\cal E}(K)\cap{\cal E}_{\rm int}}\int_\sigma (r -  \phi ) \ \ju{ g }_{\sigma,\vc n_{\sigma,K}} \ [\vu_\sigma \cdot \vc{n}_{\sigma,K}]^- {\rm d}S_x 
$$
$$
+ \sum_{K \in {\cal T}} \sum_{\sigma\in {\cal E}(K)}\int_\sigma  g (\vu -
\vu_\sigma ) \cdot \vc{n}_{\sigma,K}(\phi-r){\rm d}S_x
+ \intO{ (r - \phi) g \Divh \vu }
$$
$$
+  \sum_{\sigma\in {\cal E}^{\rm out}}\int_\sigma
g\vu_{B,\sigma}\cdot \vc n_{\sigma}(\phi-r){\rm d}S_x
+  \sum_{\sigma\in {\cal E}^{\rm in}}\int_\sigma
g\vu_{B,\sigma}\cdot \vc n_{\sigma}(\phi-r){\rm d}S_x.
$$
\end{enumerate}

\subsection{Preliminaries from mathematical analysis}

Here we recall some elements of mathematical analysis needed in the proofs. 
First is the Shaeffer's fixed point theorem, see e.g. Evans \cite[Theorem 9.2.4]{Ev}.

\begin{Lemma}\label{fpoint}{\rm [Shaeffer's fixed point theorem]}
Let $\Theta:Z\mapsto Z$ be a continuous  mapping defined on a finite dimensional normed vector space $Z$.
Suppose that the set
$$
\{ z\in Z\,|\, z=\Lambda\Theta(z),\;\Lambda\in[0,1]\}
$$
is non empty and bounded. Then there exists $z\in Z$ such that 
$$
z=\Theta(z).
$$
\end{Lemma}

The next result is a version of Div-Curl lemma, see Lemma 8.1 in \cite{AbFeNo}.

\begin{Lemma} \label{Lemma1}
Let $O = (0,T) \times \Omega$, where $\Omega \subset R^d$, $d\ge 2$ is a bounded domain. Suppose that 
\[
r_n \rightharpoonup r \ \mbox{weakly in}\ L^p(O), \ v_n \rightharpoonup v \ \mbox{ in}\ L^q(O),\ p > 1, q > 1, 
\]
and
\[
r_n v_n \rightharpoonup w \ \mbox{ in}\ L^r(O), \ r >  1.
\]
In addition, let 
\[
\partial_t r_n = \Div \vc{g}_n + h_n\ \mbox{in}\ \mathcal{D}'(O),
\]
\[
 \| \vc{g}_n \|_{L^s(O; R^d)} \aleq 1,\ s > 1,\ 
h_n \ \mbox{precompact in}\ W^{-1,z}(O),\ z > 1, 
\]
and
\[
\left\| \Grad v_n \right\|_{\mathcal{M}(\overline O; R^d)} \aleq 1 \ \mbox{uniformly for}\ n \to \infty.
\]

Then 
\[
w = r v \ \mbox{a.a. in}\ O.
\]

\end{Lemma}

Finally, the last two lemmas are well known results from convex analysis, see e.g. Lemma 2.11 and Corollary 2.2 in
Feireisl \cite{Fe}.

\begin{Lemma}\label{Lemma2}
Let $O \subset R^d$, $d\ge 2$, be a measurable set and $\{ \vc{v}_n
\}_{n=1}^{\infty}$ a sequence of functions in $L^1(O; R^M)$ such
that
\[
\vc{v}_n \rightharpoonup \vc{v} \ \mbox{ in}\ L^1(O; R^M).
\]
Let $\Phi: R^M \to (-\infty, \infty]$ be a lower semi-continuous
convex function.

Then $\Phi(\vc v):O\mapsto R$ is integrable and
$$
\int_{O} \Phi(\vc v){\rm d} x\le \liminf_{n\to\infty} \int_{O} \Phi(\vc v_n){\rm d} x.
$$
\end{Lemma}

\begin{Lemma}\label{Lemma3}

Let $O \subset R^d$, $d\ge 2$ be a measurable set and $\{ \vc{v}_n
\}_{n=1}^{\infty}$ a sequence of functions in $L^1(O; R^M)$ such
that
\[
\vc{v}_n \rightharpoonup \vc{v} \ \mbox{ in}\ L^1(O; R^M).
\]
Let $\Phi: \R^M \to (-\infty, \infty]$ be a lower semi-continuous
convex function such that $\Phi (\vc{v}_n ) \in L^1(O)$ for
any $n$, and
\[
\Phi (\vc{v}_n) \rightharpoonup \Ov{\Phi(\vc{v})}
\ \mbox{in}\ L^1(O).
\]

Then

\bFormula{FNeb4}
\Phi (\vc{v})  \leq \Ov{\Phi (\vc{v})} \ \mbox{a.e. in}\ O.
\eF

If, moreover, $\Phi$ is strictly convex on an open convex set
$U \subset R^M $, and
\[
\Phi(\vc{v}) = \Ov{\Phi (\vc{v})} \ \mbox{a.e. on}\ O,
\]
then

\bFormula{FNeb5}
\vc{v}_n  (\vc{y}) \to \vc{v} (\vc{y})
\ \mbox{for a.a.}\ \vc{y} \in \{
\vc{y} \in O \ | \ \vc{v}(\vc{y}) \in U \}
\eF
extracting a subsequence as the case may be.
\end{Lemma}

\section{Existence of solutions for the discrete problem}\label{DP}
\subsection{Renormalization of the continuity equation}\label{RC}
Using in the discrete continuity equation test function $\phi\approx \phi B'(\vr^k)$,
$\phi\in Q(\Omega)$ we obtain after an elementary but laborious calculation 
the discrete renormalized continuity equation. This result is formulated
in the following lemma.
\bLemma{rcontinuity}{\rm [Renormalized discrete continuity equation]}
Let $B\in C^1(\tn R)$. Then
\bFormula{r1}
\intO{ D_t B(\vr^k) \phi } 
+  \sum_{K\in {\cal T}}\sum_{\sigma \in {\cal E}(K)\cap{\cal E}_{\rm int}}\int_\sigma F_{\sigma,K}(B(\vr^k),\vu^k)\phi{\rm d}S_x
\eF
$$
{ + \kappa h^\omega\sum_{\sigma\in {\cal E}_{\rm int}}\int_\sigma[[\vr^k]]_\sigma[[\phi B'(\vr^k)]]_{\sigma}{\rm d}S_x}
-\intO{\phi\Big(B(\vr^k)-\vr^k B'(\vr^k)\Big)\Divh\vu^k}
$$
$$
{ + \frac 1{\Delta t}\intO{E_B(\vr^{k-1}|\vr^{k})\phi}}+
\sum_{\sigma\in {\cal E}_{\rm int}}\int_\sigma E_B(\vr_\sigma^{k,+}|\vr^{k,-}_\sigma)
|\vu^k_\sigma\cdot\vc n|{ \phi^+}{\rm d S}_x
$$
$$
+
\sum_{\sigma\in{\cal E}^{\rm out}}\int_\sigma B(\vr^{k})\vu_{B,\sigma}\cdot\vc n_\sigma\cdot\phi{\rm d}S_x+ \sum_{\sigma\in {\cal E}^{\rm in}}\int_{\sigma}
E_B(\vr_{B}|\vr^k)|\vu_{B,\sigma}\cdot\vc n_{\sigma}|\phi{\rm d}S_x
=\sum_{\sigma\in {\cal E}^{\rm in}}\int_{\sigma}
B(\vr_{B})|\vu_{B,\sigma}\cdot\vc n_{\sigma}|\phi{\rm d}S_x
$$
for all $\phi \in Q (\Omega)$. 
\eL
\subsection{Positivity of density}\label{PD}

It is enough to show that $\vr^k>0$ provided $\vr^{k-1}>0$. (Recall that $\vr_B>0$.)

We take in (\ref{r1}) $\phi=1$ and $B$  the non negative convex function
$B(\vr)=\max\{-\vr,0\}$ (resp. any of its non negative $C^1(\R)$ convex approximations
$B_\ep$ such that $B_\ep(s)\to B(s)$, $B_\ep'\to B'(s)$ for all
$s\neq 0$). We get,
$$
\intO{B_\ep(\vr^k)}\le \intO{B_\ep(\vr^{k-1})} 
+\Delta t\Big(
\sum_{\sigma\in {\cal E}^{\rm in}}\int_\sigma B_\ep(\vr_B)|\vu_{B,\sigma}
\cdot\vc n_\sigma|{\rm d}S_x+\intO{(B_\ep(\vr^k)-\vr^kB_\ep'(\vr^k))\Divh\vu^k}\Big).
$$
Seeing that $B(s)-sB'(s)=0$ for all $s\neq 0$, $B(\vr^{k-1})=0$, $B(\vr_B)=0$, we deduce from the above that $B(\vr^k)=0$, i.e. $\vr_k\ge 0$.

Suppose that $\Omega_0:=\{x\in\Omega|\vr^k(x)=0\}=\cup_{K\in\tilde{\cal T}}K$,
where $\emptyset\neq \tilde{\cal T}\subset{\cal T}$. Taking in the continuity equation (\ref{N2}) test function
$\phi=1_{\Omega_0}$, we get
$$
-\int_{\Omega_0}\vr^{k-1}{\rm d}x
=-\sum_{K\in \tilde{\cal T}}
\sum_{\sigma\in\partial\Omega_0\cap{\cal E}_{\rm int}\cap{\cal E}(K)}
\int_\sigma\vr^{k,+}_{\vc n_{\sigma,K}}[\vu_\sigma\cdot\vc n_{\sigma,K}]^-{\rm d}S_x 
+\sum_{\sigma\in\partial\Omega_0\cap{\cal E}^{\rm in}}\int_\sigma
\vr_B |\vu_{B,\sigma}\cdot\vc n_\sigma|{\rm d}S_x,
$$
where the left hand side is strictly negative while the right hand side is positive. Consequently $\Omega_0=\emptyset$.

\subsection{Existence of solutions to the numerical scheme}\label{Ex}

Our goal is to show the following: Given $0<\vr^{k-1}\in Q(\Omega)$, $\vv^{k-1}\in V_0(\Omega;R^3)$
(and the boundary data $\vr_B,\vu_B$), the algebraic system
(\ref{N2}-\ref{N3}) admits a solution $0<\vr^{k}\in Q(\Omega)$ and $\vv^{k}\in V_0(\Omega;R^3)$. 

We first prove the solvability of the linear equation (\ref{N2}) with $\vu^k=\vv^{k}+{\vu}_B$, $\vv^k\in V_0(\Omega;R^3)$ given: It reads
$$
\intO{\vr^k\phi_K} + \Delta t\sum_{K\in {\cal T}}\sum_{\sigma \in {\cal E}(K)\cap{\cal E}^{\rm int}}
\int_\sigma F_{\sigma,K}(\vr^k,\vu^k)\phi_K{\rm d}S_x
$$
$$
+
\Delta t\sum_{\sigma\in{\cal E}^{\rm out}}\int_\sigma \vr^{k}\vu_{B,\sigma}\cdot\vc n_\sigma
\phi_K{\rm d}S_x
{ + \Delta t\kappa h^\omega\sum_{\sigma\in{\cal E}_{\rm int}}\int_\sigma[[\vr^k]]_{\sigma}[[\phi_K]]_\sigma{\rm d}S_x}=F_K,
$$
where $\{\phi_K\}_{K\in {\cal T}}$ is a basis in $Q(\Omega)$ and 
$$
F_K=\intO{\vr^{k-1}\phi_K}- \Delta t\sum_{\sigma\in{\cal E}^{\rm in}}\int_\sigma \vr_B\vu_{B,\sigma}\cdot\vc n_\sigma
\phi_K{\rm d}S_x.
$$
It admits a unique solutions $\vr^k\in Q(\Omega)$ -- we denote it by $\eta(\vv^k)$-- since the corresponding homogenous system
$$
\intO{\vr^k\phi_K} + \Delta t\sum_{K\in {\cal T}}\sum_{\sigma \in {\cal E}(K)\cap{\cal E}^{\rm int}}
\int_\sigma F_{\sigma,K}(\vr^k,\vu^k)\phi_K{\rm d}S_x
$$
$$
+
\Delta t\sum_{\sigma\in{\cal E}^{\rm out}}\int_\sigma \vr^{k}\vu_{B,\sigma}\cdot\vc n_\sigma
\phi_K{\rm d}S_x
{ + \Delta t\kappa h^\omega\sum_{\sigma\in{\cal E}_{\rm int}}\int_\sigma[[\vr^k]]_{\sigma}[[\phi_K]]_\sigma{\rm d}S_x}=0
$$
admits a unique solution $\vr^k=0$ as we easily show by applying to it the reasoning of Section \ref{PD} with the
convex function $B(s)=|s|$. We already proved in Section \ref{PD} that $\eta(\vv^k)>0$ and it is easy to see that 
$\intO{\eta(\vv^k)}$ is bounded and that the map $\eta$ is continuous {from $V_0(\Omega;R^3)$ to $Q(\Omega)$.}

We now define a map $\Theta: V_0(\Omega;R^3)\to V_0(\Omega; R^3)$ by saying that for $\vc z\in V_0(\Omega;R^3)$,
$\vv^k=\Theta(\vc z)$ is a solution of the problem
\bFormula{1}
-\intO{\mathbb{S}(\nabla_h\vv^k):\Gradh \phi} =<{\cal F}(\vc z),\phi>+<{\cal S},\phi> \;\mbox{for any $\phi \in V_{0}(\Omega;R^3)$},
\eF
where
$$
<{\cal S},\phi>=\intO{\mathbb{S}(\nabla_h\vu_B):\Gradh \phi},
$$
and
$$
<{\cal F}(\vc z),\phi>=\intO{ D_t (\eta(\vc z) \widehat{\vv}^k) \cdot \widehat\phi }
+\sum_{K\in {\cal T}}\sum_{\sigma \in {\cal E}(K)\cap{\cal E}_{\rm int}}\int_\sigma F_{\sigma,K}(\eta(\vc z)\widehat{\vv}^k,\vu^k)\cdot\widehat\phi{\rm d}S_x
$$
$$
{ + \sum_{K\in {\cal T}}\sum_{\sigma \in {\cal E}(K)}\int_\sigma\eta(\vc z)\widehat\vu^k\cdot{\vc n}_{\sigma,K}\vu_B\cdot\widehat\phi
{\rm d} S_x} 
+
\sum_{\sigma\in{\cal E}^{\rm out}}\int_\sigma \eta(\vc z)\vu_{B,\sigma}\cdot\vc n_\sigma\widehat\vv^k\cdot\widehat\phi{\rm d}S_x
+
\sum_{\sigma\in {\cal E}^{\rm in}}\int_\sigma \vr_B\vu_{B,\sigma}\cdot\vc n_\sigma\widehat\vv^k\cdot\widehat\phi{\rm d}S_x
$$
$$
 -\intO{ p_h(\eta(\vc z)) \Divh \phi }
 + \kappa h^\omega \sum_{\sigma\in{\cal E}_{\rm int}}\int_\sigma[[\eta(\vc z^k)]]_{\sigma}\{\widehat\vv^k\}_{\sigma}[[\widehat\phi]]_\sigma{\rm d}S_x.
$$
Clearly, $\vu^k=\vv^k+\vu_B$, $\vr^k=\eta(\vv^k)$, with $\vv^k$ any fixed point of the map $\Theta$
solves the algebraic problem (\ref{N1}--\ref{N3}).

We apply to the map $\Theta$  Schaeffer's fixed point theorem (see Lemma \ref{fpoint}) with $Z=V_0(\Omega;R^3)$. We have only to verify that
the set 
$$
\{\vc z\in Z\,|\, \vc z\;\mbox{solves (\ref{1}) with r.h.s. $\Lambda{\cal F}(\vc z) +\Lambda{\cal S}$, $\Lambda\in [0,1]$}\} 
$$
is non empty (indeed $0$ belongs to the set) and is bounded. The boundedness, follows from the energy balance. Indeed, the same calculation
--as we effectuate later in Section \ref{SEst} to  derive the energy estimates for the numerical scheme-- performed on the modified problem (\ref{1}) { with $\vv^k=\vc z$}, yields
(\ref{ebalance}) with $\vv^k$ replaced by $\vc z$, where all terms are multiplied by $\Lambda$ except the term 
$\Delta t\intO{\mathbb{S}(\nabla_h\vc z):\Gradh \vc z}$ at its left hand side. This readily implies the boundedness of the latter set by the  argumentation
in the proof of Lemma \ref{Lebalance}.

This completes the proof of Theorem \ref{TN1}.

\section{Uniform estimates}\label{SEst}
\subsection{Energy balance}

\begin{enumerate}
 \item We take in the continuity equation (\ref{N2}) the test function $\phi=1$ and obtain
 the conservation of mass (see equation (\ref{mbalance})).
\item We shall mimick the derivation of the energy iequality in the contiuous case
(cf. \cite{KwNo}) and compute 
\begin{equation}\label{eb1}
\mbox{
(\ref{N3})$_{\phi=\vc v^k}$ $+$ (\ref{N2})$_{\phi=
-\frac 12|\widehat{\vc v}^k|^2}$ $+$ (\ref{N2})$_{\phi=H_h'( \vr^k)}$,
}
\end{equation}
where
\bFormula{Hh}
H_h(\vr)=H(\vr)+\tilde\kappa h^\eta \vr^2.
\eF
Regrouping conveniently various terms, we get the following identities:
\begin{enumerate}
 \item Contribution of the terms with the time derivatives to (\ref{eb1}):
 $$
 \intO{\Big( D_t [\vr^k \vv^k]\cdot\widehat\vv^k -\frac 12 D_t [\vr^k]|\widehat\vv^k|^2+D_t [\vr^k] H_h'(\vr^k) \Big)}
 =\intO{ D_t \Big[\frac 12\vr^k |\widehat\vv^k|^2 + H_h(\vr^k)\Big] }
 $$
 $$
 +
 {\Delta t\intO{ \frac 12\vr^{k-1} |D_t[\widehat\vv^k]|^2}+\frac 1{\Delta t}\intO{\Big(H_h(\vr^{k-1})-H_h'(\vr^k)(\vr^{k-1}-\vr^k) -H_h(\vr^k)\Big)}.
 }
$$
 \item Contribution of the convective terms to (\ref{N3})$_{\phi=\vc v^k}$ $+$ (\ref{N2})$_{\phi=
-\frac 12|\widehat{\vc v}^k|^2}$:
 $$
 \sum_{K\in {\cal T}}\sum_{\sigma \in {\cal E}(K)\cap{\cal E}_{\rm int}}\int_\sigma F_{\sigma,K}(\vr^k\widehat{\vv}^k,\vu^k)\cdot\widehat\vv^k{\rm d}S_x
-
 \frac 12\sum_{K\in {\cal T}}\sum_{\sigma \in {\cal E}(K)\cap{\cal E}_{\rm int}}\int_\sigma F(\vr^k,\vu^k)|\widehat\vv^k|^2{\rm d}S_x
 $$
 $$
=\frac 12\sum_{\sigma\in {\cal E}_{\rm int}}\int_\sigma |{\rm Up}_\sigma(\vr^k,\vu^k)|\,\ju{\widehat\vv^k}_\sigma^2.
$$
\item { Contribution of the artificial density dissipation to (\ref{N3})$_{\phi=\vc v^k}$ $+$ (\ref{N2})$_{\phi=
-\frac 12|\widehat{\vc v}^k|^2}$:
$$
\kappa h^\omega\sum_{\sigma\in {\cal E}_{\rm int}}\int_{\sigma}\Big( [[\vr^k]]_\sigma\{\widehat\vv^k\}_\sigma\cdot[[\widehat\vv^k]]_\sigma
-\frac 12 [[\vr^k]]_\sigma[[|\widehat\vv^k|^2]]_\sigma\Big){\rm d}S_x=0.
$$
}
\item Contribution of the pressure term in the momentum and of the convective term in the continuity equation
to (\ref{N3})$_{\phi=\vc v^k}$ $+$ (\ref{N2})$_{\phi=
H'(\vr^k)}$:
$$
-
\intO{ p_h(\vr^k) \Divh \vv^k } + \sum_{K\in{\cal T}}\sum_{\sigma \in {\cal E}(K)\cap{\cal E}_{\rm int}}\int_\sigma F_{\sigma,K}(\vr^k,\vu^k)H_h'(\vr^k){\rm d}S_x
$$
$$
+\sum_{\sigma\in {\cal E}^{\rm out}}\int_{\sigma}
\vr^k\vu_{B,\sigma}\cdot\vc n_{\sigma}H_h'(\vr^k){\rm d}S_x
+\sum_{\sigma\in {\cal E}^{\rm in}}\int_{\sigma}
\vr_{B,\sigma}\vu_{B,\sigma}\cdot\vc n_{\sigma}H_h'(\vr^k){\rm d}S_x
$$
$$
= \sum_{K\in{\cal T}}\sum_{\sigma \in {\cal E}(K)} \int_\sigma p_h(\vr^k)\vu_{B,\sigma}\cdot
\vc n_{\sigma,K}{\rm d}S_x
$$
$$
+\sum_{\sigma \in {\cal E}_{\rm int}}
\int_\sigma\Big(H_h(\vr^{k,-})-H_h'(\vr^{k,+})(\vr^{k,-}-\vr^{k,+})- H_h(\vr^{k,+})\Big)|\vu^k_\sigma\cdot\vc n_{\sigma}|{\rm d}S_x
$$
$$
+\sum_{\sigma \in {\cal E}^{\rm out}}\int_\sigma H_h(\vr^k)\vu_{B,\sigma}\cdot\vc n_\sigma
{\rm d}S_x + \sum_{\sigma \in {\cal E}^{\rm in}}\int_\sigma H_h(\vr_B)\vu_{B,\sigma}\cdot\vc n_\sigma
{\rm d}S_x 
$$
$$
+ \sum_{\sigma \in {\cal E}^{\rm in}}\int_\sigma \Big(H_h(\vr_B)-H'(\vr^k)(\vr_B-\vr^k)
-H_h(\vr^k)\Big)
|\vu_{B,\sigma}\cdot\vc n_\sigma|
{\rm d}S_x.
$$
\item {
 Contribution of the artificial viscosity term in the continuity equation to
to (\ref{N2})$_{\phi=
H'(\vr^k)}$ is
$$
\kappa h^\omega\sum_{\sigma\in {\cal E}_{\rm int}}\int_\sigma[[\vr^k]]_\sigma[[H_h'(\vr^k)]]_\sigma{\rm d} S_x.
$$
}
\end{enumerate}
\end{enumerate}
Putting together items 2(a)--2(e) while taking into account the evident contribution of the dissipative terms for the velocity
in the momentum equation to  (\ref{N3})$_{\phi=
\vv_k}$,
we get the energy balance. We gather the above calculations in the following lemma.
\bLemma{ebalance}{\rm [Mass conservation and energy balance]}
Suppose that the pressure satisfies assumptions (\ref{pres1}). Then any solution of $(\vr,\vu)$ of the
algebraic system (\ref{N1}--\ref{N3}) satisfies for all $m=1,\ldots,N$ the following:
\begin{enumerate}
 \item Mass conservation
 \bFormula{mbalance}
 \vr^m>0,\;
 \intO{\vr^m}+\Delta t\sum_{k=1}^m\sum_{\sigma\in{\cal E}^{\rm out}}\int_\sigma\vr^k\vu_{B,\sigma}\cdot\vc n_\sigma{\rm d}S_x=
 \intO{\vr^0}-\Delta t\sum_{k=1}^m\sum_{\sigma\in{\cal E}^{\rm in}}\int_\sigma\vr_B\vu_{B,\sigma}\cdot\vc n_\sigma{\rm d}S_x;
 \eF
\item Energy balance. There exist $\overline\vr^{k-1,k}\in Q(\Omega)$, $\overline\vr_K^{k-1,k}\in [\min\{\vr_K^{k-1},\vr_K^k\},
\max\{\vr_K^{k-1},\vr_K^k\}]$, $K\in{\cal T}$ and $\overline\vr^{k,\sigma}\in [\min\{\vr_\sigma^{-},\vr_\sigma^+\}, 
\max\{\vr_\sigma^{-},\vr_\sigma^+\}]$, $\sigma\in{\cal E}_{\rm int}$, $k=1,\ldots,N$, such that
\bFormula{ebalance}
\intO{\frac 12\vr^m |\widehat\vv^m|^2} + \intO{H_h(\vr^m) }
+\Delta t  \sum_{k=1}^m\intO{\mathbb{S}(\nabla_h\vv^k):\nabla_h\vv^k}
\eF
 $$
 +
 \frac 12\sum_{k=1}^m\intO{\Big( \vr^{k-1} |\widehat\vv^k-\widehat\vv^{k-1}|^2+H_h''(\overline\vr^{k-1,k})
  |\vr^k-\vr^{k-1}|^2\Big)
 }
 $$
 $$
 +\frac {\Delta t}2\sum_{k=1}^m\sum_{\sigma\in {\cal E}_{\rm int}}\int_\sigma |{\rm Up}_\sigma(\vr^k,\vu^k)|\,\ju{\widehat\vv^k}_\sigma^2{\rm d}S_x 
+\kappa h^\omega { {\Delta t} \sum_{k=1}^m}\sum_{\sigma\in {\cal E}_{\rm int}}\int_\sigma[[\vr^k]]_\sigma[[H_h'(\vr^k)]]_\sigma{\rm d} S_x
$$
 $$
+\frac {\Delta t} 2\sum_{k=1}^m\sum_{\sigma \in {\cal E}_{\rm int}}
\int_\sigma H_h''(\overline\vr^{k,\sigma})\ju{\vr^k}_\sigma^2|\vu^k_\sigma\cdot\vc n_{\sigma}|{\rm d}S_x
 +
 {\Delta t}\sum_{k=1}^m\sum_{\sigma\in {\cal E}^{\rm out}}\int_\sigma \vr^k
 \vu^k_\sigma\cdot\vc n_\sigma|\widehat\vv^k|^2{\rm d}S_x
$$
$$
 + {\Delta t}\sum_{k=1}^m\sum_{\sigma \in {\cal E}^{\rm out}}\int_\sigma H_h(\vr^k)\vu_{B,\sigma}\cdot\vc n_\sigma
{\rm d}S_x  
+ \Delta t\sum_{k=1}^m\sum_{\sigma \in {\cal E}^{\rm in}} E_{H_h}(\vr_B|\vr^k)
|\vu_{B,\sigma}\cdot\vc n_\sigma|
{\rm d}S_x
$$
$$
 =\intO{\Big[\frac 12\vr^0 |\widehat\vv^0|^2 + H_h(\vr^0)\Big] } 
 +
 \Delta t\sum_{k=1}^m\sum_{\sigma \in {\cal E}^{\rm in}}\int_\sigma H_h(\vr_B)|\vu_{B,\sigma}\cdot\vc n_\sigma|
{\rm d}S_x
$$
$$
-\Delta t  \sum_{k=1}^m\intO{\mathbb{S}(\nabla_h\vu_B):\nabla_h\vv^k}
 -\Delta t\sum_{k=1}^m
  \intO{ p_h(\vr^k){\rm div}_h\vu_{B}}
 $$
 $$
{
 -\Delta t\sum_{k=1}^m\intO{\vr^k\widehat\vu^k\cdot\nabla_h\vu_B\cdot\widehat\vv^k}
 - {\Delta t}\sum_{k=1}^m\sum_{\sigma\in {\cal E}^{\rm in}}\int_\sigma \vr_B\vu_{B,\sigma}\cdot\vc n_\sigma |\widehat\vv^{k}|^2{\rm d}S_x}.
$$
 \end{enumerate}
 \eL
 
We are now in position to deduce from the energy balance (\ref{ebalance}) the uniform estimates.
This will be done by using the Gronwall lemma. Before starting, it would be convenient to 
mention some properties of the pressure $p$ and of the Helmholtz function $H$ which can be deduced from assumptions 
(\ref{pres1}--\ref{pres4}).
\bRemark{pres}
\begin{enumerate}
\item We know that $H$ is strictly convex (since $p'(\vr)>0$). Consequently $p$ is strictly convex as well (since
$\overline a p-H$ is convex). We may thus suppose, without loss of generality, that both $H-\underline a  p$ and
$\overline a p-H$ are strictly convex functions.
\item It is easy to deduce from the convexity of function $\overline a p-H$  that
\bFormula{HPlow}
\exists\overline R>0,\;\forall\vr>\overline R,\;  \vr^{\gamma-1}\aleq H'(\vr) ,\;  \vr^{\gamma-1}\aleq p'(\vr)\;\mbox{where $\gamma=1+ 1/\overline a$},
\eF
see \cite[formula (2.3)]{AbFeNo}. This observation explains the value of $\gamma$ in Theorem \ref{TN2}.
\item Convexity of $\underline a p -H$ and $p(0)=H(0)=0$ yield
\bFormula{condws}
\forall \vr \in (0,\infty),\; p(\vr)\aleq \frac 1{\underline a}H(\vr)+ c\vr,\; c=\frac{\underline a p'(0)-H'(0)}{\underline a}.
\eF
This observation allows to prove later Lemma \ref{Lrelaxed2} and consequently the error estimates.
\item Since $H-\underline a p$ is convex, there are numbers $d\ge 0$, $e\ge 0$ such that
\bFormula{calH*}
\inf_{\vr >0}({\cal H}(\vr)-\underline a p(\vr))\ge 0,\;\lim_{\vr\to\infty} ({\cal H}(\vr)-\underline a p(\vr))=\infty,\;\mbox{where
${\cal H}(\vr)=H(\vr) +d\vr +e$}. 
\eF
This observation allows to derive the uniform estimates from the energy balance (\ref{ebalance}).
\end{enumerate}
\eR
With this observation at hand, for the purpose of obtaining the uniform estimates,
we can replace in the energy balance at its left hand side
$$
 \intO{H_h(\vr^m) }\;\mbox{by}\; \intO{{\cal H}_h(\vr^m) },\; {\cal H}_h(\vr)={\cal H}(\vr) +\tilde\kappa h^\eta\vr^2
$$
and add  the non negative  term
$$
{\Delta t}\sum_{k=1}^m\sum_{\sigma \in {\cal E}^{\rm out}}\int_\sigma d\vr^k\vu_{B,\sigma}\cdot\vc n_\sigma
{\rm d}S_x.
$$
Likewise we can replace at its right hand side
$$
 \intO{H_h(\vr^0) }\;\mbox{by}\; \intO{{\cal H}_h(\vr^0) }
$$
and add the term
$$ 
-{\Delta t}\sum_{k=1}^m\sum_{\sigma \in {\cal E}^{\rm in}}\int_\sigma d\vr_B\vu_{B,\sigma}\cdot\vc n_\sigma
{\rm d}S_x.
$$
Due to the balance of mass (\ref{mbalance}) the energy balance modified in this way is again satisfied as an identity.

Now, we shall write the term $\intO{{\cal H}_h(\vr^m) }$ at the left hand side of the modified energy balance as a sum
$$
\intO{({\cal H}_h(\vr^m)-\underline{a} p(\vr^m ))} + \underline a \intO{p(\vr^m) }
$$
In view of this, in order to deduce from the modified energy balance  the uniform estimates, we want to find
 a bound from above of its right hand side by the expression
$$
\frac 1\delta\int_0^{\tau_m}\intO{\Big(\vr|\widehat \vv|^2 + p(\vr) +\tilde\kappa h^\eta\vr^2\Big)}{\rm d}t+
\delta \int_0^{\tau_m}\intO{|\nabla_h\vv|^2}{\rm d}t +C, \; \tau_m=m\Delta t
$$
with some $C>0$ depending only of the data,  $T$ and $\Omega$ and on $\delta>0$,
with the goal to employ the Gronwall inequality.

 The first line contains only the data. The third term can be treated by the Cauchy-Schwartz and Young inequalities. The fourth term is easy as well.
For the fifth term, we have,
$$
\int_0^{\tau_m}\intO{\vr\widehat\vu\cdot\nabla_h\vu_B\cdot\widehat\vv}{\rm d}t
\aleq \tau_m+\int_0^{\tau_m}\intO{\vr\widehat\vv^2}{\rm d}t.
$$
  For the  last term, we write
\bFormula{enb1}
 |\sum_{\sigma\in {\cal E}^{\rm in}}\int_\sigma \vr_B\vu_{B,\sigma}\cdot\vc n_\sigma |\widehat\vv|^2{\rm d}S_x|\aleq
\sum_{\sigma\in {\cal E}^{\rm in}}\|\widehat \vv\|_{L^2(\sigma)}^2
\eF
$$
\aleq
h^{-1}\sum_{K\in {\cal T},\,K\cap{\cal E}^{\rm in}\neq\emptyset}\|\vv\|^2_{L^2(K)}
\aleq h \sum_{K\in {\cal T},\,K\cap{\cal E}^{\rm in}\neq\emptyset}\|\Grad\vv\|^2_{L^2(K)}\aleq h \|\nabla_h\vv\|^2_{L^2(\Omega)}.
$$
Indeed the first inequality is the Cauchy-Schwartz inequality for the integrals over $\sigma$. The second inequality uses the discrete trace estimate
(\ref{traceFE}) and the Jensen's inequality (\ref{JensenV}). Finally, we use one of the inequalities
(\ref{PoincareV-}) in combination with the fact that $\vv_\sigma=0$ for any $\sigma\in {\cal E}^{\rm in}$.

\subsection{Estimates}
In the following lemma we shall gather all estimates which can be deduced
from Lemma \ref{Lebalance}.
 
\bLemma{estimates}
We denote $I=[0,T]$,
$$
\overline E_0=\sup_{h\in (0,1)}E_{0,h},\;\; E_{0,h}= \intO{ \left[ \frac{1}{2} \vr^0_h | \widehat{\vu}^0_h |^2 +  H_h(\vr^0_h)\right] }.
$$
We suppose that the pressure satisfies assumption (\ref{pres1}--\ref{pres4}) and that $h=\Delta t$.

Then there exists a number $\mathfrak{d}$, $
\mathfrak{d}:= \mathfrak{d}(\|\mathfrak{u}_B\|_{C^1(\overline\Omega)},\|\mathfrak{r}_B\|_{C(\overline\Omega)}, \overline E_{0}, T,\Omega)\aleq 1$,
such that any solution $(\vr,\vu)=(\vr_{h,\Delta t}, \vu_{h,\Delta t})$ of the discrete problem (\ref{N1}--\ref{N3}) admits the following estimates:
\bFormula{e0}
\vr>0, \|\vr\|_{L^\infty(I;L^\gamma(\Omega))} \aleq 1 ,\;\|p(\vr)\|_{L^\infty(0,T;L^1(\Omega))}\aleq 1,\; \|H(\vr)\|_{L^\infty(0,T;L^1(\Omega))}\aleq 1,
\eF
\bFormula{e1}
\|\vr|\mathfrak{u}_{B}\cdot\vc n|^{1/\gamma}\|_{L^\gamma(I;L^\gamma(\partial\Omega))}\aleq 1,\;{
\|H(\vr)|\mathfrak{u}_{B}\cdot\vc n|\|_{L^1(I;L^1(\partial\Omega))}\aleq 1,}\footnote{Here we use the definition of $\Pi^V$
and the property (\ref{fit}) of the mesh.}
\eF
\bFormula{e0+}
\tilde\kappa h^{\eta/2} \|\vr\|_{L^\infty(I;L^2(\Omega))} \aleq 1 , 
\eF
\bFormula{e1+}
\tilde\kappa h^{\eta/2} \|\vr |\mathfrak{u}_{B}\cdot\vc n|^{1/2}\|_{L^2(I;L^2(\partial\Omega))}\aleq 1,
\eF
\bFormula{e2}
{\rm sup}_{\tau \in (0,T)} \| \sqrt{ \vr } \widehat{\vu} (\tau, \cdot) \|_{L^\infty(I;L^2(\Omega)} \aleq 1,
\eF
\bFormula{e3}
\|\nabla_h \vv\|_{L^2(I\times\Omega)} \aleq 1,\quad \|\vv\|_{L^2(I; L^6(\Omega))}\aleq 1, \footnote{The latter estimate follows from the first one by the discrete Sobolev inequality (\ref{SobolevV0}).}
\eF
\bFormula{e4}
\tilde\kappa h^\eta\int_0^T\intO{\Delta t |D_t\vr|^2}{\rm d}t=
\tilde\kappa h^\eta\sum_{k=1}^m \intO{   
\left| \vr^k - \vr^{k-1}  \right|^2  } \aleq 1,
\eF
\bFormula{e4+}
\int_0^T\intO{\Delta t \vr(t-\Delta t) |D_t\vv (t)|^2}{\rm d}t=
\sum_{k=1}^m \intO{    \vr^{k-1}\left| \widehat{\vv}^k - \widehat{\vv}^{k-1}  \right|^2 } \aleq 1,
\eF
\bFormula{e5}
\sum_{\sigma\in {\cal E}_{\rm int}}\int_0^T\int_\sigma |{\rm Up}_\sigma(\vr,\vu)
| \, [[\vv]]^2_\sigma{\rm d}S_x{\rm d}t \aleq 1,
\eF
\bFormula{e6}
 \tilde\kappa h^\eta\sum_{\sigma\in{\cal E}_{\rm int}} \int_0^T\int_\sigma \ju{ \vr  }_\sigma^2 | \vu_\sigma \cdot \vc{n}| {\rm dS}_x  \dt \aleq 1,
\eF
\bFormula{e7}
\kappa\tilde\kappa h^{\omega+\eta}\sum_{\sigma\in{\cal E}_{\rm int}}\int_0^T\int_\sigma  \ju{ \vr  }_\sigma^2  {\rm dS}_x  \dt \aleq 1.
\eF
\eL

\bRemark{rest}
Some other estimates can be deduced if we consider particular situations enumerated in
Remark \ref{Rtr1}. They are essential in the proof of the versions of Theorem \ref{TN2} in those particular situations.

Suppose that the pressure satisfies in addition to (\ref{pres1}--\ref{pres4}) also (\ref{pres2}). Then we have:
\begin{enumerate}
\item If in (\ref{pres2}), $\underline\pi>0$ and $\gamma\ge 2$ we have in addition to (\ref{e6}),
\bFormula{e6*}
 \sum_{\sigma\in{\cal E}_{\rm int}} \int_0^T \int_{\sigma} 
 \ju{ \vr  }_\sigma^2 | \vu_\sigma \cdot \vc{n}| {\rm dS}_x  \dt \aleq 1,
\eF
and, in addition to (\ref{e4})
\bFormula{e4*}
\int_0^T\intO{\Delta t |D_t\vr|^2}{\rm d}t=
\sum_{k=1}^m \intO{   
\left| \vr^k - \vr^{k-1}  \right|^2  } \aleq 1.
\eF

\item If  $1<\gamma\le 2$, we have in addition to (\ref{e6})
\bFormula{e6**}
 \sum_{\sigma\in{\cal E}_{\rm int}} \int_0^T \int_{\sigma} 
 \Big(\frac 1{\{\vr\}_\sigma +1}\Big)^{2-\gamma}\ju{ \vr  }_\sigma^2 | \vu_\sigma \cdot \vc{n}| {\rm dS}_x  \dt \aleq 1,
\eF
and, in addition to (\ref{e4}),
\bFormula{e4**}
\sum_{k=1}^m \intO{    \Big(\frac 1{\vr^{k-1}+\vr^{k} +1}\Big)^{2-\gamma}
\left| \vr^k - \vr^{k-1}  \right|^2  } \aleq 1.
\eF
\item If $1<\gamma<2$ and condition (\ref{pres3}) is satisfied, we have also
\bFormula{e8*}
\kappa h^\omega\sum_{\sigma\in {\cal E}_{\rm int}}\int_0^T \int_{\sigma} 
 \ju{ \vr^{\gamma/2}  }_\sigma^2  {\rm dS}_x  \dt \aleq 1,
\eF
and consequently
\bFormula{e8**}
\kappa h^{1+\omega}\int_0^T\|\vr\|_{L^{3\gamma}(\Omega)}^\gamma{\rm d}t\aleq 1. \footnote{This requires some explanation:
We deduce from (\ref{pres3}) that $[[(\vr^k)^{\gamma/2}]]^2_\sigma\aleq H''(\overline\vr^{k,\sigma})(\vr^{k,-}-\vr^{k,+})^2$ and thus
the fourth line of (\ref{ebalance}) yields, in particular, (\ref{e8*}).
The latter implies (\ref{e8**}) by virtue of Sobolev inequality (\ref{SobolevQ}) for the broken norms.}
\eF 
\end{enumerate}
\eR

\section{Consistency}
\label{NC}

Having collected all the available uniform bounds, our next task is to verify that our numerical method is \emph{consistent} with the variational formulation of the original problem.

\subsection{Continuity equation}\label{CE}

For $\phi \in C^1([0,T]\times\Ov{\Omega})$, take $\widehat\phi(t)$ as a test function in the discrete continuity equation (\ref{N2}). We see that
$$
\intO{D_t\vr\widehat\phi}= \intO{D_t\vr\phi}.
$$
Further, using the formula (\ref{Up3}), and noticing that
$$
\intO{ (\widehat\phi  - \phi ) \vr \Divh \vu } = \sum_{K\in {\cal T}} \int_K{ 
(\widehat\phi- \phi ) \vr\Div \vu } = 0
$$
(indeed $\Divh \vu(t)$ is constant on each element $K$), 
we check without difficulty using formula (\ref{Up3}) that
\[ 
\intO{  \vr \vu \cdot \Grad \phi  }= \sum_{K\in {\cal T}} 
\sum_{\sigma\in {\cal E}(K)}\int_{\sigma} (\phi-\widehat\phi) \vr(\vu - \vu_\sigma ) \cdot \vc{n}_{\sigma,K} \ {\rm dS}_x +\intO{(\widehat\phi-\phi)\vr\Divh\vu}
\]
\[
 - \sum_{K\in {\cal T}}\sum_{\sigma \in {\cal E}(K)\cap{\cal E}_{\rm int} } \int_\sigma F_{\sigma,K}(\vr, \vu)\ \widehat\phi  {\rm d}S_x + 
\sum_{K\in {\cal T}}\sum_{\sigma\in {\cal E}(K)\cap{\cal E}_{\rm int} }\int_\sigma ( \widehat\phi - \phi ) \ju{\vr}_{{\vc n}_{\sigma,K}}  \ [\vu_\sigma \cdot \vc{n}_{\sigma,K}]^- {\rm d}S_x  
\]
$$
+  \sum_{\sigma\in {\cal E}^{\rm out}}\int_\sigma
\vr\vu_{B,\sigma}\cdot \vc n_{\sigma}(\phi-\widehat\phi){\rm d}S_x
+  \sum_{\sigma\in {\cal E}^{\rm in}}\int_\sigma
\vr\vu_{B,\sigma}\cdot \vc n_{\sigma}(\phi-\widehat\phi){\rm d}S_x.
$$

Consequently, equation (\ref{N2}) rewrites
\bFormula{ce-v}
\intO{\Big(D_t\vr\phi-\vr\widehat\vu\cdot\Grad\phi\Big)} +  \int_{\Gamma^{\rm out}}
\vr\mathfrak{u}_{B}\cdot \vc n \phi{\rm d}S_x
+  \int_{\Gamma^{\rm in}}
\vr_B\mathfrak{u}_{B}\cdot \vc n\phi{\rm d}S_x=<{ R}^C_h(t),\phi> 
\eF
in $(0,T]$, where
$$
<{R}^C_h(t),\phi>=
 \sum_{K\in {\cal T}}\sum_{\sigma\in {\cal E}(K)\cap{\cal E}_{\rm int} }\int_\sigma ( \phi - \widehat\phi ) \ju{\vr}_{{\vc n}_{\sigma,K}}  \ [\vu_\sigma \cdot \vc{n}_{\sigma,K}]^- {\rm d}S_x  
$$
$$
+\sum_{K\in {\cal T}} 
\sum_{\sigma\in {\cal E}(K)}\int_{\sigma} (\widehat\phi-\phi) \vr(\vu - \vu_\sigma ) \cdot \vc{n}_{\sigma,K} \ {\rm dS}_x 
{ +\intO{\vr{\rm div}_h\vu(\phi-\widehat\phi)}}
$$
$$
{ -\kappa h^\omega\sum_{\sigma\in {\cal E}_{\rm int}}\int_\sigma [[\vr]]_\sigma [[\widehat\phi]]_\sigma {\rm d} S_x}
+ \sum_{\sigma\in {\cal E}^{\rm in}}\int_\sigma
(\vr-\vr_B)\vu_{B,\sigma}\cdot \vc n_{\sigma}(\widehat\phi-\phi){\rm d}S_x
{ + \intO{\vr(\vu-\widehat\vu)\cdot\Grad\phi}}.
$$
In the above, we have also used the fact that the mesh fits to the inflow-outflow boundary and the definition (\ref{PiV}) of the projection $\Pi^V$.

{ Alternatively, in view of (\ref{pwt+}), we can integrate the equation (\ref{ce-v}) over time  and rewrite it in the form}
\bFormula{ce-v+}
\intO{\utilde{\vr}\phi(\tau)}-\intO{{\vr}^0\phi(0)}-\int_{0}^\tau\intO{\Big(\utilde{\vr}\partial_t\phi+\vr\widehat\vu\cdot\Grad\phi\Big)}{\rm d}t
\eF
$$
+  \int_{0}^{\tau}\int_{\Gamma^{\rm out}}
\vr\mathfrak{u}_{B}\cdot \vc n\phi{\rm d}S_x{\rm d}t
+  \int_{0}^{\tau}\int_{\Gamma^{\rm in}}
\vr_B\mathfrak{u}_{B}\cdot \vc n\phi{\rm d}S_x{\rm d}t
=\int_{0}^{\tau}< R^C_{h}(t),\phi>{\rm d}t 
$$
with any $\tau\in(0,T]$ and all $\phi\in C^1([0,T]\times\overline\Omega)$,
{ where $\utilde \vr$ is defined in (\ref{pwt+}).} 


Our next goal is to estimate conveniently all terms in the remainder $<R^C_{h,\Delta t},\phi>$. To do this we shall use the tools evoked in Section \ref{Preln} and employ the bounds (\ref{e0}--\ref{e7}). 
\begin{enumerate}
\item 
 {\it The first term in $<{ R}^C_{h},\phi>$} is bounded from above by
\bFormula{cest1}
\aleq \Big[h^\eta\sum_{\sigma\in{\cal E}_{\rm int}}  \int_{\sigma} 
 \ju{ \vr  }_\sigma^2 | \vu_\sigma \cdot \vc{n}_\sigma| {\rm dS}_x \Big]^{1/2}
\times\Big[h^{-\eta}\sum_{K\in {\cal T}} 
\sum_{\sigma\in {\cal E}(K)\cap{\cal E}_{\rm int}}\int_\sigma  ( \widehat\phi - \phi )^2  |\vu_\sigma \cdot \vc{n}_{\sigma}| {\rm d}S_x\Big]^{1/2}, 
\eF
where the first term in the product is controlled by (\ref{e6}). 
As far as for the second term in the  product, we have for any $\gamma>1$:
$$
h^{-\eta}\Big|\sum_{K\in {\cal T}} 
\sum_{\sigma\in {\cal E}(K)\cap{\cal E}_{\rm int}}\int_\sigma  ( \widehat\phi - \phi )^2  |\vu_\sigma \cdot \vc{n}_{\sigma}| {\rm d}S_x\Big|
{ \aleq
h^{-\eta}\sum_{K\in {\cal T}} 
\sum_{\sigma\in {\cal E}(K)\cap{\cal E}_{\rm int}}\int_\sigma\|\widehat\phi-\phi\|_{L^\infty(\sigma)}^2 \|\vu\|_{L^{1}(\sigma)}
}
$$
\bFormula{cest2}
\aleq
h^{1-\eta}\sum_{K\in {\cal T}} \|\Grad\phi\|_{L^\infty(K)}^2 \|\vu\|_{L^{1}(K)}
\aleq h^{1-\eta} \|\Grad\phi\|^2_{L^\infty(\Omega)}\|\vu\|_{L^{1}(\Omega)},
\eF
where $\|\vu\|_{L^{1}(\Omega)}\aleq \|\vu\|_{L^{6}(\Omega)}$ is bounded
in $L^2((0,T))$. Indeed, to get the second line we have employed the H\"older and Jensen inequalities on $\sigma$, cf. (\ref{JensenV}), third line employs the trace estimates (\ref{trace}--
\ref{traceFE}) and then  the H\"older inequality. 

\item By the same token, the absolute value of {\it the second term in $<{R}^C_{h},\phi>$}  
$$
\sum_{K\in {\cal T}} 
\sum_{\sigma\in {\cal E}(K)}\int_{\sigma} (\widehat\phi-\phi) \vr(\vu - \vu_\sigma ) \cdot \vc{n}_{\sigma,K} \ {\rm dS}_x
$$
allows for any $\gamma>1$ the following estimate
\bFormula{cest3}
\aleq  h^{1-\frac\eta 2}\|h^{\eta/2}\vr\|_{L^2(\Omega)}\|\nabla_h\vu\|_{L^{2}(\Omega)}\|\Grad\phi\|_{L^\infty(\Omega)}\aleq h^{1-\frac\eta 2} \|\nabla_h\vu\|_{L^{2}(\Omega)} \|\Grad\phi\|_{L^\infty(\Omega)},
\eF
where we have used (\ref{PoincareV-}). {{\it The third term and the last term in $<{R}^C_{h},\phi>$} admit the same estimate. }
\item {\it The artificial density diffusion term in $<{ R}^C_h,\phi>$} can be processed in the similar way,
\bFormula{cestad}
 \aleq \kappa h^{1+\omega}\|\Grad\phi\|_{L^\infty(\Omega)}\sum_{\sigma\in {\cal E}_{\rm int}}\|\vr\|_{L^\gamma(\sigma)}|\sigma|^{1/\gamma'} 
\aleq \kappa h^{\omega} \|\vr\|_{L^\gamma(\Omega))} \|\Grad\phi\|_{L^\infty(\Omega)}
\eF
for any value $\gamma>1$. 

\item {\it The boundary term in  $<{ R}^C_{h}(t),\phi>$} is controlled in the same way (without loss of generality, 
we perform the calculation just for the first of them, and with $\vr$ only -instead of with $(\vr-\vr_B)$):
\bFormula{cest4}
\aleq
h\|\nabla\phi\|_{L^{\infty}(\Omega)}\sum_{K\in {\cal T},\, {\cal E}(K)\cap {\cal E}^{\rm in}\neq\emptyset}
\sum_{\sigma\in {\cal E}^{\rm in}}\|\vr\|_{L^1(\sigma)}
\aleq \|\nabla\phi\|_{L^{\infty}(\Omega)} \|\vr\|_{L^2({\cal U})}|{\cal U}|^{1/2}\aleq h^{1/2-\eta/2} \|\nabla\phi\|_{L^{\infty}(\Omega)},
\eF
where ${\cal U}=\cup_{K\in {\cal T},\,{\cal E}(K)\cap {\cal E}^{\rm in}\neq\emptyset} K$.
Indeed, in the passage from the first to the second line, we have used trace estimates (\ref{traceFE}), and for the rest the H\"older inequalities 
as well as the fact that
the Lebesgue measure $|{\cal U}|$ of ${\cal U}$ is $\approx h$.


{
\item At one place of the convergence proof, we shall need a slightly improved version of the above estimates with test functions $\phi\in L^p(0,T; W_0^{1,p}(\Omega))$
with $p$ "large". To achieve this, we observe that the final bound in (\ref{cest2}) can be replaced by $\aleq h^{1-\eta} \|\Grad\phi(t)\|_{L^p(\Omega)}\|\vu(t)\|_{L^{6}(\Omega)}$ provided $p\ge 12/5$, final bound in (\ref{cest3}) by $ \aleq h^{\frac{3(p-2)}{2p}-\frac {1+\eta}2}
 \|h^{\eta/2}\vr\|_{L^2(\Omega)} \|\vu(t)\|_{L^{6}(\Omega)} \|\Grad\phi(t)\|_{L^p(\Omega)}$ provided ${\frac{3(p-2)}{2p}-\frac {1+\eta}2}>0$
--here we must use also "negative" interpolation estimate (\ref{interG})-- and
the final bound in (\ref{cestad}) by
$\aleq \kappa h^{\omega} \|\vr(t)\|_{L^\gamma(\Omega))} \|\Grad\phi(t)\|_{L^p(\Omega)}
$ provided $\frac 1\gamma+\frac 1p>1$.
}
\end{enumerate}

We resume the above calculations in the following lemma.
\bLemma{c-consistency}{\rm [Consistency for the continuity equation]}
Let $h=\Delta t$, $\tilde\kappa=1$ and $\eta\in (0,1).$ Let the pressure satisfy the hypotheses (\ref{pres1}--\ref{pres4}). 
\begin{enumerate}
\item
Then   any numerical solution 
of problem (\ref{N1}--\ref{N3}) 
satisfies the continuity equation in the variational form (\ref{ce-v+}) with any $\phi\in C^1([0,T]\times\overline\Omega)$,  where
\bFormula{NC1}
\left|\left<{ R^C_{h}}, \phi \right> \right| \aleq h^{\alpha_C} r^C_h \|\phi \|_{L^{\infty}(0,T; W^{1,\infty}(\Omega))}, \;
\|r^C_h\|_{L^{ 4/3}((0,T))}\aleq 1
\eF
with $\alpha_C=\frac{1-\eta} 2>0$ if $\kappa=0$ and 
$\alpha_C=\min\{\frac{1-\eta} 2,\omega\}>0$ if $\kappa=1$.
{
\item There is $p_0(\eta)>1$ such that for all $p>p_0(\eta)$
\bFormula{NC1+}
\int_0^\tau\left|\left<{ R^C_{h}}, \phi \right> \right|{\rm d}t \aleq h^{\underline\alpha_C}  \|\phi \|_{L^{p}(0,T; W^{1,p}(\Omega))},
\eF
with some $\underline\alpha_C= \underline\alpha_C(p)>0$ and all $\phi\in C^1_c((0,T)\times\Omega)$.
}
\end{enumerate}

\eL

\bRemark{Rc-consistency}
{ Let $\tilde\kappa=0$ and let $p$ satisfy in addition to (\ref{pres1}--\ref{pres4}) also condition (\ref{pres2}).}
\begin{enumerate}
\item  Formulas (\ref{NC1}--\ref{NC1+}) with $\alpha_C>0$ and $\underline\alpha_C>0$  remain still valid  under condition $6/5<\gamma<2$ provided also  (\ref{pres3}) is satisfied.  In this case
$$
\alpha_C=\min\Big\{\frac 12, 
\frac{5\gamma-6}{2\gamma}, \frac{1}{\gamma'}\Big\} \;\mbox{if $\kappa=0$}
$$
and it is minimum of the above value and $\omega$ if $\kappa=1$. { Also, there is $p_0=p_0(\gamma)>1$
such that (\ref{NC1+}) holds for all $p>p_0$ with some $\underline\alpha_C=\underline\alpha_C(p)>0$.}

Indeed, 
we can use estimates (\ref{e6**}) instead of (\ref{e6}) when calculating 
(\ref{cest1}) which now yields,
$$
\aleq\Big[\sum_{\sigma\in{\cal E}_{\rm int}}  \int_{\sigma} 
 \Big(\frac 1{\{\vr\}_\sigma +1}\Big)^{2-\gamma}\ju{ \vr  }_\sigma^2 | \vu_\sigma \cdot \vc{n}_\sigma| {\rm dS}_x \Big]^{1/2}
$$
$$
\times\Big[\sum_{K\in {\cal T}} 
\sum_{\sigma\in {\cal E}(K)\cap{\cal E}_{\rm int}}\int_\sigma \Big(\{\vr\}_\sigma +1\Big)^{2-\gamma} ( \widehat\phi - \phi )^2  |\vu_\sigma \cdot \vc{n}_{\sigma}| {\rm d}S_x\Big]^{1/2} 
$$
and the estimate { (we use systematically  the "negative" interpolation (\ref{interG})),}
$$
\aleq h^{\frac 12}\|\vr\|^{\frac{2-\gamma}2}_{L^\gamma(\Omega)} \|\nabla\phi\|_{L^{\infty}(\Omega)} \|\vu\|^{1/2}_{L^6(\Omega)} A_h\aleq h^{\frac 12} B_h \|\nabla\phi\|_{L^{\infty}(\Omega)},\;
\mbox{ $B_h$  bounded in $L^{4/3}((0,T))$}
$$ 
for the first term in ${ R}^C_h$. 

Further, we use
estimate 
$$
\|\vr\|_{L^2(\Omega)}\aleq h^{\min\{0,\frac 32-\frac 3\gamma\}}\|\vr\|_{L^\gamma(\Omega)},
$$
when evaluating term corresponding to (\ref{cest3}), which gives
the bound
$$
\aleq h^{1+\min\{0,\frac 32-\frac 3\gamma\}} \|\nabla\phi\|_{L^{\infty}(\Omega)} A_h,\;\mbox{$A_h$ bounded in $L^2((0,T))$}
$$ 
in these cases.\footnote{This is the place when $\gamma>6/5$ is needed, in order
to let the power of $h$ positive. Here the situation could be still handled for $1<\gamma\le 6/5$ provided $\kappa=1$ by using estimate (\ref{e8**}). Similar terms in
the remainder of the consistency formulation of the momentum equation, however, impose the restriction $\gamma>6/5$
even if $\kappa=1$. It is therefore useless to treat this situation here.} 

 Estimate (\ref{cestad}) of the artificial diffusion term remains in force, while the the estimate of the boundary term (\ref{cest4}) can be replaced by
$$
\aleq \|\nabla\phi\|_{L^{\infty}(\Omega)} \|\vr\|_{L^\gamma({\cal U})}|{\cal U}|^{1/\gamma'}{ \aleq h^{1/\gamma'}\|\nabla\phi\|_{L^{\infty}(\Omega)}}.
$$
  
\item Similar reasoning can be carried out if $\gamma\ge 2$ and if in (\ref{pres2}), $\underline\pi>0$. In this case, in (\ref{NC1}),
$$
\alpha_C=\min\Big\{\frac 12, \frac{1}{\gamma'}\Big\}\;\mbox{if $\kappa=0$}
$$
and it is minimum of the above value and $\omega$ if $\kappa=1$.
\end{enumerate}
\eR

\subsection{Momentum equation}

The next step is to take
\[
\tilde\phi(t), \ \phi \in C^1_c([0,T]\times{\Omega}; R^3),
\]
as a test function in the discrete momentum equation (\ref{N3}). Seeing that, in accordance with (\ref{vv1}), (\ref{proj}),
\[
\intO{\mathbb{S}(\nabla_h\vv):\nabla_h\tilde\phi }=\intO{\mathbb{S}(\nabla_h\vv):\Grad\phi }
,\quad
\intO{ p_h(\vr) { \Divh} \tilde\phi } = \intO{ p_h(\vr) \Div \phi},
\]
we may rewrite  (\ref{N3})--by using  the formula (\ref{Up3})
and rearanging conveniently several terms-- in the form
\bFormula{me-v}
\intO{ D_t (\vr \widehat{\vv}) \cdot \phi } - \intO{ \vr { \widehat\vu } \otimes \widehat{\vv} : \Grad \phi } 
- \intO{ \vr {\widehat \vu } \cdot\Grad \phi\cdot \widehat{\vu}_B }
 +\intO{\vr\widehat\vu\cdot\nabla_h(\vu_B\cdot\phi)} 
\eF
$$
 +
\intO{\mathbb{S}(\nabla_h\vu):\Grad\phi }
 -\intO{ p_h(\vr) \Div \phi } 
 =<{ R}^M_{h,\Delta t},\phi>\;\mbox{in $(0, T]$},
$$
where
$$
 <{R}^M_{h,\Delta t},\phi>
 = \intO{ D_t( \vr \widehat{\vv}) \cdot (\phi - \widehat{\tilde\phi}) } + 
\sum_{K\in{\cal T}} \sum_{\sigma\in {\cal E}(K)\cap{\cal E}_{\rm int}}\int_\sigma (\phi - \widehat{\tilde\phi}) \cdot 
\ju{ \vr\widehat\vv }_{\sigma,\vc n_{\sigma,K}} \ [\vu_\sigma \cdot \vc{n}_{\sigma,K}]^- {\rm d}S_x 
$$
$$
+ \sum_{K \in {\cal T}} \sum_{\sigma\in {\cal E}(K)}\int_\sigma \vr(\widehat{\tilde\phi}-\phi)\cdot\widehat\vv (\vu -
\vu_\sigma )\cdot \vc{n}_{\sigma,K}{\rm d}S_x + \intO{ \vr(\phi - \widehat{\tilde\phi}) \cdot\widehat\vv  \Divh \vu }
$$
$$
{ + \intO{\vr(\vu-\widehat\vu)\cdot\Grad\phi\cdot\widehat\vv}+ \intO{\vr\widehat\vu\cdot\nabla_h\vu_B\cdot(\phi-\widehat{\tilde\phi})}}
{
 +\intO{\vr\widehat\vu\cdot\nabla\phi\cdot(\vu_B-\widehat\vu_B)}}
$$
$$
{ - \kappa h^\omega \sum_{\sigma\in{\cal E}_{\rm int}}\int_\sigma[[\vr^k]]_{\sigma}\{\widehat\vu^k\}_{\sigma}[[\widehat{\tilde\phi}]]_\sigma{\rm d}S_x}
+  \sum_{\sigma\in {\cal E}^{\rm in}}\int_\sigma
(\vr_B-\vr)\vu_{B,\sigma}\cdot \vc n_{\sigma}{ \widehat\vv}\cdot (\phi - \widehat{\tilde\phi})=\sum_{i=1}^9I_i 
$$

Alternatively, in view of (\ref{pwt+}), equation (\ref{me-v}) can be rewritten as follows,
\bFormula{me-v+}
\intO{\utilde{\vr\widehat\vv}\cdot\phi(\tau,x)}-\intO{{\vr^0\widehat\vv^0}\cdot\phi(0,x)}
- \int_{0}^\tau\intO{\Big[\utilde{\vr\widehat\vv}\cdot\partial_t\phi+ \Big(\vr { \widehat\vu }\otimes \widehat{\vu}+p_h(\vr)\mathbb{I}\Big):\Grad\phi\Big] }{\rm d}t
\eF
$$
+\int_{0}^\tau\intO{\vr\widehat\vu\cdot\nabla_h(\vu_B\cdot\phi)} 
+
\int_{0}^\tau\intO{\mathbb{S}(\nabla_h\vu):\Grad\phi }{\rm d}t=
\int_{0}^\tau<{ R}^M_{h,\Delta t}(t),\phi>{\rm d}t
$$
with any $\tau\in (0, T]$ and all $\phi\in C_c^1([0,T]\times\Omega,R^3)$. 


Our goal is to estimate conveniently all terms in $<R^M_{h,\Delta  t},\phi>$. As in Section \ref{CE},
we shall use the tools reported in Section \ref{Preln} and uniform estimates for  numerical solutions derived in Lemma \ref{Lestimates}.

\begin{enumerate}

\item {\it  Most terms in $<{R}^M_{h,\Delta  t},\phi>$ contain the expression $\widehat{\tilde\phi}-\phi
=\Pi^Q(\tilde\phi-\phi)+\widehat\phi-\phi$}. We notice that by virtue of (\ref{PoincareV-}), (\ref{globalV}), (\ref{errorV}), 
$$
\|\phi-\widehat{\tilde\phi}\|_{L^q(\Omega)}\aleq h\|\Grad\phi\|_{L^q(\Omega)},\; 1\le q\le\infty.
$$
\item {\it Estimate of the first term in ${R}_{h,h}^M$ (term $I_1$)}:
For the error in the time derivative, we obtain
$$
\left| \intO{ D_t (\vr \widehat{\vv}) \cdot (\phi - \widehat{\tilde\phi} ) }\right| \aleq \sqrt{h} (\Delta t)^{-1/2} A_h,\; \mbox{$A_h$ bounded in $L^2(0,T)$}.
$$
Indeed, we split the term  $\intO{ D_t (\vr^k \widehat{\vv}^k) \cdot (\phi - \widehat{\tilde\phi} ) }$ into two parts,
\bFormula{mest0}
\intO{ \sqrt{ \vr^{k-1} } \sqrt{ \vr^{k-1}} \frac{ \vv^k - \vv^{k-1} }{\Delta t} \cdot ( \phi - \widehat{\tilde\phi}) }
+ \intO{ \frac{ \vr^k - \vr^{k-1} }{\Delta t} \vv^{k} \cdot ( \phi - \widehat{\tilde\phi}) },
\eF
where, by virtue of H\"older's inequality, for any $\gamma>1$, the first term is bounded by
\bFormula{mest0+}
\aleq
h (\Delta t)^{-1/2}\| \vr^{k-1} \|_{L^\gamma(\Omega)}^{1/2} \left(\Delta t \intO{ \vr^{k} \left( \frac{ \vv^{k - 1} - \vv^{k-1} }{\Delta t} \right)^2 } \right)^{1/2}
\left\| \Grad\phi  \right\|_{L^{\frac{2 \gamma}{\gamma - 1}}(\Omega)}
\eF
while the second one is
\bFormula{mest0++}
\aleq h^{1-\eta/2}(\Delta t)^{-1/2}  \left( \Delta t\,h^\eta \intO{ \left( \frac{ \vr^k - \vr^{k-1} }{ \Delta t } \right)^2 } \right)^{1/2} \| \vv^k \|_{L^6(\Omega;R^3)}
 \| \Grad \phi \|_{L^3(\Omega)},
\eF
where the first integrals of both expressions are controlled by means of (\ref{e4}) and (\ref{e4+}).

\item {\it Estimate of the second term in ${R}_{h,h}^M$ (term $I_2$)}:
This essentially amounts to estimate two terms,
$$
I_{2,1}= \sum_{K\in{\cal T}} \sum_{\sigma\in {\cal E}(K)\cap{\cal E}_{\rm int}}\int_\sigma \vr^+_{\sigma,\vc n_{\sigma,K}}(\phi - \widehat{\tilde\phi}) \cdot 
\ju{ \widehat\vv }_{\sigma,\vc n_{\sigma,K}} \ [\vu_\sigma \cdot \vc{n}_{\sigma,K}]^- {\rm d}S_x 
$$
and
$$
I_{2,2}= \sum_{K\in{\cal T}} \sum_{\sigma\in {\cal E}(K)\cap{\cal E}_{\rm int}}\int_\sigma (\phi - \widehat{\tilde\phi}) \cdot 
\widehat\vv^-_{\sigma,\vc n_{\sigma,K}}\ju{ \widehat\vr }_{\sigma,\vc n_{\sigma,K}} \ [\vu_\sigma \cdot \vc{n}_{\sigma,K}]^- {\rm d}S_x,
$$ 
where
$$
|I_{2,1}|\aleq \Big[\sum_{K\in{\cal T}} \sum_{\sigma\in {\cal E}(K)\cap{\cal E}_{\rm int}}\int_\sigma \vr^+_{\sigma,\vc n_{\sigma,K}}  
\ju{ \widehat\vv}^2_{\sigma,\vc n_{\sigma,K}} \ |[\vu_\sigma \cdot \vc{n}_{\sigma,K}]^-| {\rm d}S_x\Big]^{1/2}\times 
$$
\bFormula{I21}
\Big[\sum_{K\in{\cal T}} \sum_{\sigma\in {\cal E}(K)\cap{\cal E}_{\rm int}}\int_\sigma \vr^+_{\sigma,\vc n_{\sigma,K}}(\phi - \widehat{\tilde\phi})^2  \ |[\vu_\sigma \cdot \vc{n}_{\sigma,K}]^-| {\rm d}S_x \Big]^{1/2}
\eF
with the first term in the product controlled by (\ref{e5}), i.e. belonging to $L^2((0,T))$, while
$$
|I_{2,2}|\aleq \Big[h^\eta\sum_{K\in{\cal T}} \sum_{\sigma\in {\cal E}(K)\cap{\cal E}_{\rm int}}\int_\sigma 
\ju{ \widehat\vr }_{\sigma,\vc n_{\sigma,K}}^2 \ |[\vu_\sigma \cdot \vc{n}_{\sigma,K}]^-| {\rm d}S_x\Big]^{1/2}\times
$$
\bFormula{I22}
\Big[h^{-\eta}\sum_{K\in{\cal T}} \sum_{\sigma\in {\cal E}(K)\cap{\cal E}_{\rm int}}\int_\sigma  (\phi - \widehat{\tilde\phi})^2 \cdot 
|\widehat\vv^-_{\sigma,\vc n_{\sigma,K}}|^2 \ |[\vu_\sigma \cdot \vc{n}_{\sigma,K}]^-| {\rm d}S_x\Big]^{1/2}
\eF
with the first term in the product controlled by (\ref{e6}), i.e. belonging also to $L^2((0,T))$. It will be therefore enough to estimate
the expressions under the second square roots of $I_{2,1}$, $I_{2,2}$, respectively. In fact, it is enough to estimate only the "leading terms"
of these expressions, where we replace $\vu$ by $\vv$.

We have
$$
h^{-\eta}\Big|\sum_{K\in{\cal T}} \sum_{\sigma\in {\cal E}(K)\cap{\cal E}_{\rm int}}\int_\sigma \vr^+_{\sigma,\vc n_{\sigma,K}}
(\phi - \widehat{\tilde\phi})^2  
\ |[\vv_\sigma \cdot \vc{n}_{\sigma,K}]^-| {\rm d}S_x \Big|
$$
\bFormula{mest1}
\aleq h^{2-\eta} \|\nabla\phi\|_{L^\infty(\Omega))}^2 \sum_{K\in{\cal T}} \sum_{\sigma\in {\cal E}(K)\cap{\cal E}_{\rm int}}
\|\vr\|_{L^2(\sigma)}\|\vv\|_{L^{2}(\sigma)}
\eF
$$
\aleq h^{1-\eta} \|\nabla\phi\|_{L^\infty(\Omega))}^2 \sum_{K\in{\cal T}} 
\|\vr\|_{L^2(K)}\|\vv\|_{L^{2}(K)} 
\aleq h^{1-\frac 32\eta} \|\nabla\phi\|_{L^\infty(\Omega))}^2 
\|h^{\eta/2}\vr\|_{L^2(\Omega)}\|\vv\|_{L^{6}(\Omega)}.
$$
Likewise
$$
\Big|\sum_{K\in{\cal T}} \sum_{\sigma\in {\cal E}(K)\cap{\cal E}_{\rm int}}\int_\sigma  (\phi - \widehat{\tilde\phi})^2 
|\widehat\vv^-_{\sigma,\vc n_{\sigma,K}}|^2 \ |[\vv_\sigma \cdot \vc{n}_{\sigma,K}]^-| {\rm d}S_x\Big|
$$
\bFormula{mest2}
\aleq h \|\Grad\phi\|^2_{L^\infty(\Omega)} \|\vv\|_{L^{3}(\Omega)}^3
\aleq
 h (\Delta t)^{-1/2}  \|(\Delta t)^{1/6}\vv\|_{L^{6}(\Omega)}^3 \|\Grad\phi\|^2_{L^\infty(\Omega)},
\eF
where  $\|(\Delta t)^{1/6}\vv\|_{L^{6}(\Omega)}^3$  is bounded in $L^1((0,T))$ { by virtue of (\ref{timeG})}.
\item {\it The upper bound of the third term in ${ R}_{h,h}^M$ (term $I_3$)} is determined by the upper  bound of
$$
|\sum_{K \in {\cal T}} \sum_{\sigma\in {\cal E}(K)}\int_\sigma \vr(\phi - \widehat{\tilde\phi})\cdot\widehat\vv (\vv -
\vv_\sigma )\cdot \vc{n}_{\sigma,K}{\rm d}S_x|,
$$
which is
\bFormula{mest3}
\aleq h \|\vr\|_{L^2(\Omega)}\|\vv\|_{L^\infty(\Omega)} \|\nabla_h\vv\|_{L^{2}(\Omega)} 
\|\Grad\phi\|_{L^\infty(\Omega)}
\eF
$$
\aleq h^{\frac 12-\frac \eta 2} \|h^{\eta/2}\vr\|_{L^2(\Omega)}\|\vv\|_{L^6(\Omega)} \|\nabla_h\vv\|_{L^{2}(\Omega)} 
\|\Grad\phi\|_{L^\infty(\Omega)},
$$
where $\|\vv\|_{L^6(\Omega)} \|\nabla_h\vv\|_{L^{2}(\Omega)}$ is bounded in $L^1((0,T))$.
\item {\it The bounds of term $I_4$ in ${R}_{h,h}^M$} are determined by the bounds of
$$
\Big|\intO{ \vr(\phi - \widehat{\tilde\phi}) \cdot\widehat\vv  \Divh \vv }\Big|.
$$
They are exactly the same as in the previous case. The same is true for the terms $I_5$--$I_7$, since they have the same structure.

\item The {\it the artificial viscosity  term} $|I_8|=\kappa h^\omega |\sum_{\sigma\in{\cal E}_{\rm int}}$ $\int_\sigma[[\vr^k]]_{\sigma}\{\widehat\vv^k\}_{\sigma}[[\widehat\phi]]_\sigma{\rm d}S_x|$
 is bounded by
\bFormula{mestad}
\aleq\kappa h^\omega\|\vr\|_{L^2(\Omega)}\|\vv\|_{L^2(\Omega)}
\|\nabla\phi\|_{L^\infty(\Omega)} \aleq \kappa h^{\omega-\eta/2}\|h^{\eta/2}\vr\|_{L^2(\Omega)}\|\vv\|_{L^2(\Omega)}
\|\nabla\phi\|_{L^\infty(\Omega)}
\eF
where  $\|\vv\|_{L^2(\Omega)}\aleq\|\vv\|_{L^6(\Omega)}$ is bounded in 
$L^2((0,T))$.

%
 \item {\it The last term in ${ R}_{h,h}^M$ to be evaluated is the boundary term ($I_9$)} whose decay with $h$ is determined by 
$$
\sum_{\sigma\in {\cal E}^{\rm in}}\int_\sigma
\vr\vu_{B,\sigma}\cdot \vc n_{\sigma}\widehat\vv\cdot (\widehat{\tilde\phi}-\phi) 
{\rm d}S_x.
$$
We have, by the same reasoning as in (\ref{cest4}),  the bound
\bFormula{mest4}
\aleq h\|\vr\|_{L^2(\Omega)}\|\nabla_h \vv\|_{L^2(\Omega)}\| \|\Grad\phi\|_{L^\infty(\Omega)}
\aleq
h^{1-\eta/2} \|h^{\eta /2}\vr\|_{L^2(\Omega)} |\nabla_h \vv\|_{L^2(\Omega)} 
\|\Grad\phi\|_{L^\infty(\Omega)}.
\eF
\end{enumerate}

\bLemma{m-consistency}{\rm [Consistency for the momentum equation]}
Let $h=\Delta t$, $\tilde\kappa=1$ and $\eta\in (0,2/3).$ Let the pressure satisfy the hypotheses { (\ref{pres1}--\ref{pres4}).} Then   any numerical solution 
of problem (\ref{N1}--\ref{N3}) 
satisfies the momentum equation in the variational form (\ref{me-v+}) with any $\phi\in C^1([0,T]\times\overline\Omega)$, where
\bFormula{NM1}
\left| \left<  R^M_{h,h} (t), \phi \right> \right| \aleq h^{\alpha_M} r^M_h \|\phi \|_{L^{\infty}(0,T;W^{1,\infty}(\Omega))},\;\|r^M_h\|_{L^1((0,T))}
\aleq 1
\eF
with $\alpha_M=\frac 12-\frac 34\eta >0$ if $\kappa=0$ and $\alpha_M=\min\{\omega-\frac \eta 2, \frac 12-\frac 34\eta\}$ if $\kappa=0$, $\eta<2\omega$. 
\eL

{
\bRemark{Rm-consistency}
{ Let $\tilde\kappa=0$ and let $p$ satisfy in addition to (\ref{pres1}--\ref{pres4}) also condition (\ref{pres2}).}
\begin{enumerate}
\item   Formula (\ref{NM1}) with $\alpha_M>0$ remains still valid under condition $6/5<\gamma<2$ provided  (\ref{pres3}) is satisfied.  In this case,
if $\kappa=0$ and $\gamma>3/2$,
$$
\alpha_M=\min\Big\{\frac 14, \frac{5\gamma-6}{2\gamma}, \frac{2\gamma-3}{2\gamma} \Big\}
$$
and if $\kappa=1$, $6/5<\gamma\le 3/2$ and $\omega\in (0,2(\gamma-1))$,
$$
\alpha_M=\min\Big\{\frac 1 4, \frac{5\gamma-6}{2\gamma},1-\frac 1{2\gamma}-\frac{1+\omega}{2\gamma},\omega \Big\}.
$$
Here are some hints (we suppose $\Delta t=h$ and we use systematically the "negative interpolation" (\ref{interG}), (\ref{timeG})): 
\begin{enumerate}
\item
When estimating the second term in formula (\ref{mest0}), the control (\ref{e4}) must be replaced by the control (\ref{e4**}). Under condition $h=\Delta t$, this
amounts to replace (\ref{mest0++}) by
$$
\aleq h^{-\frac 12}  \left[ \Delta t\,\intO{\Big(\frac{1}{\vr^k+\vr^{k-1}+1}\Big)^{2-\gamma} \left( \frac{ \vr^k - \vr^{k-1} }{ \Delta t } \right)^2 } \right]^{1/2} \Big[\intO{(\vr^k+\vr^{k-1}+1)^{2-\gamma}|\vv^k|^2|\widehat{\tilde\phi}-\phi|^2}\Big]^{1/2}
$$
$$
\aleq h^{\frac 12}\|\vr\|^{\frac {2-\gamma} \gamma}_{L^\gamma(\Omega)} A_h\|\vv^k\|_{L^6(\Omega)} \|\Grad\phi\|_{L^\infty(\Omega)},
$$
provided $\gamma\ge 6/5$, where $A_h\|\vv^k\|_{L^6(\Omega)}$ is bounded in $L^1((0,T))$.
\item Decomposition (\ref{I21}) remains unchanged and yields the bound
$$
\aleq h^{\frac 12} \|\vr\|_{L^\gamma(\Omega)}\|\vu\|_{L^6(\Omega)} A_h \|\Grad\phi\|_{L^\infty(\Omega)}
$$
where $A_h \|\vu\|_{L^6(\Omega)}$ is bounded in $L^1((0,T))$, cf. (\ref{e3}), (\ref{e5}), {and (\ref{mest2}).}

Decomposition (\ref{I22}) takes the form,
$$
\Big[\sum_{K\in{\cal T}} \sum_{\sigma\in {\cal E}(K)\cap{\cal E}_{\rm int}}\int_\sigma \Big( \frac 1{\{\vr\}_\sigma+1}\Big)^{2-\gamma} 
\ju{ \widehat\vr }_{\sigma,\vc n_{\sigma,K}}^2 \ |[\vu_\sigma \cdot \vc{n}_{\sigma,K}]^-| {\rm d}S_x\Big]^{1/2}\times
$$
$$
\Big[\sum_{K\in{\cal T}} \sum_{\sigma\in {\cal E}(K)\cap{\cal E}_{\rm int}}\int_\sigma (\{\vr\}_\sigma+1)^{2-\gamma}  (\phi - \widehat{\tilde\phi})^2 \cdot 
|\widehat\vv^-_{\sigma,\vc n_{\sigma,K}}|^2 \ |[\vu_\sigma \cdot \vc{n}_{\sigma,K}]^-| {\rm d}S_x\Big]^{1/2},
$$
where the first term is controlled by (\ref{e6**}). 
Under assumption $\Delta t=h$, the expression under  the second square root  (where we take $\vv$ on place of $\vu$
without loss of generality)
admits the bound
$$
\aleq h^{\frac 12+\min\{0,\frac 9{p_m}-\frac 32\}}\|\vr\|^{2-\gamma}_{L^\gamma(\Omega)}  \|(\Delta t)^{1/6}\vv\|_{L^{6}(\Omega)}^3 \|\Grad\phi\|^2_{L^\infty(\Omega)},
$$
where $\frac 3{p_m}+\frac{ 2-\gamma}\gamma=1$ and $\|(\Delta t)^{1/6}\vv\|_{L^{6}(\Omega)}^3$ is bounded in $L^1((0,T))$ and the power of $h$ is positive provided $\gamma>6/5$. \footnote{This is the place, where we cannot
get a positive power of $h$ for $\gamma\le 6/5$ even if we take $\kappa=1$ and employ estimate (\ref{e8**}). This is not due to the boundary conditions. Indeed, the same limitation is observed for the  no-slip boundary conditions.}

\item The estimate (\ref{mest3}) is modified as follows,
$$
\aleq h\|\vr\|_{L^\infty(\Omega)}^{1/2}\|\vr\widehat\vv^2\|_{L^1(\Omega)}^{1/2}\|\nabla_h\vv\|_{L^2(\Omega)}\|\Grad\phi\|_{L^\infty(\Omega)}
\aleq h^{1-\frac 3{2\gamma}} \|\nabla_h\vv\|_{L^2(\Omega)}\|\Grad\phi\|_{L^\infty(\Omega)},
$$
where $\|\nabla_h\vv\|_{L^2(\Omega)}$ is bounded in $L^2((0,T))$ and the power of $h$ is positive provided $\gamma>3/2$.

Alternatively, if $\kappa=1$, we have in this case the bound
$$
\aleq h^{1-\frac 1{2\gamma}-\frac{1+\omega}{2\gamma}}\Big(h^{\frac{ 1+\omega}\gamma}\|\vr\|_{L^{3\gamma}(\Omega)}\Big)^{1/2}
\|\vr\widehat\vv^2\|_{L^1(\Omega)}^{1/2} \|\nabla_h\vv\|_{L^2(\Omega)}\|\Grad\phi\|_{L^\infty(\Omega)},
$$
where $\Big(h^{\frac{ 1+\omega}\gamma}\|\vr\|_{L^{3\gamma}(\Omega)}\Big)^{1/2}
\|\vr\widehat\vv^2\|_{L^1(\Omega)}^{1/2} \|\nabla_h\vv\|_{L^2(\Omega)}$ is bounded in $L^1((0,T))$ and the power of $h$
is positive provided $\omega\in (0, 2(\gamma-1))$. In the above, we have used, in particular, estimate (\ref{e8**}).

\item The estimate (\ref{mestad}) now reads
$$
\aleq\kappa h^\omega\|\vr\|_{L^\gamma(\Omega)}\|\vv\|_{L^6(\Omega)}
\|\nabla\phi\|_{L^\infty(\Omega)},
$$
provided $\gamma\ge 6/5$.

\item Finally, by the argumentation of (\ref{cest4}), estimate (\ref{mest4}) can be modified as
$$
\aleq h^{\frac{5\gamma-6}{6\gamma}} \|\vr\|_{L^\gamma(\Omega)} \|\vv\|_{L^{6}(\Omega)} 
\|\Grad\phi\|_{L^\infty(\Omega)}.
$$
\end{enumerate}
\item Similar reasoning can be carried out if $\gamma\ge 2$ and in (\ref{pres2}), $\underline\pi>0$.
In this case, we can take also advantage  of estimates (\ref{e6*}) and (\ref{e4*}). We obtain, at least,
$$
\alpha_M=\frac 12\;\mbox{if $\kappa=0$}\;\mbox{and}\; \alpha_M=\min\{\frac 12,\omega\}\;\mbox{if $\kappa=1$}.
$$
\end{enumerate}
\eR
}

\subsection{Energy balance}

Neglecting several positive terms at the left hand side and using the definition of { $\Pi^V$ and (\ref{fit}) in order to replace $\vu_{B,\sigma}=\tilde{\mathfrak{u}}_{B,\sigma}$ by $\mathfrak{u}_B$ in the
remaining boundary integrals,} we can rewrite the energy identity (\ref{ebalance}) 
in the following form,
\bFormula{ee-v}
\Big[\intO{\Big(\frac 12\vr |\widehat\vv|^2 + H_h(\vr)\Big) }\Big]_0^\tau +\int_0^\tau\intO{\mathbb{S}(\nabla_h\vu):\nabla_h\vv}{\rm d}t
\eF
$$
+\int_0^\tau\int_{\Gamma^{\rm out}} H_h(\vr){\mathfrak{u}}_{B}\cdot\vc n
{\rm d}S_x {\rm d}t
 \aleq
  -\int_0^\tau\int_{\Gamma^{\rm in}} { H_h(\vr_B)}{\mathfrak{u}}_{B}\cdot\vc n
{\rm d}S_x {\rm d}t
$$
$$
- \int_0^\tau\intO{\Big(\vr\widehat\vu\otimes\widehat\vu+p_h(\vr)\mathbb{I}\Big):\nabla_h\vu_B}{\rm d}t
+\int_0^\tau\intO{\vr\widehat\vu\cdot\nabla_h\vu_B\cdot\widehat\vu_B}{\rm d}t  
 + { R}_h^E[\vr,\vu],
 $$
where
$$
\Big[{ R}_h^E[\vr,\vu]\Big](\tau)=
\int_\tau^{\tau_m}\Big[\int_{\Gamma^{\rm in}} { H_h(\vr_B)}|\vu_{B}\cdot\vc n|
{\rm d}S_x { + \int_{\Gamma^{\rm out}}  (d\vr+e)|\vu_{B}\cdot\vc n|
{\rm d}S_x} 
$$
$$
- \intO{\Big(\vr\widehat\vu\otimes\widehat\vu+p_h(\vr)\mathbb{I}\Big):\nabla_h\vu_B}
+\intO{\vr\widehat\vu\cdot\nabla_h\vu_B\cdot\widehat\vu_B}\Big]{\rm d}t,\;
\tau\in(\tau_{m-1},\tau_m],\; m=1,\ldots,N ; 
$$
whence
\bFormula{enb3}
\Big|{ R}_h^E[\vr,\vu]\Big|\aleq h,\;{ \mbox{ where $\Delta t=h$} }
\eF
by virtue of (\ref{e0}--\ref{e3}), cf. also { (\ref{mbalance}), (\ref{enb1}), (\ref{calH*}).}

\section{Convergence to a dissipative solution}\label{SConv}

We denote by $[\vr_h,\vu_h]$, $h>0$ a sequence of numerical solution to  the scheme (\ref{N1}--\ref{N3}). We want to show that
 there is a  subsequence with weak limit $[{\mathfrak{r}},\mathfrak{u}]$, and there is a positively semi-definite tensor measure
$\mathfrak{R}$ (which we shall construct as well), such that the couple $[{\mathfrak{r}},\mathfrak{u}]$ and 
the associated $\mathfrak{R}$ is a dissipative solution of the continuous problem (\ref{NS1}--\ref{inB}) in the sense of Definition
\ref{DD1}.
\subsection{Weak limits. Continuity equation}
Recalling regularity (\ref{ru0}--\ref{ruB}) of the initial and boundary data, we deduce from (\ref{errorQ}--\ref{PoincareV-}), (\ref{globalV}) (\ref{errorV}),
\bFormula{co0-}
\vr_{B,h}\to\mathfrak{r}_B\;\mbox{in $L^q(\Omega)$ and $L^q(\partial\Omega)$, $1\le q\le\infty$},
\eF 
$$
\widehat\vu_{B,h}\to \mathfrak{u}_B,\;\vu_{B,h}\to\mathfrak{u}_B,\;\nabla_h\vu_{B,h}\to\Grad\mathfrak{u}_B\;\mbox{in $L^q(\Omega)$, $1\le q\le \infty$}, 
$$
$$
\vr^0_h\to \mathfrak{r}_0\;\mbox{in $L^q(\Omega)$}, \;\vu^0_h\to\mathfrak{u}_0\;\;\mbox{in $L^q(\Omega)$},\;
1\le q\le\infty.
$$

By virtue of estimate { (\ref{e3}) and Jensen's inequlaity (\ref{JensenV})},
\bFormula{co0+}
\|\widehat\vv_h\|_{L^2(0,T;L^{6}(\Omega))}\aleq 1.
\eF
{Consequently,\footnote{All convergences in this section hold for a chosen subsequence of the
original sequence; for the sake of simplicity, we do not relabel.}
\bFormula{co1+}
\widehat\vv_h\rightharpoonup\mathfrak {v},\;
\vv_h\rightharpoonup\mathfrak{v}
\;\mbox{in $L^2(0,T;L^6(\Omega))$}
\eF
again  by virtue of (\ref{e3})}, where the limit of both sequences is the same by virtue of (\ref{globalV}). { Further},
$$
\intO{\vv_h\Grad\phi}=-\intO{\nabla_h\vv_h\phi}+I_h +J_h,\;\phi\in C^1((0,T)\times\overline\Omega)),
$$
where
{
$$
I_h:=\sum_{K\in {\cal T}}\sum_{\sigma\in {\cal E}(K)}\int_\sigma\vv_h\vc n_{\sigma,K}\phi{\rm d}S_x=
\sum_{K\in {\cal T}}\sum_{\sigma\in {\cal E}(K)\cap{\cal E}_{\rm int}}\int_\sigma(\vv_h-\vv_{h,\sigma})\vc n_{\sigma,K}(\phi-\phi_\sigma){\rm d}S_x,\quad
J_h=
 \int_{\partial\Omega}\vv_h\vc n\phi{\rm d}S_x
$$
admit the bounds
\bFormula{dod*}
 |I_h|\aleq h\|\nabla_h\vv_h\|_{L^2(0,T;L^2(\Omega))}\|\phi\|_{L^2(0,T;W^{1,2}(\Omega))}, \;
|J_h|\aleq h^{1/3}\|\nabla_h\vv_h\|_{L^2(0,T;L^2(\Omega))}\|\phi\|_{L^2(0,T;W^{1,2}(\Omega))}
\eF
by virtue of the H\"older inequality, trace estimates (\ref{trace}--\ref{traceFE}), 
the first inequality in (\ref{PoincareV-}), and in the second estimate also
the standard Sobolev imbedding $W^{1,2}(\Omega)\hookrightarrow L^6(\Omega)$,  the fact that $\vv_h\in V_0(\Omega;R^3)$ and that $|\cup_{K\cap {\cal E}_{\rm ext}\neq \emptyset} K|\aleq h.$}
We deduce from this and from (\ref{e3}), on one hand,
\bFormula{co2}
\nabla_h\vv_h\rightharpoonup\Grad\mathfrak{v}\;\mbox{in $L^2(0,T;L^2(\Omega))$},\;\mbox{and $\mathfrak{v}\in
L^2(0,T;W^{1,2}_0(\Omega))$,}
\eF
and on the other hand, in particular,
\bFormula{co3}
\|\Grad\vv_h\|_{L^2(0,T;W^{-1,2}(\Omega))}\aleq 1,\quad \|\Grad\widehat\vv_h\|_{L^2(0,T;W^{-1,2}(\Omega))}\aleq 1,
\eF
where the very latter bound is derived from the previous one by virtue of (\ref{globalV}).

According to (\ref{e0}), (\ref{e2}),
\bFormula{co0}
\|\vr_h\widehat\vv_h\|_{L^\infty(0,T;L^{\frac{2\gamma}{\gamma+1}}(\Omega))}\aleq 1.
\eF
It is the consequence of (\ref{pwt+}) and (\ref{e0}) that
\bFormula{etildevr}
\utilde{\vr_h}\;
\mbox{ is bounded in}\; L^\infty(0,T;L^\gamma(\Omega)).
\eF

Coming back with (\ref{co0}) and (\ref{e0}) to the continuity equation (\ref{ce-v+}), we find,
by virtue of Lemma \ref{Lc-consistency}, that for all $\phi\in C^1_c(\Omega)$,
$$
\intO{\utilde{\vr_h}(t)\phi}=A^\phi_h(t)+B^\phi_h(t),
$$
where $t\mapsto A^\phi_h(t)$ is equi-bounded and equi-continuous in $C([0,T])$ while $B^\phi_h\to 0$ in $C([0,T])$. Consequently,
by density of $C^1_c(\Omega)$ in $L^{\gamma'}(\Omega)$ and by the Arzela-Ascoli type argument, we get
\bFormula{co4}
\utilde{\vr_h}\to\mathfrak{r}\;\mbox{in $C_{\rm weak}([0,T];L^\gamma(\Omega))$}.
\eF
We have also, by virtue of (\ref{e0}),
\bFormula{co4+}
{\vr_h}\rightharpoonup_*\mathfrak{r}\;\mbox{in $L^\infty(0,T;L^\gamma(\Omega))$}.
\eF
{ The limits (\ref{co4}), (\ref{co4+}) are the same, since $\|\utilde{\vr_h}-\vr_h\|_{L^1((0,T)\times\Omega)}\aleq (\Delta t)^{1/2} h^{-\eta/2}$ due to
(\ref{pwt+}) and (\ref{e4}).}

{Due to (\ref{co0}),
\bFormula{dodano}
\vr_h\widehat\vv_h\rightharpoonup_*\mathfrak{s}\;\mbox{in $L^\infty(0,T;L^{\frac{2\gamma}{\gamma+1}}(\Omega))$},\;\;\mathfrak{s}=\mathfrak{r}\mathfrak{v}
\eF
where the latter identity follows from 
Lemma \ref{Lemma1}, where
we set $r_n\approx\vr_h$, $v_n\approx\widehat\vv_h$, $\vc g_n\approx\vr_h\widehat\vu_h$, $h_n\approx R^C_{h}$. (Indeed,
one easy checks by using (\ref{NC1+}), (\ref{co0}), (\ref{co1+}), (\ref{co3}), (\ref{co4}), that assumptions of the Lemma \ref{Lemma1}
are satisfied.\footnote{This is the only point, where we need "additional" $L^p$-estimate (\ref{NC1+}).})
}

Finally, due to estimate (\ref{e3})
\bFormula{crhob}
\vr_h\rightharpoonup\mathfrak{r}\;\mbox{in $L^\gamma(0,T;L^\gamma(\partial\Omega,|\mathfrak{u}_B\cdot\vc n|{\rm d}S_x))$}.
\eF

At this stage we can pass to the limit in the consistency formulation (\ref{ce-v+}) of the continuity equation. 
Due to Lemmma \ref{Lc-consistency}
$$
\int_0^\tau<R^C_{h},\phi>{\rm d}t\to 0
$$
and
 we obtain 
the weak formulation (\ref{D2}).

\subsection{Limit in the momentum equation}

According to  (\ref{e0}), (\ref{e0+}), 
\bFormula{co5}
\|p_h(\vr_h)\|_{L^\infty(0,T;L^1(\Omega))}\aleq 1.
\eF
In view of (\ref{pwt+}) and (\ref{co0}),
\bFormula{etildevrv}
\utilde{\vr_h\widehat\vv_h}\;
\mbox{is bounded in}\; \mbox{$L^\infty(0,T;L^{\frac{2\gamma}
{\gamma+1}}(\Omega))$.}
\eF

Coming back with (\ref{e2}--\ref{e3}), (\ref{co5}) to the momentum equation (\ref{me-v+}), we obtain by the same arguments as in (\ref{co4}), employing now Lemma \ref{Lm-consistency} (instead of Lemma \ref{Lc-consistency}),
\bFormula{co6}
\utilde{{\vc q}_h}\to{\mathfrak{q}}\;\mbox{in $C_{\rm weak}([0,T];L^{\frac{2\gamma}{\gamma+1}}(\Omega))$},\; {\vc q}_h:=\vr_h\widehat\vv_h,
\eF
where
$$
{\mathfrak{q}}=\mathfrak{r}\mathfrak{v}\;\mbox{a.e. in $(0,T)\times\Omega$}.
$$
Indeed, the latter identity can be deduced from (\ref{dodano}) and from
the inequality
$$
\|\utilde{\vr_h\widehat\vv_h}-\vr_h\widehat\vv_h\|_{L^1((0,T)\times\Omega)}\aleq (\Delta t)^{1/2}h^{-\eta/2}
$$
which follows from (\ref{pwt+}) and (\ref{e4}--\ref{e4+}).

Due to (\ref{co0-}) and (\ref{co6}), also,
\bFormula{co6+}
\utilde{{\vc m}_h}\to{\mathfrak{m}}\;\mbox{in $C_{\rm weak}([0,T];L^{\frac{2\gamma}{\gamma+1}}(\Omega))$},\;{\vc m}_h:=\vr_h\widehat\vu_h
\eF
where
\bFormula{co6++}
{\mathfrak{m}}=\mathfrak{r}\mathfrak{u}\;\mbox{a.e. in $(0,T)\times\Omega$}\,\mbox{and}\,\mathfrak{u}=\mathfrak{v}+\mathfrak{u}_B.
\eF

We introduce
$$
\mathbb{E}(r,\vc z)=
\left\{
\begin{array}{c}
\frac{\vc{ z}\otimes\vc{ z}}r,\;\mbox{if $\vc z\in R^3$, $r>0$}\\
0\; \mbox{if $\vc z\in R^3$, $r=0$}\\
+\infty\; \mbox{if $\vc z\in R^3$, $r<0$}
\end{array}
\right\},
$$
$$
E^{(\xi)}(r,{\vc z}):=\xi^T\mathbb{E}(r,\vc z)\xi
=\left\{
\begin{array}{c}
\frac{|\vc z\cdot\xi|^2}r,\;\mbox{if $\vc z\in R^3$, $r>0$}\\
0\; \mbox{if $\vc z\in R^3$, $r=0$}\\
+\infty\; \mbox{if $\vc z\in R^3$, $r<0$}
\end{array}
\right\}.
$$
In the above $\xi\in R^3$. We observe that $(r,\vc z)\mapsto E^{(\xi)}(r,\vc z)$ is a lower semi-continuous convex function on $R^4$ and that
$$
\mathbb{E}(\vr_h,\widehat{\vc m}_h)=\vr_h\widehat\vu_h\otimes\widehat\vu_h.
$$
whence, recalling estimate (\ref{e2}), we obtain, by a sequential version of the Banach-Alaoglu-Bourbaki theorem
\bFormula{c7}
\mathbb{E}(\vr_h,\widehat{\vc m}_h)\rightharpoonup_* \overline{\mathbb{E}(\mathfrak{r},\mathfrak{m})}\;\mbox{in 
$L^\infty(0,T; {\cal M}(\overline\Omega; R_{\rm sym}^{3\times 3}))$}
\eF
while\footnote{Here and hereafter $\overline{f(\mathfrak{r},\mathfrak{v},\Grad\mathfrak{v},\ldots)}$
denotes a star-weak limit of the sequence $f(\vr_h,\vv_h,\nabla_h\vv_h,\ldots)$.}
\bFormula{c8}
\overline{\xi^T\mathbb{E}(\mathfrak{r},\widehat{\mathfrak m})\xi}-{E}^{(\xi)}(\mathfrak{r},\mathfrak{m})\ge 0,\quad\mbox{and ${E}^{(\xi)}(\mathfrak{r},\mathfrak{m})\in L^\infty(0,T;L^1(\Omega))$}.
\eF
by virtue of Lemma \ref{Lemma2}. Consequently, seeing (\ref{co6++}),
$$
\mathbb{E}(\mathfrak{r},\mathfrak{m})=\mathfrak{r}\mathfrak{u}\otimes\mathfrak{u}\in L^\infty(0,T;L^1(\Omega)).
$$
We reserve the same treatment to the sequence $p_h(\vr_h)\mathbb{I}$. Due to (\ref{co5}), and since, by (\ref{pres4}), (\ref{ph}), $p_h$ is convex, we have,
\bFormula{c9}
p_h(\vr_h)\mathbb{I}\rightharpoonup_* \overline {p(\mathfrak{r})}\mathbb{I} \;\mbox{in $L^\infty(0,T; {\cal M}(\overline\Omega; R_{\rm sym}^{3\times 3}))$},
\eF
where
$$
\overline {p(\mathfrak{r})}-p(\mathfrak{r})\ge 0,\;p(\mathfrak{r})\in L^\infty(0,T; L^1(\Omega)).
$$

Last but not least, we introduce 
\bFormula{co10}
\mathfrak{R}:=\Big(\overline{\mathbb{E}(\mathfrak{r},\mathfrak{m}}) +\overline{p(\mathfrak{r})}\mathbb{I}\Big)-
\Big(\mathbb{E}(\mathfrak{r},\mathfrak{m}) + p(\mathfrak{r})\mathbb{I}\Big);
\eF
in view of (\ref{c8}--\ref{c9}),
$$
\mbox{for all $\xi\in R^3$},\; \xi^T\mathcal{ R}\xi \in L^\infty(0,T,{\cal M}^+(\overline\Omega))\;\mbox{i.e., $ \mathfrak{R}\in L^\infty(0,T;{\cal M}^+ (\overline\Omega; R_{\rm sym}^{3\times 3}))$}.
$$

Finally, due to Lemma \ref{Lm-consistency},
$$
\int_0^\tau<R^M_{h,h},\phi>{\rm d}t\to 0.
$$

Now, we are ready to pass to the limit in the momentum equation (\ref{me-v+}), in order to obtain
\bFormula{me-d}
\intO{{\mathfrak{q}}\cdot\phi(\tau,x)}-\intO{\mathfrak{r}_0\mathfrak{v}_0\cdot\phi(0,x)}
\eF
$$
- \int_0^\tau\intO{\Big(\mathfrak{r}\mathfrak{v}\cdot\partial_t\phi+ \Big(\mathfrak{r} { \mathfrak{u} } \otimes \mathfrak{u}+p(\mathfrak{r})\mathbb{I}\Big) : \Grad \phi\Big) }{\rm d}t
$$
$$
 +\int_0^\tau\intO{\mathfrak{r}\mathfrak{u}\cdot\nabla(\mathfrak{u}_B\cdot\phi)} {\rm d}t
 +
\int_0^\tau\intO{ \mathbb{S}(\Grad\mathfrak{u}):\Grad\phi}{\rm d}t=
\int_0^\tau\int_\Omega\Grad\phi:{\rm d}\mathfrak{R}(t){\rm d}t.
$$
According to (\ref{D2}),
$$
\int_0^\tau\intO{\mathfrak{r}\mathfrak{u}\cdot\Grad(\mathfrak{u}_B\cdot\phi)}{\rm d}t=-\int_0^\tau\intO{\mathfrak{r}\mathfrak{u}_B\cdot\partial_t\phi}{\rm d}t
+\Big[\intO{\mathfrak{r}\mathfrak{u}_B\cdot\phi}\Big]_0^\tau
$$
we obtain from (\ref{me-d}) the formulation (\ref{D3}) of the momentum equation.

\subsection{Limit in the  energy balance}\label{LEI}
Due to estimate (\ref{e2}) 
$$
{\rm Tr}[{\mathbb{E}}(\vr_h,\vr_h\widehat\vv_h)]=\frac{|\vr_h\widehat\vv_h|^2}{\vr_h}\rightharpoonup_*\overline {{\rm Tr}[{\mathbb{E}}(\mathfrak{r},{\mathfrak{q}})]}:=
\overline{\Big(\frac {{\mathfrak{q}}^2}{\mathfrak{r}}\Big)}\;\mbox{in $L^\infty(0,T;{\cal M} (\overline\Omega))$},
$$
where, for almost all $\tau\in (0,T)$,
$$
\overline{\Big(\frac {{\mathfrak{q}}^2}{\mathfrak{r}}\Big)}(\tau)\ge {\rm Tr}[\mathbb{E}(\mathfrak{r}(\tau),{\mathfrak{q}}(\tau))]=\mathfrak{r}{\mathfrak{v}}^2(\tau)
$$
by virtue of Lemma \ref{Lemma2}. Likewise, by the same token, due to estimate (\ref{e0})
$$
H_h(\vr_h)\rightharpoonup_* \overline{H(\mathfrak{r})}\;\mbox{in $L^\infty(0,T;{\cal M} (\overline\Omega))$},\; \overline{H(\mathfrak{r})}\ge H(\mathfrak{r}).
$$
We can thus introduce
$$
\mathfrak{E}= \Big[\frac 12\overline{\Big(\frac {{\mathfrak{q}}^2}{\mathfrak{r}}\Big)}+\overline{H(\mathfrak{r})}\Big]
- \Big[\frac 12\mathfrak{r}{\mathfrak{v}}^2+{H(\mathfrak{r})}\Big]\in L^\infty(0,T;{\cal M}^+(\overline\Omega)).
$$
Clearly, in definition  of $\mathfrak{E}$, one can replace $\mathfrak{q}$ by $\mathfrak{m}$ and $\mathfrak{v}$ by $\mathfrak{u}$.
Consequently, due to the structural property (\ref{pres4}) of the pressure $p$,
\bFormula{co11}
0\le D{\rm Tr}(\mathfrak{R})\le \mathfrak{E}\;\mbox{where $D=\min\{1/2, {\underline a}/3\}$}.
\eF
{
Finally, we use (\ref{crhob}) in conjonction with  Lemma \ref{Lemma2} to show
$$
\int_0^\tau\int_{\Gamma^{\rm out}} H(\mathfrak{r})\mathfrak{u}_B\cdot\vc n{\rm d}S_x
{\rm d}t\le \liminf_{h\to 0} \int_0^\tau\int_{\Gamma^{\rm out}} H(\vr_h)\mathfrak{u}_B\cdot\vc n{\rm d}S_x{\rm d}t.
$$
}
We can now pass to the limit in the  consistency formulation (\ref{ee-v}) of the energy balance, in order to get energy inequality
in the form (\ref{D4}).

This finishes the proof of the first part of Theorem \ref{TN2}.

\section{Strong convergence to a strong solution}\label{SSt}

Let $[\vr_h,\vu_h]$ be any subsequence (not relabeled) of the sequence
$[\vr_h,\vu_h]$. Than it contains a subsequence converging weakly (in the sense (\ref{A1}--\ref{A4})) to a dissipative solution $[\mathfrak{r}, \mathfrak{u}]$
with defect $\mathfrak{ R}$.

Let $[r,\vc U]$ be a strong solution of the problem (\ref{NS1}--\ref{inB}) emanating from the same initial and boundary data
as the dissipative solution $[\mathfrak{r},\mathfrak{u}]$. Then according to Lemmma \ref{WUT1}, $\mathfrak{r}=r$, $\mathfrak{u}=\vc U$
and the Reynolds defect $\mathfrak{R}=0$.

We have one one hand: Any subsequence of the whole sequence $[\vr_h,\vu_h]$ 
admits a weakly convergent subsequence -- in the sense (\ref{A1}--\ref{A3}) --converging to $[r,\vc U]$;  whence the whole sequence weakly converges to
$[r,\vc U]$.

On the other hand $\mathfrak{R}=0$ which implies immediately the convergence
(\ref{A4+}--\ref{A5}) by virtue of (\ref{co10}), (\ref{co11}).

This finishes the proof of the first item in the second part of Theorem \ref{TN2}.

\section{Relative energy for the numerical solution}\label{RE}

The relative energy is a basic tool for showing explicit convergence rate  of the strong convergence of numerical solutions to a strong solution
(provided that the latter exist). Our goal is to evaluate the time evolution of 
\[
{\cal E}_h \left(\vr, \widehat\vv \ \Big|\ r, \vc V\right)= {\cal E} \left(\vr, \widehat\vv \ \Big|\ r, \vc V\right) +h^\eta\intO{(\vr-r)^2}
\]
where the relative energy functional ${\cal E}(\cdot,\cdot|\cdot,\cdot)$ is introduced in (\ref{calE}), $[\vr, \vu]$ is a numerical solution (of system (\ref{N1}--\ref{N3})) and $[r, \vc V]$ are test functions in the class (\ref{tests}).
To achieve this goal we will be mimicking the proofs from the "continuous case", see \cite{FeJiNo}, \cite{KwNo}, \cite{AbFeNo} using the consistency formulations
(\ref{ce-v+}), (\ref{me-v+}) and (\ref{ee-v}) of the continuity, momentum and energy balance laws. Before starting, to this end, we have, however, make compatible the consistency formulations of the continuity and momentum equations with the formulation (\ref{ee-v}) of the energy inequality.

In view (\ref{pwt}--\ref{pwt+}), we deduce from (\ref{ce-v+}),
\bFormula{ce-v++}
\intO{{\vr}\phi(\tau)}-\intO{{\vr}^0\phi(0)}-\int_{0}^\tau\intO{\Big(\vr\partial_t\phi+\vr\widehat\vu\cdot\Grad\phi\Big)}{\rm d}t
+  \int_{0}^{\tau}\int_{\Gamma^{\rm out}}
\vr\mathfrak{u}_{B}\cdot \vc n\phi{\rm d}S_x{\rm d}t
\eF
$$
+  \int_{0}^{\tau}\int_{\Gamma^{\rm in}}
\vr_B\mathfrak{u}_{B}\cdot \vc n\phi{\rm d}S_x{\rm d}t
=<\Pi_{h,\Delta t}^C(\tau),\phi>,\;\tau\in (0,T],
$$
with  all $\phi\in C^1([0,T]\times\overline\Omega)$,
where
$$
<\Pi^C_{h,\Delta t}(\tau),\phi>=
\intO{\vr(\tau)(\phi(\tau)-\phi(\tau_m))}+
\int_0^{\tau_m}<R^C_{h,\Delta t},\phi> {\rm d}t
$$
$$
+\int_0^{\tau_m}
\intO{(\utilde\vr-\vr)\partial_t\phi}
{\rm d}t 
+ \int_{\tau}^{\tau_m}\Big[
\intO{\Big(\vr\partial_t\phi+\vr\widehat\vu\cdot\Grad\phi\Big)} 
$$
$$
-\int_{\Gamma^{\rm out}}
\vr\mathfrak{u}_{B}\cdot \vc n\phi{\rm d}S_x-\int_{\Gamma^{\rm in}}
\vr_B\mathfrak{u}_{B}\cdot \vc n\phi{\rm d}S_x\Big]{\rm d}t,\; \tau\in (\tau_{m-1},\tau_m], \;m=1,\ldots, N.
$$
We have 
$$
\int_0^{\tau_m}
\intO{|\utilde\vr-\vr| \,|\partial_t\phi|}
{\rm d}t\aleq \Delta t\sum_{k=1}^m\|\vr^k-\vr^{k-1}\|_{L^2(\Omega)}
\|\partial_t\phi\|_{L^\infty((0,T)\times\Omega)}
$$
$$
\aleq (\Delta t)^{1/2} \Big(\sum_{k=1}^m\|\vr^k-\vr^{k-1}\|^2_{L^2(\Omega)}\Big)^{1/2}\aleq (\Delta t)^{1/2}h^{-\eta/2}
$$
by virtue of (\ref{pwt+}), where the last inequality follows from (\ref{e4}).
Employing in estimating of the other terms (\ref{e0+}), (\ref{e1+}), (\ref{e2}) and (\ref{NC1}) we get
\bFormula{ce-v++*}
\Big|<\Pi^C_{h,h}(\tau),\phi>\Big|\aleq h^{\alpha_C}\|\phi,\Grad\phi,\partial_t\phi\|_{L^\infty((0,T)\times\Omega)}.
\eF
Concerning the momentum equation, we deduce from (\ref{me-v+}) by the same token,
\bFormula{me-v++}
\intO{{\vr\widehat\vv}\cdot\phi(\tau,x)}-\intO{{\vr^0\widehat\vv^0}\cdot\phi(0,x)}
\eF
$$
- \int_{0}^\tau\intO{\Big[\vr\widehat\vv\cdot\partial_t\phi+ \Big(\vr { \widehat\vu }\otimes \widehat{\vu}+p_h(\vr)\mathbb{I}\Big):\Grad\phi\Big] }{\rm d}t
+\int_{0}^\tau\intO{\vr\widehat\vu\cdot\nabla_h(\vu_B\cdot\phi)}{\rm d}t
$$
$$ 
+
\int_{0}^\tau\intO{\mathbb{S}(\nabla_h\vu):\Grad\phi }{\rm d}t=
<\Pi^M_{h,\Delta t}(\tau),\phi>,\;\tau\in (0,T],
$$
with  all $\phi\in C_c^1([0,T]\times\Omega,R^3)$, where
\bFormula{cm-v++*}
<\Pi^M_{h,\Delta t}(\tau),\phi>=
\intO{\vr\widehat\vv(\tau)\cdot(\phi(\tau)-\phi(\tau_m))}+
\int_0^{\tau_m}<R^M_{h,\Delta t},\phi>{\rm d}t 
\eF
$$
+\int_0^{\tau_{m}}\intO{(\utilde{\vr\widehat\vv}-\vr\widehat\vv)\partial_t\phi}{\rm d}t
+\int_{\tau}^{\tau_m}\Big[\intO{\Big((\vr { \widehat\vu }\otimes \widehat{\vu}+p_h(\vr)\mathbb{I}):\Grad\phi
-\vr\widehat\vu\cdot\nabla_h(\vu_B\cdot\phi)\Big)}
$$
$$
+ \intO{\mathbb{S}(\nabla_h\vu):\Grad\phi}
\Big] {\rm d}t;\; \tau\in (\tau_{m-1},\tau_m], \;\tau_m=m\Delta t, \;m=1,\ldots, N;
$$
whence, using (\ref{e0+}), (\ref{e1+}), (\ref{e2}), (\ref{e4}--\ref{e4+}) and (\ref{NM1}), we get
\bFormula{me-v++*}
\Big|<\Pi^M_{h,h}(\tau),\phi>\Big|\aleq h^{\alpha_M}\|\phi,\Grad\phi,\partial_t\phi\|_{L^\infty((0,T)\times\Omega)}.
\eF

\subsection{Relative energy inequality}

 We calculate,
\begin{center}
(\ref{ee-v}) + (\ref{ce-v++*})$_{\phi=\frac 12|\vc V|^2- H_h'(r)}$ $+$ (\ref{me-v++*})$_{\phi=-\vc V}$
\end{center}
\begin{enumerate}
\item We denote $\vc U=\vc V+\mathfrak{u}_B$ and recall that $\vu=\vv+\tilde{\mathfrak{u}}_B$.
\item The above expression gives the following inequality
\begin{equation} \label{RR4}
\left[ \intO{\mathcal{E}_h\left( \vr, \widehat\vv \ \Big|\ r, \vc V \right) } \right]_{t = 0}^{ t = \tau} + 
\int_0^\tau \intO{ \mathbb{S}(\nabla_h\vu):\nabla_h(\vu-\vc U) } \dt 
\end{equation} 
$$
+\int_0^\tau \int_{\Gamma_{\rm out}} \left[ H_h(\vr) - H_h'(r) \vr \right]  \mathfrak{u}_B \cdot \vc{n} \ {\rm d} S_x \dt
+\int_0^\tau \int_{\Gamma_{\rm in}} \left[ H_h(\vr_B) - H_h'(r) \vr_B \right]  \mathfrak{u}_B \cdot \vc{n} \ {\rm d} S_x \dt
$$
$$
\leq 
- \int_0^\tau \intO{ \Big[ \vr \widehat\vu \cdot \partial_t \vc U  + \vr \widehat\vu\cdot\Grad\vc U\cdot\widehat\vu 
+ p_h(\vr) \Div \vc U  \Big] } \dt + \int_0^\tau \intO{ \partial_t p_h(r) } \dt
$$
$$
+\int_0^\tau \intO{ \Big[ \vr \partial_t\left( \frac{1}{2} |\vc U|^2 -  H_h'(r) \right)  + 
\vr \widehat\vu \cdot \Grad \left( \frac{1}{2} |\vc U|^2  -  H_h'(r) \right) \Big] } \dt 
$$
$$
+\int_0^\tau{\cal R}_h^R[\vr,\vv,r,\vc V,\mathfrak{u}_B](t){\rm d}t
+ { R}^E_h[\vr,\vu](\tau)
$$
where
$$
{\cal R}_h^R[\vr,\vv,r,\vc V,\mathfrak{u}_B]=
<\Pi^C_{h,h}, \frac 12|\vc V|^2- H_h'(r)>- <\Pi^M_{h,h},\vc V> 
$$
$$
+\int_\Omega\Big(\mathbb{S}(\nabla_h\vu):\nabla_h(\tilde{\mathfrak{u}}_B-\mathfrak{u}_B)+
p_h(\vr)\Divh(\mathfrak{u}_B-\tilde{\mathfrak{u}}_B)+\vr\partial_t\vc U\cdot (\widehat{\tilde{\mathfrak{u}}}_B-\mathfrak{u}_B)
$$
$$
+\vr\widehat\vu\cdot\nabla_h\mathfrak{u}_B\cdot(\widehat{\tilde{\mathfrak{u}}}_B-\tilde{\mathfrak{u}}_B)
+\vr\widehat\vu\cdot\Grad\vc V\cdot (\tilde{\mathfrak{u}}_B-\mathfrak{u}_B)
+\vr\widehat\vu\cdot\nabla_h(\tilde{\mathfrak{u}}_B -\mathfrak{u}_B)\cdot(\vc V-\widehat\vv)\Big){\rm d}x.
$$
\item
Regrouping conveniently several terms in the above expression, we get
\begin{equation} \label{RR5}
\left[ \intO{{\cal E}_h\left( \vr, \widehat\vv \ \Big|\ r, \vc V \right) } \right]_{t = 0}^{ t = \tau} + 
\int_0^\tau \intO{ \mathbb{S}(\nabla_h\vu):\nabla_h(\vu-\vc U)} \dt 
\end{equation}
$$ 
{ +\int_0^\tau \int_{\Gamma_{\rm out}} \left[ H_h(\vr) - H_h'(r) (\vr - r) - H_h(r)  \right]  \vu_B \cdot \vc{n} \ {\rm d} S_x \dt}
$$
$$
{ +\int_0^\tau \int_{\Gamma_{\rm in}} \left[ H_h(\vr_B) - H_h'(r) (\vr_B - r) - H_h(r)  \right]  \vu_B \cdot \vc{n} \ {\rm d} S_x \dt}
$$
$$
\aleq
- \int_0^\tau \intO{ \vr (\vc U - \widehat\vu) \cdot \Grad \vc U\cdot (\vc U - \widehat\vu)  } \dt
$$
$$
- \int_0^\tau \intO{ \Big[ p_h(\vr) - p_h'(r) (\vr - r) - p_h(r) \Big] \Div \vc U } \dt   
$$ 
$$
+ \int_0^\tau \intO{ \frac{\vr}{r} (\vc U - \widehat\vu) \cdot \Big[ \partial_t (r \vc U)  +  \Div (r \vc U \otimes \vc U) 
  + \Grad p_h(r) \Big] } \dt
	$$
	$$
	+ \int_0^\tau \intO{ \left( \frac{\vr}{r} (\widehat\vu - \vc U) \cdot \vc U +
	p_h'(r)\left( 1 - \frac{\vr}{r} \right) \right) \Big[ \partial_t r   +  \Div (r \vc U) 
  \Big] } \dt
	$$
	$$
	+\int_0^\tau{\cal R}_h^R[\vr,\vv,r,\vc V,\mathfrak{u}_B](t){\rm d}t
+ {R}^E_h[\vr,\vu](\tau),
	$$
	\item We shall continue to re-arranging the inequality (\ref{RR5}) and write
\begin{equation} \label{RR5+}
\left[ \intO{{\cal E}_h\left( \vr, \widehat\vv \ \Big|\ r, \vc V \right) } \right]_{t = 0}^{ t = \tau} + 
\int_0^\tau \intO{ \mathbb{S}(\nabla_h(\vu-\tilde{\vc U})):\nabla_h(\vu-\tilde\vc U)} \dt 
\end{equation}
$$ 
{ +\int_0^\tau \int_{\Gamma_{\rm out}} \left[ H_h(\vr) - H_h'(r) (\vr - r) - H_h(r)  \right]  \vu_B \cdot \vc{n} \ {\rm d} S_x \dt}
$$
$$
{ +\int_0^\tau \int_{\Gamma_{\rm in}} \left[ H_h(\vr_B) - H_h'(r) (\vr_B - r) - H_h(r)  \right]  \vu_B \cdot \vc{n} \ {\rm d} S_x \dt}
$$
$$
\aleq
- \int_0^\tau \intO{ \vr (\vc V - \widehat\vv) \cdot \Grad \vc U\cdot (\vc V - \widehat\vv)  } \dt
$$
$$
- \int_0^\tau \intO{ \Big[ p_h(\vr) - p_h'(r) (\vr - r) - p_h(r) \Big] \Div \vc U } \dt   
$$ 
$$
+ \int_0^\tau \intO{ \frac{\vr}{r} (\vc U - \widehat\vu) \cdot \Big[ \partial_t (r \vc U)  +  \Div (r \vc U \otimes \vc U) 
  + \Grad p(r) \Big] } \dt
	-\int_0^\tau\intO{\mathbb{S}(\Grad{\vc U}):\nabla_h(\vu-\tilde\vc U)}{\rm d}t
	$$
	$$
	+ \int_0^\tau \intO{ \left( \frac{\vr}{r} (\widehat\vu - \vc U) \cdot \vc U +
	p_h'(r)\left( 1 - \frac{\vr}{r} \right) \right) \Big[ \partial_t r   +  \Div (r \vc U) 
  \Big] } \dt
	$$
	$$
	+\int_0^\tau{\mathbb{R}}_h^R[\vr,\vv,r,\vc V,\mathfrak{u}_B](t){\rm d}t
+ {R}^E_h[\vr,\vu](\tau),
	$$
	where
	$$
	{\mathbb{R}}_h^R[\vr,\vv,r,\vc V,\mathfrak{u}_B]={\cal R}_h^R[\vr,\vv,r,\vc V,\mathfrak{u}_B] +\intO{\mathbb{S}(\nabla_h\vu):\Gradh(\vc U-\tilde{\vc U})} -\intO{\mathbb{S}(\Gradh(\tilde{\vc U}-\vc {U})):\nabla_h(\vu-\tilde\vc U)}
	$$
	$$
	+\intO{(\vr\widehat{\tilde{\mathfrak{u}}}_B-\mathfrak{u}_B)\cdot\Grad
	\vc U\cdot(\vc U-\widehat\vu)}
	+\intO{\vr(\widehat\vv-\vc V)\cdot\Grad
	\vc V\cdot(\mathfrak{u}_B-\widehat{\tilde{\mathfrak{u}}}_B)}
	+2 h^\eta\intO{ {\vr}\Grad r\cdot (\vc U - \widehat\vu)}.
	$$

	\item The strategy we want to apply recommends to estimate the right hand side of (\ref{RR5})--where $[r,\vc U]$ is a strong solution of (\ref{NS1}--\ref{inB})-- via the relative energy functional. Seeing the quantities being compared in the relative energy functional, we shall still rewrite (\ref{RR5+}) as follows:
\begin{equation} \label{RR5++}
\left[ \intO{{\cal E}_h\left( \vr, \widehat\vv \ \Big|\ r, \vc V \right) } \right]_{t = 0}^{ t = \tau} + 
\int_0^\tau \intO{ \mathbb{S}(\nabla_h(\vu-\tilde{\vc U})):\nabla_h(\vu-\tilde\vc U)} \dt 
\end{equation}
$$ 
{ +\int_0^\tau \int_{\Gamma_{\rm out}} \left[ H_h(\vr) - H_h'(r) (\vr - r) - H_h(r)  \right]  \vu_B \cdot \vc{n} \ {\rm d} S_x \dt}
$$
$$
{ +\int_0^\tau \int_{\Gamma_{\rm in}} \left[ H_h(\vr_B) - H_h'(r) (\vr_B - r) - H_h(r)  \right]  \vu_B \cdot \vc{n} \ {\rm d} S_x \dt}
$$
$$
\aleq
- \int_0^\tau \intO{ \vr (\vc V - \widehat\vv) \cdot \Grad \vc U\cdot (\vc V - \widehat\vv)  } \dt
$$
$$
- \int_0^\tau \intO{ \Big[ p_h(\vr) - p_h'(r) (\vr - r) - p_h(r) \Big] \Div \vc U } \dt   
$$ 
$$
+ \int_0^\tau \intO{1_{\vr\ge 2\overline r}(\vr) \frac{\vr}{r} (\vc U - \widehat\vu) \cdot \Big[ \partial_t (r \vc U)  +  \Div (r \vc U \otimes \vc U) 
  + \Grad p(r) \Big] } \dt
	$$
	$$
	+\int_0^\tau\intO{1_{\vr\ge 2\overline r}(\vr)(\hat\vu-\vc U)\cdot\Div(\mathbb{S}(\Grad{\vc U}))}{\rm d}t
	$$
	$$
+ \int_0^\tau \intO{1_{\vr\le 2\overline r}(\vr) \frac{\vr}{r} (\tilde\vc U - \vu) \cdot \Big[ \partial_t (r \vc U)  +  \Div (r \vc U \otimes \vc U) 
  + \Grad p(r) \Big] } \dt
	$$
	$$
	+\int_0^\tau\intO{1_{\vr\le 2\overline r}(\vr)(\vu-\tilde\vc U)\cdot\Div(\mathbb{S}(\Grad{\vc U}))}{\rm d}t
	$$
	$$
	+ \int_0^\tau \intO{ \left( \frac{\vr}{r} (\widehat\vu - \vc U) \cdot \vc U +
	p_h'(r)\left( 1 - \frac{\vr}{r} \right) \right) \Big[ \partial_t r   +  \Div (r \vc U) 
  \Big] } \dt
	$$
	$$
	+\int_0^\tau{{R}}_h^R[\vr,\vv,r,\vc V,\mathfrak{u}_B](t){\rm d}t
+ {R}^E_h[\vr,\vu](\tau),
	$$
	where
	$$
	{{R}}_h^R[\vr,\vv,r,\vc V,\mathfrak{u}_B]={\mathbb{R}}_h^R[\vr,\vv,r,\vc V,\mathfrak{u}_B] 
	$$
	$$
	+\intO{1_{\vr\ge 2\overline r}(\vr)\Div(\mathbb{S}(\Grad{\vc U}))\cdot(\vu-\hat\vu+\vc U-\tilde{\vc U})}
	$$
	$$
+ \int_0^\tau \intO{1_{\vr\le 2\overline r}(\vr) \frac{\vr}{r} (\vc U-\tilde{\vc U} + \vu-\widehat\vu) \cdot \Big[ \partial_t (r \vc U)  +  \Div (r \vc U \otimes \vc U) 
  + \Grad p(r) \Big] } \dt
	$$
	$$
	-\sum_{K\in {\cal T}}\sum_{\sigma\in {\cal E}(K)}\int_\sigma\vc n_{\sigma, K}\cdot\mathbb{S}(\Grad\vc U)\cdot(\vu-\tilde{\vc U})
	{\rm d} S_x.
	$$
	In the above, we have, among others, integrate by parts in $\intO{\mathbb{S}(\Grad{\vc U}):\nabla_h(\vu-\tilde\vc U)}$
	and then "dispatch"  the volume integral conveniently, to the sets $\{\vr\le 2\overline r\}$ and $\{\vr> 2\overline r\}$.
	We recall that $\overline r$ is defined in (\ref{C}).
	
Since $\vu-\tilde{\vc U}=\vv-\tilde{\vc V}$,  the latter term in ${{R}}_h^R[\vr,\vv,r,\vc V,\mathfrak{u}_B]$ is equal
to
$$
-\sum_{K\in {\cal T}}\sum_{\sigma\in {\cal E}(K)}\int_\sigma\vc n_{\sigma, K}\cdot\mathbb{S}(\Grad\vc U)\cdot(\vv-\tilde{\vc V})
	{\rm d} S_x	
$$
and thus its absolute value is $\aleq \sqrt h$. Indeed,	for any $\vc w\in V_0(\Omega;R^3)$
$$
\sum_{K\in {\cal T}}\sum_{\sigma\in{\cal E}(K)}\int_{\sigma}
\vc n_{\sigma,K}\otimes\Grad\vc U\otimes\vc w{\rm d}S_x
=
\sum_{K\in {\cal T}}\sum_{\sigma\in{\cal E}(K)\cap {\cal E}_{\rm int}}\int_{\sigma}
\vc n_{\sigma,K}\otimes(\Grad\vc U-(\Grad\vc U)_\sigma)\otimes(\vc w-\vc w_\sigma){\rm d}S_x
$$
$$
+ \sum_{K\in {\cal T}}\sum_{\sigma\in{\cal E}(K)\cap {\cal E}_{\rm ext}}\int_{\sigma}
\vc n_{\sigma,K}\otimes\Grad\vc U\otimes\vc w{\rm d}S_x
$$
where we have used the fact that
both  mean values of $\vc w$  and $\Grad\vc U$ are continuous over each face $\sigma\in {\cal E}_{\rm int}$. Consequently, 
$$
\Big|\sum_{K\in {\cal T}}\sum_{\sigma\in {\cal E}(K)\cap {\cal E}_{\rm int}}\int_{\sigma}
\vc n_{\sigma_K}\otimes \Grad\vc U\otimes\vc w{\rm d}S_x\Big|\aleq h \|\Grad^2\vc U\|_{L^\infty(\Omega)}
\|\nabla_h\vc w\|_{L^2(\Omega)},
$$
$$
\Big|\sum_{K\in {\cal T}}\sum_{\sigma\in{\cal E}(K)\cap {\cal E}_{\rm ext}}\int_{\sigma}
\vc n_{\sigma,K}\otimes\Grad\vc U\otimes\vc w{\rm d}S_x\Big|\aleq
\sqrt h\|\Grad \vc U\|_{L^\infty(\Omega)}
\|\nabla_h\vc w\|_{L^2(\Omega)},
$$
where we have adapted the reasoning from (\ref{dod*}) (see also (\ref{cest3}) and (\ref{cest4})). 
	
	Further, revisiting (\ref{RR4}), the form of $\mathbb{R}_h^R$ in (\ref{RR5+}), recalling that $\intO{p_h(\vr){\rm div}_h(\mathfrak{u}_B-\tilde{\mathfrak{u}}_B)}=0$, cf. (\ref{vv1}),
	using (\ref{NC1}), (\ref{NM1}), (\ref{enb3}) and employing (\ref{PoincareV-}), (\ref{errorV}), (\ref{globalV}) together
	with estimates proved in Lemma \ref{Lestimates}, we deduce
\bFormula{NR1}
|\int_0^\tau{{R}}_h^R[\vr,\vv,r,\vc V,\mathfrak{u}_B](t){\rm d}t|\aleq C h^{\alpha_{ R}},\quad
R^E_h[\vr,\vu](\tau)\aleq  h,\; \alpha_R=\min\{\alpha_C,\alpha_M,\eta\}
	\eF
with some number
		$C>0$ which depends on the strong solution $[r,\vc U]$ as indicated in (\ref{C}).
\end{enumerate}

We have shown the following result. 

\begin{Lemma}{\rm [{ Relative energy inequality for numerical solutions}]} \label{RRP1} 
Let $[\vr, \vu]$ be a numerical solution  of the algebraic system (\ref{N1}--\ref{N3}) in the setting (\ref{range1}) with the pressure 
satisfying (\ref{pres1}--\ref{pres4})
and the initial and boundary conditions verifying (\ref{ibc}--\ref{ruB}). Suppose that the couple $[r,\vc V]$ 
belongs to the class (\ref{tests}). Let $\vc U=\vc V+\mathfrak{u}_B$.  Then the relative energy inequality (\ref{RR5++}) 
holds for all $0 \leq \tau \leq T$. The remainders ${{R}}_h^R$, ${R}_h^E$ verify estimates (\ref{NR1}).
\end{Lemma}

\subsection{Error estimates}
Now, we use in the relative energy inequality (\ref{RR5++}) as test functions a strong solution $[r,\vc U]$
of problem (\ref{NS1}--\ref{inB}) in the class (\ref{tests}).

In view of (\ref{NS1}) we can rewrite (\ref{RR5++}) in the form
\begin{equation} \label{RR6}
\left[ \intO{{\cal E}_h\left( \vr, \widehat\vv \ \Big|\ r, \vc V \right) } \right]_{t = 0}^{ t = \tau} + 
\int_0^\tau \intO{ \mathbb{S}(\nabla_h(\vu-\tilde{\vc U})):\nabla_h(\vu-\tilde{\vc U})} \dt 
\end{equation}
$$ 
{ +\int_0^\tau \int_{\Gamma_{\rm out}} \left[ H_h(\vr) - H_h'(r) (\vr - r) - H_h(r)  \right]  \vu_B \cdot \vc{n} \ {\rm d} S_x \dt}
$$
$$
\aleq  
- \int_0^\tau \intO{ \vr (\vc V - \widehat\vv) \cdot \Grad \vc U\cdot (\vc V - \widehat\vv)  } \dt
$$
$$
- \int_0^\tau \intO{ \Big[ p_h(\vr) - p_h'(r) (\vr - r) - p_h(r) \Big] \Div \vc U } \dt   
$$ 
$$
+\int_0^\tau\intO{1_{\vr\le 2\overline r}(\vr)\frac{\vr-r} r(\tilde{\vc U}-\vu)\cdot{\rm div}_x\mathbb{S}(\Grad\vc U)}{\rm d}t
+\int_0^\tau\intO{1_{\vr\ge 2\overline r}(\vr)\frac{\vr-r} r(\vc U-\widehat\vu)\cdot{\rm div}_x\mathbb{S}(\Grad\vc U)}{\rm d}t + h^{\alpha_R},
$$
where we have used the inequality
$$
\Big|\int_0^\tau \int_{\Gamma_{\rm in}} \left[ H_h(\vr_B) - H_h'(r) (\vr_B - r) - H_h(r)  \right]  \vu_B \cdot \vc{n} \ {\rm d} S_x \dt\Big|\aleq h^2
$$
(recall that $\vr_B=\widehat{\mathfrak{ r}}_B$).

At this stage, it is convenient to recall a simple but in our context important
algebraic lemma:

\bLemma{relaxed2}
Let $0<a<b<\infty$ and let p satisfiy (\ref{pres1}--\ref{pres4}).
Then there exists a number $c=c(a,b)>0$ such that for all $\vr\in [0,\infty)$ and $r\in [a,b]$,
\begin{equation}\label{E1}
E(\vr|r)\ge c(a,b)\Big( 1_{{\cal O}_{\rm res}}(\vr) + \vr 1_{{\cal O}_{\rm res}}(\vr)
+1_{{\cal O}_{\rm res}}(\vr)p(\vr)+ (\vr-r)^2 1_{{\cal O}_{\rm ess}}(\vr)\Big),
\end{equation}
where $E$ is defined in (\ref{calE}) and
${\cal O}_{\rm ess}=[a/2, 2b]$, ${\cal O}_{res}= [0,\infty)\setminus {\cal O}_{\rm ess}$.
\eL
	
	With Lemma \ref{Lrelaxed2} it is rudimentary to see that
	$$
	\Big|\int_0^\tau \intO{ \vr (\vc V - \widehat\vv) \cdot \Grad \vc U\cdot (\vc V - \widehat\vv)  } \dt\Big|\aleq C \int_0^\tau{\cal E}(\vr,\widehat\vv| r,\vc V){\rm d}t
	$$
	and, under assumptions (\ref{pres1}--\ref{pres4}) on the pressure, also 
	$$
	\Big|\int_0^\tau \intO{ \Big[ p(\vr) - p'(r) (\vr - r) - p(r) \Big] \Div \vc U } \dt\Big| \aleq C \int_0^\tau{\cal E}(\vr,\widehat\vv| r,\vc V){\rm d}t,
	$$
	We observe that, under (\ref{pres1}--\ref{pres4}), also 
	$$
	\Big|\int_0^\tau\int_{\vr\ge 2\overline r}\frac{\vr-r} r(\vc U-\hat\vu)\cdot{\rm div}_x\mathbb{S}(\Grad\vc U){\rm d}x{\rm d}t\Big|
	\aleq C \Big(\int_0^\tau{\cal E}(\vr,\widehat\vv| r,\vc V){\rm d}t + h\Big)
	$$
	while
	$$
	\Big|\int_0^\tau\int_{\vr\le 2\overline r}\frac{\vr-r} r(\tilde{\vc U}-\vu)\cdot{\rm div}_x\mathbb{S}(\Grad\vc U){\rm d}x{\rm d}t\Big|\aleq 
	\delta \int_0^\tau\|\tilde{\vc U}-\vu\|_{L^2(\Omega)}^2{\rm d}t +c(\delta)C \int_0^\tau{\cal E}(\vr,\widehat\vv| r,\vc V){\rm d}t 
	$$
	In the last four estimates $\delta>0$, $c=c(\delta)>0$ and $C>0$ depends on the strong solution at most as indicated in (\ref{C})
	and $\overline r$ is defined in (\ref{C}).

	Finally recalling the inequality, cf. (\ref{SobolevV0}),
	$$
	\|\tilde{\vc U}-\vu\|_{L^2(\Omega)}^2 \aleq \|\nabla_h(\tilde{\vc U}-\vu)\|_{L^2(\Omega)}^2,
	$$
	and  estimates  (\ref{NR1}),
	we deduce from (\ref{RR5++}),
$$
\left[ \intO{{\cal E}_h\left( \vr, \widehat\vv \ \Big|\ r, \vc V \right) } \right]_{t = 0}^{ t = \tau} + 
\int_0^\tau\Big(  \|\tilde{\vc U}-\vu\|_{L^2(\Omega)}^2 + \|\nabla_h(\tilde{\vc U}-\vu)\|_{L^2(\Omega)}^2\Big){\rm d}t
$$
$$
\aleq C\Big(h^{\alpha_R}+ \int_0^\tau{\cal E}_h(\vr,\widehat\vv| r,\vc V){\rm d}t\Big).
$$

Applying to (\ref{RR6}) the Gronwall argument, we get the following lemma
\begin{Lemma}\label{Error}{\rm [Error estimates for the numerical solution]}
Let assumptions of Lemma \ref{RRP1} be satisfied. Let $[\vr_h,\vu_h]$ be a numerical solution  of the algebraic system (\ref{N1}--\ref{N3}) in the setting (\ref{range1})
emanating from the initial and boundary value data (\ref{ibc}--\ref{ruB}). 
Suppose that the couple $[r,\vc U]$, $\vc U=\vc V+\mathfrak{u}_B$ belongs to the class (\ref{tests}) and represents a strong solution of the problem (\ref{NS1}--\ref{inB}). Then
$$
\left[ \intO{{\cal E}_h\left( \vr_h, \widehat\vv_h \ \Big|\ r, \vc V \right) } \right]_{t = 0}^{ t = \tau} + 
\int_0^\tau\Big(  \|\tilde{\vc U}-\vu_h\|_{L^2(\Omega)}^2 + \|\nabla_h(\tilde{\vc U}-\vu_h)\|_{L^2(\Omega)}^2\Big){\rm d}t\aleq C h^{\alpha_R},
$$
where $\alpha_R=\min\{\alpha_C,\alpha_M, \eta\}>0$ with $C>0$ dependent on the strong solution in the way indicated in (\ref{RR5}).
\end{Lemma}

This finishes the proof of Theorem \ref{TN2}.

\bRemark{rfinal}
In view of estimates evoked in Remark \ref{Rrest}, and due to Remarks \ref{RRc-consistency}, \ref{RRm-consistency} and (\ref{enb3}),
Lemma \ref{Error} remains true with $\alpha_R=\min\{\alpha_C,\alpha_M\}>0$ also for the numerical solutions of the scheme (\ref{N1}--\ref{N3}) { with $\tilde\kappa=0$
and $p$ in the class (\ref{pres1}--\ref{pres4}) in the situations described by conditions (\ref{pres2}),(\ref{range2}) resp. (\ref{pres2}), (\ref{pres3}), (\ref{range3}), resp. 
(\ref{pres2}), (\ref{pres3}), (\ref{range4}).} We let the details for the interested reader.
\eR

\section{Concluding remarks}\label{CR}
In this section, we wish to mention some open problems related to the numerical approximations of the Navier-Stokes equations
with non-homogenous boundary data.
\begin{enumerate}
\item In this paper, we consider only the case when the computational domain $\Omega_h=\cup_{K\in {\cal T}_h}K_h$ coincides with the physical domain $\Omega$.
If this is not the case, one must work on so called unfitted meshes (see e.g. Babu\v ska \cite{Ba}), i.e. $\Omega_h\neq \Omega$. For the no-slip boundary conditions, 
a solution has been proposed  in \cite{FeHoMaNo}, \cite{FeKaMi} but for the non-homogenous boundary data, this problem is clearly more involved. Namely, the quality of approximation of $\Omega$ by $\Omega_h$ and of $\partial\Omega$
by $\partial\Omega_h$ will play a preponderant role, and remains to be  determined.
\item Karper's scheme \cite{Ka} is a collocalized scheme (the density and velocity are discretized on the same "primal" mesh). Quite often, the
staggered schemes (the density resp. the pressure are discretized on a primal mesh while the velocity components on a "bi-dual" meshes) --as e.g. the Marker and Cell (MAC) finite difference schemes, \cite{HaAm}, \cite{GaHeLaMa}, \cite{GHMN-MAC}, \cite{FeGa}, \cite{GaGaLaHe} or staggered discretizations combined with Rannacher-Turek finite elements or Crouzeix-Raviart finite elements \cite{HeLaMiTh}-- are computationally more efficient. It would be certainly of interest to
extend the convergence results to dissipative solutions to these types of schemes. 
\item Since the proof of convergence to dissipative solutions is relatively weakly  scheme-structure dependent it would be of great interest to determine
as weak as possible universal  criteria imposed on the numerical scheme in order to obtain the convergence to the dissipative solutions. The universal
description of several staggered schemes in one unique formalism provided in \cite{HeLaMiTh} may be a starting point to approach this task.
\item Gallouet et al. established in \cite{GaHeMaNo} the error estimates for the Karper's scheme and in \cite{GHMN-MAC} for the MAC scheme 
under assumption $\gamma>3/2$ without any correlation between the time step $\Delta t$ and the size of the space discretization $h$ by using a different approach,
which consists in mimicking the proof of the weak-strong uniqueness known from  the continuous case (see \cite{FeJiNo}). For "large" values of $\gamma$, these estimates seem to be optimal judging from what can be obtained for the sole continuity equation with the same regularity of the transporting velocity field, see
\cite{MaNo} for the discussion about the optimality. It is certainly of interest to try to prove similar results with the general non homogenous boundary
data. The weak strong uniqueness principle established in this boundary value setting in the continuous case in \cite{KwNo} and later in \cite{AbFeNo} encourages 
the attempts in this direction.
\item Weak solutions to the problem (\ref{NS1}--\ref{inB}) are known to exist for the adiabatic coefficients $\gamma>3/2$, see \cite{FNP} for the
no-slip case and \cite{ChJiNo}, \cite{ChoNoY}, \cite{KwNo} for the general inflow-outflow boundary conditions. The proofs use, among others, very specific 
mathematical tools closely related to the structure of the system -- as. e.g. compensated compactness and various commutator lemmas-- the tools whose numerical counterpart is usually not available. The consistency formulation of the balance laws allows to bring the numerical solution close to a weak solution modulo
a remainder, which in the case of the convergence to weak solutions must have a "better quality" than in the proofs of the convergence to the dissipative solutions.
T. Karper provided in 2013 in \cite{Ka} a convergence proof of the Karper's scheme to a weak solution under assumption $\gamma>3$ in the case of the no-slip boundary conditions. To prove the same for the general inflow-outflow boundary data is definitely an imminent open problem whose solutions is of independent interest.
\item Last but not least some more complex hydrodynamical models with the similar
structure of convective terms as e.g. models of compressible fluids with non-linear stress, see \cite{AbFeNo}, fluids of compressible polymers, see S\"uli et al. \cite{BaSu}, \cite{FeLuSu} or multi-fluid models with differential closure, see \cite{NoSCM}, could be treated
on the basis of the methodology introduced in this paper. These studies would  certainly be of
a non negligible interest.

\end{enumerate}

\end{document}